\documentclass[nonblindrev]{informs3_noinforms}

\usepackage{enumerate}
\newtheorem{observation}{Observation}

\newcommand{\pr}{\mathbb{P}}
\newcommand{\E}{\mathbb{E}}

\newcommand{\ignore}[1]{\relax}

\renewcommand{\theequation}{\thesection.\arabic{equation}}

\usepackage[square,numbers]{natbib}
\usepackage[colorlinks=true,breaklinks=true,bookmarks=true,urlcolor=blue,
     citecolor=blue,linkcolor=blue,bookmarksopen=false,draft=false]{hyperref}

\def\EMAIL#1{\href{mailto:#1}{#1}}
\def\URL#1{\href{#1}{#1}}         

\TheoremsNumberedThrough     

\EquationsNumberedThrough    

\begin{document}



\RUNTITLE{Simple and explicit bounds for multi-server queues with $\frac{1}{1-\rho}$ scaling}

\TITLE{Simple and explicit bounds for multi-server queues with $\frac{1}{1-\rho}$ scaling}

\ARTICLEAUTHORS{%
\AUTHOR{\bf Yuan Li}
\AFF{Amazon}
\AUTHOR{\bf David A. Goldberg}
\AFF{Cornell ORIE, \EMAIL{dag369@cornell.edu}, \URL{https://people.orie.cornell.edu/dag369/}}
}
\RUNAUTHOR{{Yuan Li and David A. Goldberg}}


\ABSTRACT{%
We consider the FCFS $GI/GI/n$ queue, and prove the first simple and explicit bounds that scale as $\frac{1}{1-\rho}$ under only the assumption that inter-arrival times have finite second moment, and service times have finite $2+\epsilon$ moment for some $\epsilon > 0$.  Here $\rho$ denotes the corresponding traffic intensity.  Conceptually, our results can be viewed as a multi-server analogue of Kingman's bound.  Our main results are bounds for the tail of the steady-state queue length and the steady-state probability of delay.  The strength of our bounds (e.g. in the form of tail decay rate) is a function of how many moments of the service distribution are assumed finite.  Our bounds scale gracefully even when the number of servers grows large and the traffic intensity converges to unity simultaneously, as in the Halfin-Whitt scaling regime.  Some of our bounds scale better than $\frac{1}{1-\rho}$ in certain asymptotic regimes.  In these same asymptotic regimes we also prove bounds for the tail of the steady-state number in service.
\\\\Our main proofs proceed by explicitly analyzing the bounding process which arises in the stochastic comparison bounds of \citet{GG13} for multi-server queues.  Along the way we derive several novel results for suprema of random walks and pooled renewal processes which may be of independent interest.  We also prove several additional bounds using drift arguments (which have much smaller pre-factors), and point out a conjecture which would imply further related bounds and generalizations.  We also show that when all moments of the service distribution are finite and satisfy a mild growth rate assumption, our bounds can be strengthened to yield explicit tail estimates decaying as $O\big(\exp(-x^{\alpha})\big)$, with $\alpha \in (0,1)$ depending on the growth rate of these moments.
}



\maketitle

%


\section{Introduction.}\label{Introduction}
The multi-server queue with independent and identically distributed (i.i.d.) inter-arrival and service times, and first-come-first-serve (FCFS) service discipline, is a fundamental object of study in Operations Research and Applied Probability.  Its study was originally motivated by the design of telecommunication networks in the early 20th century (\citet{Janssen08}).  Since that time the model has found many additional applications across a wide range of domains (\citet{Worth09}).  
\\\indent A key result for $GI/GI/1$ (i.e. single server) queues came in the seminal 1962 paper of John Kingman (\citet{Kingman62a}), in which a simple and explicit upper bound was given for the steady-state expected waiting time $\E[W]$ (and thus also for $\E[L]$ by Little's Law, see e.g. \citet{Bertsimas95}).  This bound, now referred to as \textbf{Kingman's bound}, states the following (with $A$ a random variable with the same distribution as an inter-arrival time, and $S$ a random variable with the same distribution as a service time): $\E[W] \leq \frac{\sigma^2_A + \sigma^2_S}{2 \E[A]} \times \frac{1}{1 - \rho}, \E[L] \leq \frac{\sigma^2_A + \sigma^2_S}{2 (\E[A])^2} \times \frac{1}{1 - \rho}$.  Here $\sigma_A$ (respectively $\sigma_S$) denotes the standard deviation of the inter-arrival (respectively service) distribution, $\E[A]$ denotes the mean inter-arrival time, and $\rho$ denotes the traffic intensity.  For $GI/GI/1$ queues, $\rho = \frac{\E[S]}{\E[A]}$.  We now show that one can rewrite Kingman's bound for $\E[L]$ in terms of only $\rho$ and certain quantities derived from the distributions of $A$ and $S$ which are invariant under scaling $A$ or $S$ by a constant.  In particular, for a general random variable (r.v.) $X$ with finite mean and variance, let $\mu_X \stackrel{\Delta}{=} \frac{1}{\E[X]}$, and $c_X$ denote the coefficient of variation (i.e. $c_X \stackrel{\Delta}{=} \mu_X \sigma_X$).  Note that for a given r.v. X, $c_X$ is invariant under scaling X by a constant.  Then it follows from some straightforward algebra that Kingman's bound for $\E[L]$ may be equivalently written as follows: $\E[L] \leq \frac{1}{2} \big( c^2_A + \rho^2 c^2_S \big) \times \frac{1}{1 - \rho}$.  Importantly, note that since $c^2_A$ and $c^2_S$ are invariant under scaling $A$ and $S$ by any constants, the above bound is in some sense insensitive / unchanged / scale-free as one passes to heavy-traffic.  The term $\frac{1}{1 - \rho}$ dictates how $\E[L]$ scales as $\rho \uparrow 1$ in a broad sense.  Note that the following slight weakening of Kingman's bound (which is essentially equivalent as $\rho \uparrow 1$),
\begin{equation}\label{king12}
\E[L] \leq \frac{1}{2} \big( c^2_A + c^2_S \big) \times \frac{1}{1 - \rho},
\end{equation}
captures the fundamental essence of Kingman's bound in heavy traffic.  As most performance metrics of the general $GI/GI/1$ queue have no simple closed-form solution, this combination of simplicity, accuracy, and scalability has made Kingman's bound very popular in the analysis of queueing systems.   In the same paper, Kingman established that (under appropriate technical conditions) this bound becomes tight as $\rho \uparrow 1,$ and the bound was later tightened in \citet{Whitt.82, Whitt.82b, Daley.77}.
\\\indent Soon after Kingman's seminal work, other authors had begun using the tools of weak convergence to attempt to extend this analysis to the more complicated multi-server queue.  Indeed, Kingman himself conjectured such a result in 1965 for a sequence of $GI/GI/n$ queues which approaches heavy-traffic (with $n$ held fixed) as the traffic intensity $\rho$ (defined as $\frac{\mu_A}{n \mu_S}$ for $GI/GI/n$ queues) approaches 1 (\citet{K65}).  Such a weak convergence result was proven in \citet{Kollerstrom74}.  Namely, it was proven in \citet{Kollerstrom74} that if one considers an appropriate sequence of $GI/GI/n$ queues in heavy-traffic (indexed by their traffic intensity $\rho$, with $W_{\rho}$ and $L_{\rho}$ the corresponding steady-state waiting time and queue length), with $n$ held fixed as $\rho \uparrow 1$, then (letting $\Rightarrow$ denote weak convergence and $E_1$ denote a mean one exponentially distributed r.v.)
\begin{equation}\label{king1weakmulti}
\lbrace (1 - \rho) W_{\rho} , \rho \uparrow 1 \rbrace \Rightarrow \frac{\E[A]}{2} \times \big( c^2_A + c^2_S \big) \times E_1\ \  ;\ \ \lbrace (1 - \rho) L_{\rho} , \rho \uparrow 1 \rbrace \Rightarrow \frac{1}{2} \big( c^2_A + c^2_S \big) \times E_1.
\end{equation}
Related results are also proven in \citet{Borovkov65,IW.70b,Loulou73,Kollerstrom79,Kennedy72,Nagaev.70,Nagaev.70b}.  \textbf{However, all these results and bounds were premised on keeping the number of servers fixed while} $\rho \uparrow 1.$  Other approaches (based on techniques including stochastic comparison) ran into similar obstacles (see e.g. \citet{K70, Wolff87, Mori.75, Scheller03,Scheller06,RSW13}).  
\\\indent That said, there has been some partial progress towards a non-asymptotic multi-server analogue of Kingman's bound (with all prefactors independent of the number of servers), with early progress for restricted asymptotic regimes and distribution classes appearing in \citet{Makino.69,Brumelle73,Suzuki70,Arjas78,Oliver.74,Seshadri.96}.  Later progress was made using a range of techniques, including convexity (e.g. \citet{Rol.76,W80,Whitt.83,Daley.84}), other modified service disciplines (e.g. \citet{Chawla17,Smith81,HB.09}), large deviations theory (e.g. \citet{Sadowsky91}), and robust optimization (e.g. \citet{Whitt.84,Whitt.84b,Whitt.84c,Bandi15}).  \textbf{However, to date, none of those approaches have been able to make substantial progress towards proving a non-asymptotic bound such as} (\ref{king12}) \textbf{for general multi-server queues.}    
\\\indent Other recent works have focused on proving bounds on the error of heavy-traffic approximations in the Halfin-Whitt asymptotic scaling regime (\citet{DDG14,BD15,BD16,BD20,BDF15,Gur13,Gur14,MMR.98,Janssen08b,Janssen11,Jin21}).  However, outside the case of Markovian service times, all of these results suffer from the presence of non-explicit constants, which may depend on the underlying service distribution in a complicated and unspecified way. Furthermore, in the Halfin-Whitt setting, the relevant limiting quantities themselves generally have no explicit representation (\citet{R.09,AR16,GM08}).  Also, at least regarding multi-server queues, essentially all such results are restricted to a particular asymptotic regime (although more universal results are known for single-server systems, see \citet{huang16,BDF15,BD16,BD20,Gaunt20}).
\\\indent Lyapunov function and drift arguments have also been used to yield bounds (\citet{GZ06,GM08,GS12,DDG14,hokstad85,grosof2018srpt, scully2020gittins, grosof2021finite}).  However, these works generally have additive error terms that scale with the number of servers, making them not amenable to proving bounds which hold in heavy-traffic regimes where $\rho$ and $n$ vary together (as in the Halfin-Whitt regime), and/ or involve non-explicit constants or restrict to certain asymptotic regimes.  An exception is the very interesting recent work \citet{wang2021zero}, which uses Lyapunov function arguments to provide a general bound for the steady-state probability of delay (s.s.p.d.), i.e. the probability that all servers are busy, in certain multi-server systems with hyper-exponentially distributed service times and Markovian inter-arrival times.  Although that work considers a more general type of queueing model in which jobs can occupy more than one server, and is thus incomparable to our own, when restricted to our setting of FCFS GI/GI/n queues, it implies the following bound.  For systems in which service times are hyper-exponentially distributed (as a mixture of $K$ exponentials) and inter-arrival times are exponentially distributed, the s.s.p.d. is at most $\frac{1}{1 - \rho}(\frac{3K}{\sqrt{n}} + \frac{1}{n})$.  Related techniques were used in \citet{hong2021sharp} to provide asymptotic upper and lower bounds for $\E[L]$ scaling like Kingman's bound for certain sequences of queues (up to some non-explicit constants) in a super-Halfin-Whitt scaling regime, in which $(1 - \rho) \sqrt{n} \log(n) \rightarrow 0.$  The authors also show that for such sequences of multi-server queues, but in a sub-Halfin-Whitt regime in which $(1 - \rho) \frac{\sqrt{n}}{\log(n)} \rightarrow \infty,$ the s.s.p.d. is asymptotically at most $\exp\big(- C n (1 - \rho)^2\big)$ for some non-explicit $C > 0.$  
\\\indent For further discussion regarding the massive literature on this problem, we refer the interested reader to the surveys \citet{Doig57,Ovuworie80,Whitt.93}.  We especially point the reader to \citet{Whitt.93} for an overview of the many heuristic approximations available in the literature.  We note that this body of work is further complicated by several mistakes in the literature, as discussed in \citet{Daley.77,Daley.84,Daley.97,Wolff87}.  Indeed, even the popular textbook \citet{Gross.11} seems to incorrectly claim the bound (\ref{king12}) for general multi-server queues in its Section 7.1.3.  

\subsubsection{Why $\frac{1}{1 - \rho}$?}
In the above discussion, we several times referenced the fact that certain bounds ``scaled as $\frac{1}{1 - \rho}$", or ``did not scale correctly" because they did not scale as $\frac{1}{1 - \rho}$.  It is of course reasonable to ask why, and in what precise sense, $\frac{1}{1 - \rho}$ should be the bar.  There are at least two fundamental justifications here.  First, it follows from well-known results for the M/M/n queue (\citet{HW.81}) that for any fixed $B > 0$, there exist $\zeta_B \in (0,1)$, independent of $n$ and $\rho$, such that any $M/M/n$ queue for which $\rho \in (1 - B  n^{-\frac{1}{2}}, 1)$ satisfies $\zeta_B \times \frac{1}{1 - \rho} \leq \E[L] \leq \frac{1}{1 - \rho}$.  Thus, in a fairly general sense, this is the correct scaling for the $M/M/n$ queue in heavy-traffic.  Second, this is the correct scaling in classical heavy-traffic, as well as the Halfin-Whitt and non-degenerate slowdown scaling regimes (\citet{Kollerstrom74,GG13,AR16,Atar11}), and in single-server queues.  As referenced above, this is also the correct scaling for certain sequences of queues in a super-Halfin-Whitt regime (\citet{hong2021sharp}).  It has also been known for some time that this is the correct scaling when service times are deterministic, as in that setting certain lower and upper bounds proved earlier by Kingman coincide asymptotically (\citet{K70}).  More broadly, several known lower bounds exhibit such a scaling across multiple heavy-traffic regimes (e.g. \citet{K70, GG13, G16, hong2021sharp}).  Indeed, the $\frac{1}{1 - \rho}$ scaling is a guiding meta-principle throughout much of the literature on multi-server queues.  

\subsubsection{Summary of state-of-the-art.} In summary, the question of whether a multi-server analogue of Kingman's bound exists remains an open problem despite over 50 years of research.  This includes not just whether the exact bound (\ref{king12}) holds, but whether any bound representable as a simple function of a few normalized moments (of $A$ and $S$) multiplied by $\frac{1}{1-\rho}$ holds.  Over the years, Daley has several times lamented on this state of affairs, and we refer the reader to \citet{Daley.77,Daley.84,Daley.97,Allen.14} for some directly related discussion.  There Daley conjectures that such a bound should indeed be possible (and even that (\ref{king12}) should hold).  However, in spite of Daley's optimism, other works bring this into question.  For example, the results of \citet{Gupta10} prove that two queues whose inter-arrival and service times have the same first two moments can still have very different mean waiting times.  Similarly, the known results for queues in the Halfin-Whitt regime (e.g. \citet{AR16,GM08,DDG14}) suggest that simple bounds may not hold in that regime.  Indeed, those works show that the limiting behavior of the steady-state waiting time may depend in a very complex way on the underlying service distribution.  That stands in contrast to the classical heavy-traffic setting in which only one simple limiting behavior is possible (as dictated by (\ref{king1weakmulti})).  Moreover, when service times only have few finite moments, deriving uniform tail bounds is very subtle even in the single-server setting (\citet{Olvera.11,Olvera.11b}).  In light of such results, it is unclear whether a simple bound which scales as $\frac{1}{1-\rho}$ across different notions of heavy-traffic and depends only on a few normalized moments even exists.  

\subsection{Our contribution.}
In this paper, we use stochastic comparison arguments (combined with several novel bounds for associated random walks) to prove the first such Kingman-like bound for general multi-server queues, only requiring the assumption that $\E[A^2] < \infty$ and $\E[S^{2 + \epsilon}] < \infty$ for some $\epsilon > 0$.  Our bounds for the steady-state queue length and probability of delay are simple, explicit, and scale as a simple function of a few normalized moments (of the inter-arrival and service distributions) multiplied by $\frac{1}{1-\rho}$, regardless of the particular notion of heavy-traffic considered (and including both the classical and Halfin-Whitt scalings).  Some of our bounds scale better than $\frac{1}{1-\rho},$ and in these same asymptotic regimes we also prove bounds for the the tail of the steady-state number in service.  We also prove several additional bounds using drift arguments (which have much smaller pre-factors).  In addition, we prove much stronger tail bounds under the assumption that all moments of the service distribution are finite and satisfy a mild growth rate assumption.

\subsection{Outline of paper.}
The remainder of our paper proceeds as follows.  In Section\ \ref{Mainresults}, we state our main results and provide a few illustrative implications and extensions.  More precisely, we state our most central results Theorem\ \ref{mastertheoremsmallr} and \ref{mastertheoremlarger} (bounds for the tail of the queue length and steady-state probability of delay under minimal assumptions) in Subsection\ \ref{mainmainsec}.  Then, we state a handful of illustrative implications and extensions in Subsection\ \ref{additionalresultssec}, and provide a separate outline of those results in Subsubsection\ \ref{suboutlinesecc}.  Let us point out that this includes our stronger results (under mild growth assumptions on the moments of $S$) in Subsubsection\ \ref{strongersec}.  We also provide a discussion of the prefactors arising in our results, and some of the limitations of our results, in Subsection\ \ref{additionaldiscusssec}.  We state our results derived using simplified drift arguments (with no large prefactors), as well as an intuitive conjecture which would imply even stronger such results, in Section\ \ref{noprefacsec}.    
\\\\Section\ \ref{mainproofsec} is devoted to the proofs of our most central results Theorems\ \ref{mastertheoremsmallr} and \ref{mastertheoremlarger}, which constitute the majority of the technical analysis of the manuscript.  As the proofs are somewhat involved, we proceed as follows to improve readability.  First, we sketch a high-level outline of the proof in Subsection\ \ref{highleveloutlinesec}.  Second, we provide a more detailed proof (but still without most technical details), containing all of the most important auxiliary results and main flow of logic (albeit in many cases without their proofs), in Subsection\ \ref{abitmoredetailsec}.  Third, we provide many of the most important technical details of the proofs (but still with many of the finer subarguments omitted) in the technical appendix Section\ \ref{moredetailsec}.  Finally, we defer many of the finer subarguments of these proofs to the supplemental appendix Section\ \ref{appsec}.  
\\\\We prove our bounds for the steady-state probability of delay and number of busy servers in Section\ \ref{sspdsec}, and provide some concluding remarks and directions for future research in Section\ \ref{concsec}.  Let us also point out that in addition to providing many of the finer subarguments of the proofs of our main results Theorems\ \ref{mastertheoremsmallr} and \ref{mastertheoremlarger}, and the proofs of our results with no large prefactors based on simple drift arguments, our supplamental appendix also includes : implications of our main results for higher order moments in Subsection\ \ref{highermomentsec}; implications of our main results for queues in the Halfin-Whitt regime (and an open question of \citet{Chawla17}) in Subsection\ \ref{HalfinChawlasec}; a more in-depth discussion of the prefactors arising in our main results in Subsection\ \ref{rnotsmallsec}; and a sketch of a plausible approach to generalizing our main results to the network setting in Subsection\ \ref{appnetworksec}.

\section{Main Results.}\label{Mainresults}
\subsection{Notation.}
Let us fix an arbitrary FCFS $GI/GI/n$ queue with inter-arrival times having the same distribution as r.v. $A,$ and service times having the same distribution as r.v. $S$, and denote this queueing system by ${\mathcal Q}^n.$  Let ${\mathcal N}_o$ (respectively ${\mathcal A}_o)$ denote an ordinary renewal process with renewals distributed as S (respectively A), and $N_o(t) \big($respectively $A_o(t) \big)$ the number of renewals in $[0,t]$.  In general, we will use script font (e.g. ${\mathcal N}_o, {\mathcal A}_o$) to refer to the corresponding stochastic process, with notation such as $N_o(t), A_o(t)$ referring to the associated counting process evaluated at a particular time.  Let $\lbrace {\mathcal N}_{e,i}, i \geq 1 \rbrace \bigg($respectively $\lbrace {\mathcal N}_{o,i}, i \geq 1 \rbrace \bigg)$ denote a mutually independent collection of equilibrium (respectively ordinary) renewal processes with renewals distributed as $S$; ${\mathcal A}_e$ an independent equilibrium renewal process with renewals distributed as $A$; and $\lbrace N_{e,i}(t), i \geq 1 \rbrace \bigg($respectively $\lbrace N_{o,i}(t), i \geq 1 \rbrace \bigg)$ and $A_e(t)$ the respective number of renewals in $[0,t]$.  We also let ${\mathcal N}_e$ be a separate independent equilibrium renewal process with renewals distributed as $S$, with $N_e(t)$ the corresponding number of renewals in $[0,t]$.
\\\indent Here we recall that an equilibrium renewal process is one in which the first renewal interval is distributed as the equilibrium distribution of $S$.  For a r.v. $X$, recall that a r.v. $R$ is distributed according to the equilibrium distribution of $X$ if $\pr( R > y) = \frac{1}{\E[X]} \int_y^{\infty} \pr( X > z ) dz$ for all $y > 0$.   For a r.v. $X$, we let $R(X)$ denote a r.v. distributed as the equilibrium distribution of $X$.  We note that such a process captures the long-run behavior of a renewal process, since under quite general assumptions the ``time until next renewal" in a renewal process with renewals distributed as $X$ converges in distribution (as time grows large) to $R(X)$.  Noting that in heavy traffic any given server of a multi-server queue behaves like a renewal process for long stretches of time (as there is some job waiting in queue to replace any job that completes service), it is intuitive that (at least in heavy-traffic) the residual service time of a busy server would have the same distribution (at least approximately).  Interestingly, it can be shown that this is true generally for $GI/GI/n$ queues, i.e. under mild technical conditions the steady-state residual work on a busy server has the equilibrium distribution of a service time (see e.g. \citet{hokstad85}).  It is also well-known that the same phenomena manifests in infinite-server models and loss models with Markovian arrival processes (\citet{Sev57,Eick93}).
\\\indent Let $\lbrace A_i, i \geq 1 \rbrace$ (respectively $\lbrace S_i, i \geq 1 \rbrace$) denote the sequence of inter-event times in ${\mathcal A}_o$ (respectively ${\mathcal N}_o)$.  Let us evaluate all empty summations to zero, and all empty products to unity; and as a convention take $\frac{1}{\infty} = 0$ and $\frac{1}{0} = \infty$.  For an event ${\mathcal E}$, let $I({\mathcal E})$ denote the corresponding indicator function.  Unless stated otherwise, all processes should be assumed right-continuous with left limits (r.c.l.l.), as is standard in the literature.  For our results involving steady-state queue lengths, we will generally require that the total number of jobs in ${\mathcal Q}^n$ (number in service + number waiting in queue) converges in distribution (as time goes to infinity, independent of the particular initial condition) to a steady-state r.v. $Q^n(\infty)$.  As a shorthand, we will denote this assumption by saying ``$Q^n(\infty)$ exists", and refer the interested reader to \citet{AS08} for a discussion of technical conditions ensuring this property holds.  We will adopt a parallel convention when talking about the steady-state waiting time of an arriving job.  Namely, we will generally require that the distribution of the waiting time (in queue, not counting time in service) of the $k$th arrival to the system converges in distribution (as $k \rightarrow \infty$, independent of the particular initial condition) to a steady-state r.v. $W^n(\infty)$.  As a shorthand, we will denote this assumption by saying ``$W^n(\infty)$ exists".  Supposing that $Q^n(\infty)$ exists, let $L^n(\infty)$ denote a r.v. distributed as the steady-state number of jobs waiting in queue, i.e. $L^n(\infty)$ is distributed as $\max\big(0, Q^n(\infty) - n \big)$.  
\\\indent In addition, for some of our results (i.e. those based on simple drift arguments) we will require that the (sorted) vector representing the residual service times of the set of jobs currently in service converges in distribution (independent of the particular initial condition) to a steady-state random vector $\overline{W}^n_{\textrm{service}}(\infty)$, and denote this by saying ``$\overline{W}^n_{\textrm{service}}(\infty)$ exists".  In such a setting, we let $\textrm{Num}^n_{\textrm{service}}(\infty)$ denote the corresponding steady-state number in service $\big($i.e. number of non-zero components of $\overline{W}^n_{\textrm{service}}(\infty)\big)$, and $\textrm{Work}^n_{\textrm{service}}(\infty)$ denote the corresponding steady-state amount of work in service (i.e. sum of components of $\overline{W}^n_{\textrm{service}}(\infty)$).  In these settings we will also require that the total amount of work in system (remaining work of those in service + service times of those in queue) converges to a steady-state r.v. $\textrm{Work}^n(\infty)$ (again independent of initial conditions), and denote this by saying ``$\textrm{Work}^n(\infty)$ exists".  
\\\indent For a general r.v. $X$, let $X^+$ denote $\max(0,X)$.  For $k \geq 1$, let $\rho_k \stackrel{\Delta}{=} \frac{\mu_A}{k \mu_S}$.  Whenever there is no ambiguity as regards a particular $GI/GI/n$ system, we will let $L(\infty),W(\infty),Q(\infty),\overline{W}_{\textrm{service}}(\infty),\textrm{Work}_{\textrm{service}}(\infty),\textrm{Num}_{\textrm{service}}(\infty),\textrm{Work}(\infty),\rho$ denote \\$L^n(\infty), W^n(\infty), Q^n(\infty),\overline{W}^n_{\textrm{service}}(\infty),\textrm{Work}^n_{\textrm{service}}(\infty),\textrm{Num}^n_{\textrm{service}}(\infty),\textrm{Work}^n(\infty),\rho_n$.  Note that for any $GI/GI/n$ queue, one can always rescale both the service and inter-arrival times so that $\E[S] = \mu_S = 1$, without changing either $\rho$ or the distribution of $Q^n(\infty)$.  As doing so will simplify (notationally) several arguments and statements, sometimes we impose the additional assumption that $\E[S] = \mu_S = 1$, and will point out whenever this is the case.  For $x > 0$, we let $\Gamma(x) \stackrel{\Delta}{=} \int_0^{\infty} t^{x-1} e^{-t} dt$ denote the standard gamma function (see e.g. \citet{Batir17}).  Finally, we will sometimes use the standard Bachman-Landau (i.e. ``big-O") asymptotic notation to describe the growth rate of functions, often to informally build intuition for more formal and explicit statements and bounds.  Let us recall that two functions $f,g$ of the same parameter (say r) are said to satisfy the asymptotic relation $f = O(g)$ if there exists an absolute finite constant $C > 0$ s.t. $f(r) \leq C \times g(r)$ for all $r$ (over some appropriate unbounded domain).  Similarly, the relation $f = \Omega(g)$ indicates that there exists an absolute finite constant $c > 0$ s.t. $f(r) \geq c \times g(r)$ for all $r$ (again over an appropriate domain).  Also, the relation $f = \Theta(g)$ indicates that both $f = O(g)$ and $f = \Omega(g)$.  We note that these notations can sometimese be composed with other functions.  Thus for example the statement $f = r^{O(r)}$ would indicate that there exists $C > 0$ s.t. $f(r) \leq r^{C \times r}$ for all $r$ in some appropriate domain, while $f = r^{\Omega(r)}$ would indicate that there exists $c > 0$ s.t. $f(r) \geq r^{c \times r}$ for all $r$ in that domain.
\subsection{Main results.}\label{mainmainsec}
Our main results are the following novel, explicit, and general tail bounds for multi-server queues, which scale as $\frac{1}{1 - \rho}$, along with corresponding bounds for the steady-state probability of delay.  Our bounds only require that $\E[A^2] < \infty$ and $\E[S^{2+\epsilon}] < \infty$ for some $\epsilon > 0$, although in general the more moments of $S$ assumed finite the stronger the bounds become.  Let $r^* \stackrel{\Delta}{=} \sup\lbrace r : \E[S^r] < \infty \rbrace$, where we note that $r^*$ may equal $\infty$.  As our bounds scale quite differently as $r^* \downarrow 2$ and $r^* \uparrow \infty$, we state our results by breaking into two cases : $r^* \leq 2.5$ and $r^* > 2.5$. 

\begin{theorem}[Tail bounds when $r^* \leq 2.5$, i.e. $S$ has few finite moments]\label{mastertheoremsmallr}
Suppose that for a $GI/GI/n$ queue with inter-arrival times having the same distribution as r.v. $A,$ and service times having the same distribution as r.v. $S$, the following is true : (1)\ $\E[A^2] < \infty$; (2)\ $r^* \in (2,2.5]$; (3)\ $\mu_A < n \mu_S$; (4)\ $Q(\infty)$ exists.  Then for all $x > 0$, $\pr\big(  L(\infty) \geq \frac{x}{1 - \rho} \big)$ is at most
\begin{eqnarray*}
\ &\ &\ \inf_{r \in (2,r^*)} \Bigg( 3 \times 10^{19} \times \E[(S \mu_S)^2] \times \bigg( \big(\E[(S \mu_S)^2]\big)^{r-1}  + \E[(S \mu_S)^r] \bigg) \times (\frac{r}{2}-1)^{-(r+1)} \times x^{-\frac{r}{2}} \Bigg) \nonumber
\\&\ \ \ &\ \ +\ \ 1.1 \times \exp\bigg( - .0225 \big( \E[ (A \mu_A)^2 ] \big)^{-1} x \bigg);
\end{eqnarray*}
and the steady-state probability of delay (s.s.p.d.), $\pr\big( Q(\infty) \geq n \big)$, is at most
\begin{eqnarray*}
\ &\ &\ \inf_{r \in (2,r^*)} \Bigg( 4 \times 10^{20} \times \E[(S \mu_S)^2] \times \bigg( \big( \E[(S \mu_S)^2] \big)^{r-1} + \E[(S \mu_S)^r] \bigg) \times (\frac{r}{2}-1)^{-(r+1)} \times \big(n (1 - \rho^2) \big)^{-\frac{r}{2}} \Bigg) \nonumber
\\&\ \ \ &\ \ +\ \ 1.1 \times \exp\bigg( - .0028 \big( \E[ (A \mu_A)^2 ] \big)^{-1} n (1 - \rho)^2 \bigg).
\end{eqnarray*}
\end{theorem}

\begin{theorem}[Tail bounds when $r^* > 2.5$, i.e. $S$ has more finite moments]\label{mastertheoremlarger}
Suppose that for a $GI/GI/n$ queue with inter-arrival times having the same distribution as r.v. $A,$ and service times having the same distribution as r.v. $S$, the following is true : (1)\ $\E[A^2] < \infty$; (2)\ $r^* > 2.5$ (where $r^*$ may equal $\infty$); (3)\ $\mu_A < n \mu_S$; (4)\ $Q(\infty)$ exists.  Then for all $x > 0$, $\pr\big(  L(\infty) \geq \frac{x}{1 - \rho} \big)$ is at most
\begin{eqnarray*}
\ &\ &\ \inf_{r \in [2.5,r^*)} \Bigg( 2 \times 10^4 \times \E[(S \mu_S)^2] \times \big(10^6\big)^r \times \bigg( \big(\E[(S \mu_S)^2]\big)^{r-1} \times r^{2.5 r} + r^{1.5 r} \times \E[(S \mu_S)^r] \bigg) \times x^{-\frac{r}{2}} \Bigg)
\\&\ \ \ &\ \ +\ \ 1.1 \times \exp\bigg( - .0225 \big( \E[ (A \mu_A)^2 ] \big)^{-1} x \bigg);
\end{eqnarray*}
and the steady-state probability of delay (s.s.p.d.), $\pr\big( Q(\infty) \geq n \big)$, is at most
\begin{eqnarray*}
&&\inf_{r \in [2.5,r^*)} \Bigg( 2 \times 10^4 \times \E[(S \mu_S)^2] \times \big(10^7\big)^r \times \bigg( (\E[(S \mu_S)^2])^{r-1} \times r^{2.5 r} + r^{1.5 r} \times \E[(S \mu_S)^r] \bigg) \times \big( n(1 - \rho)^2 \big)^{-\frac{r}{2}} \Bigg)
\\&\ \ \ &\ \ +\ \ 1.1 \times \exp\bigg( - .0028 \big( \E[ (A \mu_A)^2 ] \big)^{-1} n (1-\rho)^2 \bigg).
\end{eqnarray*}
\end{theorem}
\ \\The proof of our bounds for $L(\infty)$ in Theorems\ \ref{mastertheoremsmallr} - \ref{mastertheoremlarger} constitute the largest part of our technical analysis.  The proofs are first sketched at a high-level in Section\ \ref{highleveloutlinesec}, then in greater depth in Section\ \ref{abitmoredetailsec}, with additional details appearing in the technical appendix Section\ \ref{moredetailsec} and supplemental appendix Section\ \ref{appsec}.  The proof of our bounds for the s.s.p.d. in Theorems\ \ref{mastertheoremsmallr} - \ref{mastertheoremlarger} appears in Section\ \ref{sspdsec}.
\\\\Some additional comments are in order.
\begin{itemize}
\item The $\inf$ terms appearing in Theorems\ \ref{mastertheoremsmallr} and \ref{mastertheoremlarger} can be replaced by evaluating the associated expression at any $r$ in the given range, yielding an explicit bound with decay rate $x^{-\frac{r}{2}}$ for any such $r$ (where for any given $x$ there is a trade-off between $x^{-\frac{r}{2}}$ and the pre-factor).  We illustrate this explicitly (for $r = 3$) in Corollary\ \ref{mastertheorem0ccc1} below, and later (in Theorem\ \ref{tailboundedmoments}) ``optimize this bound" (solving for the optimal $r$ for any given $x$) under the assumption that $\E[(S \mu_S)^r]$ satisfies a mild growth rate assumption, yielding a much stronger tail bound.
\item Our bounds have a very different behavior depending on whether $r^*$ is very close to 2 or $r^*$ is significantly larger than 2.  The choice of setting a cutoff at 2.5 was somewhat arbitrary, and simply to make the statement of our results more clear.  Note that the relevant function of $r$ appearing within the corresponding infimum diverges (albeit in different ways) as $r \downarrow 2$ and as $r \rightarrow \infty$.  Due to the fact that (in Theorem\ \ref{mastertheoremlarger}) one takes the $\inf$ of the resulting expression over $r \in [2.5, r^*)$, the divergence as $r \rightarrow \infty$ is not as fundamental of a problem (as one can, for each $x$, apply the bound for any $r \in [2.5,r^*))$.  In contrast, as $r^* \downarrow 2$, the $\inf$ does not remedy the situation, and we leave it as a very interesting open question whether such a degradation is fundamental or merely an artifact of our approach.  Let us again point out that in the single-server case, no such degradation occurs, and one need only assume $\E[S^2] < \infty$.
\item Note that as one assumes the existence of more moments for $S$, the bounds generally become tighter, as the infimum is over a larger range.  In our analysis we prove a bound for each assumption of the form $\E[S^r] < \infty$, and the infimum thus arises since $\E[S^r] < \infty \rightarrow \E[S^{r'}] < \infty$ for all $r' \in [0,r]$.  
\item Let us point out that our analysis actually implies that in Theorem\ \ref{mastertheoremlarger} (i.e. the setting $r^* > 2.5$), we could have taken the infimum to also include the bounds of Theorem\ \ref{mastertheoremsmallr} (i.e. the setting $r^* \leq 2.5$), in which case the infimum would have been over a more complex piecewise function (we have not stated the results this way to improve readability).
\item Note that if $r^* < \infty$ and $\E[S^{r^*}] < \infty$, then the (easily verified) continuity of our bounds implies one can ``plug in" $r^*$ to derive a bound with tail decay rate $x^{-\frac{r^*}{2}}$, even though in principle $r^*$ itself is excluded from the range over which the $\inf$ is taken.  
\item Note the asymmetry of our bounds in $A$ (only finite second moment required, term involving $A$ exhibits exponential decay) and $S$ (finite $r^* > 2$ moment required, term involving $S$ exhibits power law decay depending on $r^*$).  Related discrepancies appear in past results on e.g. existence of moments for the queue length in single and multi-server queues (see e.g. \citet{Kiefer56, Scheller03, Scheller06}).  Intuitively, this manifests because an occasional very large inter-arrival time actually helps the system in some sense, while a large service time will cause the queue to build.  More formally, in our proof we first use a union bound (Lemma\ \ref{2partbound}) to separate our analysis into a term involving the arrivals and a term involving the services, and then observe that the term involving the arrivals contains the supremum of a random walk in which all jumps up are uniformly bounded (even if $A$ itself is not).
\item The precise relationship between existence of moments of the queue length / tail decay rate, existence of moments of $A,S$, and scaling in $\frac{1}{1-\rho}$ remains a very interesting open question.  We note that even for the single-server queue, such questions can become very subtle when $S$ has a sufficiently heavy tail, and to our knowledge such simple and explicit bounds for the tail of the queue length have not appeared before in the literature even in the single-server case.  We refer the interested reader to \citet{Olvera.11,Olvera.11b,Whitt.00,Abate94} for some related discussion, and to \citet{Gaunt20,K64,huang16,Kollerstrom81} for related results in the single-server setting.  
We note that the tail decay rate implied by our main results have several natural implications for the existence and scaling of higher order moments of $L(\infty)$, and for completeness we state and discuss such an implication in Section\ \ref{highermomentsec} of the supplemental appendix.
\item Note that the bound for the s.s.p.d. appearing in \citet{wang2021zero} for systems with hyper-exponentially distributed service times (with $K$ components) and Markovian inter-arrival times has a related dependence on $n(1 - \rho)^2$ as our bounds for the s.s.p.d., as the bound in that work equals $3 K \big(n (1 - \rho)^2 \big)^{-\frac{1}{2}} + \big(n (1 - \rho) \big)^{-1}$.  In our results, the inter-arrival and service times may be from a general distribution, and the demonstrated decay can be an arbitrarily high power of $n(1 - \rho)^2.$  Our optimized bounds (in Theorem\ \ref{tailboundedmoments} below) under stronger assumptions lead to a much faster decay in $n(1 - \rho)^2.$
\item Note that for multi-server queues in the Halfin-Whitt scaling regime, $n(1 - \rho)^2$ is exactly the square of the spare capacity parameter $B$.  That the s.s.p.d. would grow small as $n(1 - \rho)^2$ grows large is thus consistent with past results in the Halfin-Whitt scaling regime, in which the s.s.p.d. grows small as the spare capacity parameter $B$ grows large (see e.g. \citet{HW.81,G16}).  For completeness, we state and discuss some implications of our main results for queues in the Halfin-Whitt regime in Section\ \ref{HalfinChawlasec} of the supplemental appendix.

\item Let us also note that although the the prefactor arising in the bounds of Theorem\ \ref{mastertheoremlarger} involves large constants, these constants (and terms scaling only exponentially in $r$) are asymptotically dominated by a term scaling roughly as $r^{2.5 r} + r^{1.5 r} \times \E[S^r]$.  As $\E[S^r]$ will scale as $r^{\Theta(r)}$ for many $S$ of interest, this fact will allow us to ``optimize" our bounds (by selecting the best $r$ for each $x$) to yield much stronger results.  These results appear later in Section\ \ref{strongersec}, and we include an in-depth discussion of the $r^{\Omega(r)}$ scaling, its necessity in closely related bounds, and interesting related open questions in Section\ \ref{necessarysec}.  Let us also note that we take great care in our results and analysis to separate out and treat the terms that scale as $r^{\Omega(r)}$ or $\E[S^r]$ (which will have the same $r^{\Omega(r)}$ scaling in many cases), as these terms will dominate our bounds asymptotically (in contrast to terms scaling as $c^r$ or $(\E[S^2])^r$).
\end{itemize}
\subsection{Additional implications of main results.}\label{additionalresultssec}
We now present several implications and extensions of our main results for illustrative purposes.
\subsubsection{Outline of our presentation of additional implications of main results.}\label{suboutlinesecc}
We now briefly overview the additional implications of our main results which we will present in the sections below.  In Section\ \ref{illustratesec}, we state two explicit and concrete bounds implied by our main results for illustrative purposes (which do not involve an infimum over $r$), as well as the implied bounds for the expected queue length (which scale as $\frac{1}{1-\rho}$).  In Section\ \ref{strongersec}, we actually compute the infima appearing in Theorem\ \ref{mastertheoremlarger} under additional assumptions on the moments of $S$ (which will hold for all but very heavy-tailed service distributions), which implies much stronger tail bounds.  In Section\ \ref{betterthansec}, we show that our results actually imply a scaling better than $\frac{1}{1-\rho}$ in certain asymptotic regimes.  In Section\ \ref{numbusybsec}, we show that our results also imply bounds for the number of busy servers.

\subsubsection{Illustrative corollaries and bounds for the expected queue length.}\label{illustratesec} An important point is that Theorems\ \ref{mastertheoremsmallr} and \ref{mastertheoremlarger} also imply bounds for the expected queue length $\E[L(\infty)],$ by integrating the bounds for $\pr\big(  L(\infty) \geq \frac{x}{1 - \rho} \big)$ (also using the fact that any probability is always bounded by 1, so the divergence of the bounds as $x \downarrow 0$ is not a problem).  To most accurately state the implied bounds for $\E[L(\infty)]$ (which scale as $\frac{1}{1-\rho}$, as in Kingman's bound), the relevant expression should be the integral of an infimum over $r$ (to derive a bound for each $x$, which is then integrated).  The associated statement is somewhat cumbersome, and thus we do not state it here out of considerations of readability.  Instead, we present two illustrative implications of our main results, both of which follow from Theorems\ \ref{mastertheoremsmallr} - \ref{mastertheoremlarger} and simple integration.
\\\\First, in the corollary below we concretely illustrate our results (also the implied bound for $\E[L(\infty)]$) under the assumption that $\E[S^3] < \infty$ (where higher moments may or may not exist).  These results follow by replacing the infimum (over $r$) in Theorem\ \ref{mastertheoremlarger} by the choice $r = 3$ (implicitly using continuity as $r = 3$ would not technically be included in the infimum in the borderline case $r^* = 3$).
To be clear, Theorem\ \ref{mastertheoremlarger} actually implies a stronger result when $\E[S^3] < \infty$ (due to the infimum).  However, the infimum can be a bit hard to interpret without additional assumptions (although becomes quite interpretable under additional assumptions on the moment sequence of $S$, as in our Theorem\ \ref{tailboundedmoments}).

\begin{corollary}[Illustration of main results when $S$ has finite third moment] \label{mastertheorem0ccc1}
Suppose that for a $GI/GI/n$ queue with inter-arrival times having the same distribution as r.v. $A,$ and service times having the same distribution as r.v. $S$, the following is true : (1)\ $\E[A^2] < \infty$; (2)\ $\E[S^3] < \infty$; (3)\ $\mu_A < n \mu_S$; (4)\ $Q(\infty)$ exists.  Then for all $x > 0$, $\pr\big(  L(\infty) \geq \frac{x}{1 - \rho} \big)$ is at most
\begin{eqnarray*}
\ &\ &\ 8 \times 10^{25} \times \E[(S \mu_S)^2] \times \bigg( \big(\E[(S \mu_S)^2]\big)^2 + \E[(S \mu_S)^3] \bigg) \times x^{-1.5}
\\&\ \ \ &\ \ +\ \ 1.1 \times \exp\bigg( - .0225 \big( \E[ (A \mu_A)^2 ] \big)^{-1} x \bigg);\end{eqnarray*}
the s.s.p.d. is at most
\begin{eqnarray*}
\ &\ &\ 8 \times 10^{28} \times \E[(S \mu_S)^2] \times \bigg( \big(\E[(S \mu_S)^2]\big)^2 + \E[(S \mu_S)^3] \bigg) \times \big( n (1 - \rho)^2 \big)^{-1.5}
\\&\ \ \ &\ \ +\ \ 1.1 \times \exp\bigg( - .0028 \big( \E[ (A \mu_A)^2 ] \big)^{-1} n (1 - \rho)^2 \bigg);
\end{eqnarray*}
and $\E\big[ L(\infty) \big]$ is at most
$$\Bigg( 1.61 \times 10^{26} \times \E[(S \mu_S)^2] \times \bigg( \big(\E[(S \mu_S)^2]\big)^2 + \E[(S \mu_S)^3] \bigg) + 
49 \E[ (A \mu_A)^2 ] \Bigg) \times \frac{1}{1-\rho}.$$
\end{corollary}

Of course, there was nothing special about the number 3 in Corollary\ \ref{mastertheorem0ccc1}, and an analogous result can be easily derived for any number strictly greater than 2.
\\\\As a second illustration, we now state the bounds implied by Theorem\ \ref{mastertheoremsmallr} when one only assumes that $\E[S^{2 + \epsilon}] < \infty$ for some small $\epsilon$ (where higher moments may or may not exist).  Here our result follows directly from Theorem\ \ref{mastertheoremsmallr} and some straightforward algebra, along with the easily verified fact that $(\frac{1}{\epsilon})^{\epsilon} < 1.5$ for all $\epsilon > 0$.

\begin{corollary}[Illustration of main results when $S$ has finite $2 + \epsilon$ moment]\label{corsmallr}
Suppose that for a $GI/GI/n$ queue with inter-arrival times having the same distribution as r.v. $A,$ and service times having the same distribution as r.v. $S$, the following is true : (1)\ $\E[A^2] < \infty$; (2)\ there exists $\epsilon \in (0,.5)$ s.t. $\E[S^{2+\epsilon}] < \infty$; (3)\ $\mu_A < n \mu_S$; (4)\ $Q(\infty)$ exists.  Then for all $x > 0$, $\pr\big(  L(\infty) \geq \frac{x}{1 - \rho} \big)$ is at most
\begin{eqnarray*}
\ &\ &\ 5.1 \times 10^{20} \times \E[(S \mu_S)^2] \times \bigg( \big(\E[(S \mu_S)^2]\big)^{1 + \epsilon}  + \E[(S \mu_S)^{2 + \epsilon}] \bigg) \times (\frac{1}{\epsilon})^3 \times x^{-(1+\frac{\epsilon}{2})} \nonumber
\\&\ \ \ &\ \ +\ \ 1.1 \times \exp\bigg( - .0225 \big( \E[ (A \mu_A)^2 ] \big)^{-1} x \bigg);
\end{eqnarray*}
the s.s.p.d. is at most
\begin{eqnarray*}
\ &\ & 6.8 \times 10^{21} \times \E[(S \mu_S)^2] \times \bigg( \big( \E[(S \mu_S)^2] \big)^{1 + \epsilon} + \E[(S \mu_S)^{2 + \epsilon}] \bigg) \times (\frac{1}{\epsilon})^3 \times \big(n (1 - \rho^2) \big)^{-(1+\frac{\epsilon}{2})}  \nonumber
\\&\ \ \ &\ \ +\ \ 1.1 \times \exp\bigg( - .0028 \big( \E[ (A \mu_A)^2 ] \big)^{-1} n (1 - \rho)^2 \bigg);
\end{eqnarray*}
and $\E\big[ L(\infty) \big]$ is at most
$$
\Bigg(
2.1 \times 10^{21} \times \E[(S \mu_S)^2] \times \bigg( \big(\E[(S \mu_S)^2]\big)^{1 + \epsilon}  + \E[(S \mu_S)^{2 + \epsilon}] \bigg) \times (\frac{1}{\epsilon})^4
+
49 \E[ (A \mu_A)^2 ] \Bigg) \times \frac{1}{1-\rho}.$$
\end{corollary}
As mentioned before, we leave it as a very interesting open question whether the divergence of these bounds (in $\epsilon$ as $\epsilon \downarrow 0$) is fundamental or merely an artifact of our approach.  
\\\\ Let us point out that our main results also imply bounds for the higher moments of $L(\infty)$ (scaling as an appropriate power of $\frac{1}{1-\rho}$), and for completeness we provide such a bound in the supplemental appendix Section\ \ref{highermomentsec}.

\subsubsection{Stronger explicit bounds when $r^* = \infty$ and $S$ is not very heavy tailed.}\label{strongersec} We now show that when all moments of $S$ exist, and this moment sequence satisfies a mild technical condition (which in general will hold for all but very heavy-tailed distributions), we can derive much tighter bounds.  These results follow by actually computing the infima appearing in Theorem\ \ref{mastertheoremlarger}, i.e. optimizing the bound (over $r$) for each fixed $x$.
\\\indent To motivate the precise conditions we will impose on the moments of $S$, let us first make an observation pointing out that all but very heavy-tailed distributions satisfy a natural growth condition on their moment sequence.  Namely, under this mild condition one will have $\E[S^r] = r^{O(r)}$.  The observation follows from standard results for the moments of a Weibull distribution and some straightforward algebra, and for completeness we include a proof in the supplemental appendix Section\ \ref{optrsec}.
\begin{observation}[If $S$ is not very heavy tailed then moment sequence scales as $r^{O(r)}$]\ \label{limitedgrowth}
Suppose there exist $a,b > 0$ and $c \leq 1$ s.t. $P(S \mu_S > x) \leq a \exp( - b x^c)$ for all $x > 0$.  Then 
$\E[(S\mu_S)^r] \leq 1.5 \times a \times \big( (b c)^{-\frac{1}{c}}\big)^r \times r^{\frac{1}{c} r}$ for all $r \geq 2$.  It follows that as long as there exists $\alpha > 0$ s.t. the tail of $S \mu_S$ decays (asymptotically) at least as fast as $\exp\big(-x^{\alpha})$, then there will exist constants $a',b',c' > 0$ s.t. $\E[(S \mu_S)^r] \leq a' \times b'^r \times r^{c' r} = r^{O(r)}$ for all $r \geq 2$.  These bounds thus apply (for example) if $S$ has exponentially decaying tails (in which case $c = 1$), or even much heavier tails with decay rate comparable to that of a heavier-tailed Weibull distribution (in which case $c \in (0,1)$).  We note that the case $P(S \mu_S > x) \sim a \exp(- b x^c)$ for $c \in (0,1)$ puts the associated queueing model beyond the scope of several past results which assume that $\sup_{t > 0} \E[S \mu_S - t | S \mu_S > t] < \infty$ (\citet{Downey.91,grosof2021finite}).
\end{observation}

The next result shows that if the moments of S can be bounded as in Observation\ \ref{limitedgrowth}, our results imply a much stronger tail bound for $GI/GI/n$ queues.  We include a proof in the supplemental appendix Section\ \ref{optrsec}.

\begin{theorem}[Stronger bounds if moments of S scale as $r^{O(r)}$]\label{tailboundedmoments}
Suppose there exist $a,b,c \geq 1$ s.t. $\E[(S \mu_S)^r] \leq a \times b^r \times r^{c r}$ for all $r \geq 2$.  Suppose also that the assumptions of Theorem\ \ref{mastertheoremlarger} hold.  Let $f \stackrel{\Delta}{=} 4 \times 10^4 \times a, g \stackrel{\Delta}{=} 10^7 \times a \times b^3 \times 4^c, \delta \stackrel{\Delta}{=} \frac{1}{e} \times (1.5 + c) \times (\frac{1}{g})^{\frac{1}{1.5 + c}}.$
\\\\Then for all $x > 0$, $\pr\big(  L(\infty) \geq \frac{x}{1 - \rho} \big)$ is at most 
$$f \times \exp\big( - \delta x^{\frac{1}{3 + 2 c}} \big) + 1.1 \times \exp\bigg( - .0225 \big( \E[ (A \mu_A)^2 ] \big)^{-1} x \bigg);$$
and the s.s.p.d. is at most
$$f \times \exp\big( - \delta \big(n (1 - \rho)^2 \big)^{\frac{1}{3 + 2 c}} \big) + 1.1 \times \exp\bigg( - .0028 \big( \E[ (A \mu_A)^2 ] \big)^{-1} n (1 - \rho)^2 \bigg).$$
\end{theorem}

Let us note that although $f$ and $\delta$ involve large constants, $c$ does not, and will in general be a small integer (depending on the actual growth rate of the moments of $S$).  Thus (for example) if $S$ has exponentially decaying tails, then $c = 1$ and $3 + 2 c = 5$.  Thus in this case, our results imply an explicit tail decay rate which scales asymptotically as $\exp(- x^{\frac{1}{5}})$ (note that the term arising from the arrival process will always have an asymptotically faster decay rate).  Let us also point out that although for illustrative purposes we have focused here on the assumption that $\E[(S \mu_S)^r] \leq a \times b^r \times r^{c r}$, analogous results (with different tail decay properties) could be derived under any growth rate assumption on the relevant moment sequence.  In addition, we note that these results regarding the s.s.p.d. (as well as our other results for the s.s.p.d.), when applied to queues in the Halfin-Whitt regime, imply progresss on an open question posed in \citet{Chawla17} on distribution-independent bounds for queues in the Halfin-Whitt regime.  We defer a more formal discussion of the connection to \citet{Chawla17} to the supplemental appendix Section\ \ref{HalfinChawlasec}.

\subsubsection{Better than $\frac{1}{1-\rho}$ scaling in certain asymptotic regimes.}\label{betterthansec}
We now observe that by integrating the minimum of our bound for $\pr\big(  L(\infty) \geq \frac{x}{1 - \rho} \big)$ and our bound for the s.s.p.d. (which also yields a bound for the tail of $L(\infty)$), we can obtain a scaling better than $\frac{1}{1 - \rho}$.  For illustrative purposes and clarity of exposition, we state this result under the assumption that $\E[S^3] < \infty$, although analogous results (which decay as different powers of $n(1-\rho)^2$) can be derived for any $r$ s.t. $\E[S^r] < \infty$,  and even stronger results could be derived under assumptions analogous to those of Theorem\ \ref{tailboundedmoments}.  For completeness, we include a proof in the supplemental appendix Section\ \ref{meanwaitcorbeyondsec}.  To make a comparison to our previous bound for $E[L(\infty)]$ appearing in Corollary\ \ref{mastertheorem0ccc1} easier, let us define $a_1 \stackrel{\Delta}{=} 1.61 \times 10^{26} \times \E[(S \mu_S)^2] \times \bigg( \big(\E[(S \mu_S)^2]\big)^2 + \E[(S \mu_S)^3] \bigg), a_2 \stackrel{\Delta}{=} 49 \E[ (A \mu_A)^2 ] $.  Note that our previous Corollary\ \ref{mastertheorem0ccc1} asserted that under appropriate assumptions, $E[L(\infty)] \leq (a_1 + a_2) \times \frac{1}{1-\rho}$.

\begin{corollary}[Better than $\frac{1}{1 - \rho}$ scaling]\label{meanwaitcorbeyond}
Under the same assumptions as Corollary\ \ref{mastertheorem0ccc1}, and supposing also that $n(1 - \rho)^2 \geq 10^6 \times \big(\E[(A \mu_A)^2]\big)^2$, it holds that $\E\big[ L(\infty) \big] \leq  (a_1 + a_2) \times \frac{1}{1-\rho} \times 1000 \big( n(1-\rho)^2\big)^{-.5}.$  
\end{corollary}

Corollary\ \ref{meanwaitcorbeyond} goes beyond the $\frac{1}{1 - \rho}$ scaling, with an additional correction term $1000 \big( n(1-\rho)^2\big)^{-.5}$ which converges to 0 as $n (1 - \rho)^2$ grows large.  Such a bound can be interpreted as a generalization of the fact that in an M/M/n queue, $\E[L(\infty)] = \pr\big( Q(\infty) \geq n \big) \times \frac{\rho}{1 - \rho}$, with the term $\big( n (1 - \rho)^2 \big)^{-.5}$ acting as a proxy for $\pr\big( Q(\infty) \geq n \big)$.  We note that such a term will be significant in certain asymptotic scaling regimes, e.g. in the Halfin-Whitt regime when the spare capacity parameter $B$ is large, or in the quality-driven scaling regime for multi-server queues (in which $n \rightarrow \infty$ for a fixed $\rho$, see e.g.  \citet{Borst04,Mag18,Whitt.82d,HW.81}).

\subsubsection{Bound for number of busy servers in certain scaling regimes.}\label{numbusybsec}
In certain scaling regimes, the primary metric of interest will be the number of busy servers, as with high probability there is no queue.  Our next result provides bounds in exactly this setting, and we include a proof in Section\ \ref{sspdsec}.

\begin{theorem}[Bound for number of busy servers]\label{numbusybbb}
Under the same assumptions as Theorem\ \ref{mastertheoremlarger}, and supposing also that $\rho \in [\frac{1}{n}, 1 - \frac{2}{n}],$ the following is true.  For all $x \in \bigg[ 4, \frac{1}{2} \min\big(\sqrt{\frac{\mu_A}{\mu_S}}, \frac{n - \frac{\mu_A}{\mu_S}}{ \sqrt{\frac{\mu_A}{\mu_S}} }\big) \bigg]$, it holds that $\pr\big( \textrm{Num}_{\textrm{service}}(\infty) \geq \frac{\mu_A}{\mu_S} + x \sqrt{\frac{\mu_A}{\mu_S}} \big)$ is at most
\begin{eqnarray*}
\ &\ &\ \inf_{r \in [2.5,r^*)} \Bigg( 2 \times 10^4 \times \E[(S \mu_S)^2] \times \big(4 \times 10^6\big)^r \times \bigg( \big(\E[(S \mu_S)^2]\big)^{r-1} \times r^{2.5 r} + r^{1.5 r} \times \E[(S \mu_S)^r] \bigg) \times x^{-r} \Bigg)
\\&\ \ \ &\ \ +\ \ 1.1 \times \exp\bigg( - .00015 \big( \E[ (A \mu_A)^2 ] \big)^{-1} x^2 \bigg);
\end{eqnarray*}
\end{theorem}

We note that in the quality-driven regime (i.e. $n \rightarrow \infty$ for fixed $\rho$), for any given $M > 0$, $\min\big(\sqrt{\frac{\mu_A}{\mu_S}}, \frac{n - \frac{\mu_A}{\mu_S}}{ \sqrt{\frac{\mu_A}{\mu_S}} }\big)$ will be greater than $M$ for all sufficiently large $n$.  Furthermore, $\frac{\mu_A}{\mu_S} + x \sqrt{\frac{\mu_A}{\mu_S}}$ is at most $n$ for $x$ in the same range.  Thus for any fixed $r^*$ and all sufficiently large $x$ (with ``sufficiently large" a function of $r^*$ and the moments of $S \mu_S, A \mu_A$), Theorem\ \ref{numbusybbb} will yield meaningful (i.e. strictly less than unity) bounds on the probability that the number of busy servers exceeds $\frac{\mu_A}{\mu_S} + x \sqrt{\frac{\mu_A}{\mu_S}}$ in the quality-driven regime.  Such a result is consistent with the known Poisson approximations for the number in service in these scaling regimes (see e.g. \citet{Borst04,Mag18,Whitt.82d,HW.81}).  By a similar logic, our results can also yield meaningful (i.e. strictly less than unity) bounds for the number of busy servers in the Halfin-Whitt scaling regime when the spare capacity parameter $B$ is sufficiently large.  We further note that a stronger result could again be proven under the same moment assumptions as Theorem\ \ref{tailboundedmoments}, although we do not formalize that here.

\subsection{Additional discussion of prefactors and ``universality" of our main results.}\label{additionaldiscusssec}
In this section we present some additional discussion of whether our main results provide ``universal" bounds (spoiler - they do not), and provide some additional context and supporting results to help explain the (large) prefactors appearing in our main results. 

\subsubsection{Do these results provide ``universal" bounds?}\label{isuniversal}
One can ask whether our bounds are in a meaningful sense ``universally accurate" across all scalings of $n$ and $\rho$.  Although our bounds do represent a step in that direction, they unfortunately fall short in this regard.  
\begin{itemize}
\item The tail decay rate implied by our bounds (also our Theorem\ \ref{tailboundedmoments}) does not match known results in several cases where the asymptotic behavior of the tail of $L(\infty)$ is well-understood (including for single-server queues).
\item It follows from known results for queues in the Halfin-Whitt regime (\citet{HW.81, G16}), super-Halfin-Whitt regime (\citet{hong2021sharp}), and light traffic regimes (\citet{Burman83,Daley.92,Gupta11}), that a faster decay in $n(1-\rho)^2$ holds for the s.s.p.d. in certain scaling regimes.  
\item Our bounds are not applicable in a very light traffic regime where $\rho \downarrow 0$ for fixed $n$, as they do not converge to zero as $\rho \downarrow 0$.  
\item The large prefactors appearing in our bounds render them meaningless for certain parameter ranges, and no doubt extremely conservative in many cases (e.g. one can compare our bounds to known results for the $M/M/1$ queue).  
\end{itemize}

\subsubsection{On the large prefactors appearing in our bounds.}\label{necessarysec}
Let us now address the proverbial ``elephant in the room" - namely, the massive prefactors in our results.  These prefactors involve both large numerical constants, as well as functions of $r$ whose asymptotics (at least for large $r$) are dominated by terms of the form $r^{c r}$ for some small constant $c$.  More precisely, $c = \max(2.5, 1.5 + \gamma)$ for some explicit $\gamma$ depending on the growth rate of the moments of $S$, where $\gamma = 1$ if $S$ has exponentially decaying tails, and will equal $\alpha$ if the tail of $S$ decays as $\exp(-x^{\frac{1}{\alpha}})$ as discussed in Observation\ \ref{limitedgrowth}.  It is natural to ask whether such a large prefactor is fundamental, or merely an artifact of our approach.  We provide some relevant discussion here, and include a more detailed discussion in Section\ \ref{rnotsmallsec} of the supplemental appendix.
\\\indent First, let us point out that for any given queueing system, the prefactors appearing in our bounds may be very loose (for example comparing our bounds to known results for $M/M/1$ queues).  Second, regarding the explicit numerical constants (e.g. $10^4$ or $10^6$) appearing in our bounds, these arise not from any one or two particular aspects of the proof, but instead from the composition of many bounds and results from the literature.  Although we did take some care to try and control these constants in our analysis, it is likely that they could be further improved by an even more careful analysis, or possibly by a completely different type of analysis (see Section\ \ref{noprefacsec}).  
\\\indent The asymptotic scaling (in $r$) of the terms and prefactors appearing within the infima arising in our main results leads to some subtle and interesting questions.  It is easy to see that in our main results, those functions of $r$ scale as $r^{O(r)}$ for large $r$ (so long as $\E[(S \mu_S)^r]$ scales as $r^{O(r)}$).  Although the fact that one takes an infimum over these terms means that these prefactors do not in general appear in the bound for any given $x$, they do appear in the bound if one insists on acheiving the highest possible tail decay rate (as a function of $r$) under minimal assumptions (i.e. how many moments of $S$ exist).   For example, our main result Theorem\ \ref{mastertheoremlarger} directly implies the following (we formally prove this in Section\ \ref{rnotsmallsec} of the supplemental appendix).
\begin{corollary}\label{rnotsmalltheorema}
For each $r > 2.5$, there exists a finite constant $c_r$ (depending only on $r$) s.t. for all $GI/GI/n$ queues in which $\E[A^r] < \infty, \E[S^r] < \infty$, $\mu_A < n \mu_S$, and $Q(\infty)$ exists, it holds that for all $x > 0$, $P\big( L(\infty) \geq \frac{x}{1 - \rho}\big) \leq c_r \bigg( \big(\E[(A \mu_A)^2]\big)^r \times \E[(A \mu_A)^r] + \big(\E[(S \mu_S)^2]\big)^r \times \E[(S \mu_S)^r] \bigg) \times x^{-\frac{r}{2}}$.  In addition, one can take $c_r < 4 \times 10^4 \times \big(10^6 \big)^r \times r^{2.5 r}$ for all $r > 2.5$.  
\end{corollary}
Note that $4 \times 10^4 \times \big(10^6 \big)^r \times r^{2.5 r}$ is asymptotically dominated by the $r^{2.5 r}$ term.  As for any $\delta > 0$ the function $r^{\delta r}$ grows quite rapidly (faster than any exponential in $r$), it is reasonable to ask whether such a growth is fundamental in the prefactors arising in e.g. Corollary\ \ref{rnotsmalltheorema}, or is merely an artifact of our analysis.  It turns out that such a scaling is indeed fundamental, as indicated by the following theorem which we prove in Section\ \ref{rnotsmallsec} of the supplemental appendix.
\begin{theorem}\label{rnotsmalltheoremb}
For each $r > 2.5$, let $\underline{c_r}$ denote the infimum of all constants $c_{r}$ for which the bound of Corollary\ \ref{rnotsmalltheorema} holds.  Then there exists an absolute constant $\epsilon > 0$ s.t. $\underline{c_{r}} \geq \epsilon \times r^{\frac{1}{9} r}$ for all $r > 16$.
\end{theorem} 
The intuition for the above is actually quite straightforward, and arises from the fact that when $A,S$ are uniformly bounded it holds that $\E[A^r],\E[S^r]$ scale only exponentially in $r$, while in a single-server queue the $r$th moment of the steady-state queue length (with inter-arrival times distributed as $A$ and service times distributed as $S$) scales as $r^{\Omega(r)}$, inherited from the scaling of the moments of an exponential distribution.  
\\\indent In the proofs of our main results, the prefactors scaling as $r^{\Omega(r)}$ arise largely from our bounds for the higher order central moments of (pooled) renewal processes, i.e. $\E\big[|N_e(t) - \mu_S t|^r\big]$ and $\E\big[|\sum_{i=1}^n N_{e,i}(t)  - \mu_S n t|^r\big]$, often inherited from other results in the literature.  We also show in Section\ \ref{rnotsmallsec} of the supplemental appendix that the relevant prefactors in these intermediate bounds must indeed have a $r^{\Omega(r)}$ scaling.
\\\indent It turns out that the lower bound analysis used in our proof of Theorem\ \ref{rnotsmalltheoremb} does not apply to our actual bounds from our main result Theorem\ \ref{mastertheoremlarger}, with the discrepancy arising from the fact that in our actual bounds there is a term of the form $1.1 \times \exp\bigg( - .0225 \big( \E[ (A \mu_A)^2 ] \big)^{-1} x \bigg)$ instead of a term of the form $\big(\E[(A \mu_A)^2]\big)^r \times \E[(A \mu_A)^r] \times x^{-\frac{r}{2}},$ and we provide further details in Section\ \ref{rnotsmallsec} of the supplemental appendix.  More broadly, the required prefactors could in principle be quite sensitive to the exact form of the desired bound.  This may be especially true if the desired bounds are to capture the ``best possible" tail decay rate.  For example, it is shown in \citet{Abate94} that in this case the relevant prefactor may depend on explicit terms arising in the asymptotics of the tail of the service time distribution, which may not be boundable (even in principle) in terms of the moments of $S$.  Also, as noted before, such uniform tail bounds are known to be very complex even in the single-server setting, when $S$ has a heavy tail (see e.g. \citet{Olvera.11,Olvera.11b}).
\\\indent We leave the question of deriving tighter bounds, possibly under different assumptions and/or in the study of qualitatively different types of bounds, as a very interesting direction for future research.  We take some steps along these lines with our results in Section\ \ref{noprefacsec}.
\section{Towards bounds with no large prefactors.}\label{noprefacsec}
Although it remains an interesting open question whether bounds analagous to our main results (including for the mean queue length) in $GI/GI/n$ queues exist with no large prefactors, here we provide some partial evidence that this may indeed be possible.  In particular, we prove such bounds in three related settings of interest : (1) the s.s.p.d. when the arrival process is Markovian; (2) the tail of the number of busy servers when the arrival process is Markovian; and (3) the mean queue length when the arrival process is Markovian and there are Markovian abandonments.  We also provide an intuitive conjecture which would imply such a bound for the mean queue length in $M/GI/n$ queues.  Our conjecture is closely inspired by a known relationship for the workload in $M/GI/n$ queues, and is similar to related conjectures appearing previously in the literature.
\\\indent Our proofs of these results are fundamentally different from those of our other results.  These results follow from a simplified version of drift arguments very similar to those appearing in \citet{wang2021zero,hong2021sharp, hokstad85, scully2020gittins, grosof2021finite}.  We note that several of these works were primarily focused on the study of more general models in which jobs can utilize more than one server and/or use other service disciplines.  We defer all proofs to Section\ \ref{noprefactorsec0} of the supplemental appendix.
\\\indent Our first result, for the s.s.p.d., is as follows.  We defer the proof to Section\ \ref{noprefactorsec1} of the supplemental appendix.
\begin{theorem}\label{sspdnoprefactor}
Consider an $M/GI/n$ queue ${\mathcal Q}^n$ with Markovian inter-arrival times having the same distribution as r.v. $A,$ and service times having the same distribution as r.v. $S$ satisfying $\E[S] < \infty$.  Suppose also that $\mu_A < n \mu_S$, and that $Q(\infty)$ and $\overline{W}^n_{\textrm{service}}(\infty)$ exist.  Then the s.s.p.d. is at most $\frac{1}{2} \sqrt{\rho} \big( n(1-\rho)^2 \big)^{-\frac{1}{2}}.$
\end{theorem}
\ \\The bound also has a similar qualitative dependence on $n(1-\rho)^2$ as several of our previous results.  In contrast to our previous results for the s.s.p.d., here no assumption on $S$ is needed except that $\E[S] < \infty,$ and the bound converges to zero as $\rho \downarrow 0$ for fixed $n.$  However, the bound does not decay faster in $n(1 - \rho)^2$ under the assumption of additional moments (on $S$), and requires inter-arrival times to be Markovian.  
\\\indent Our second result, for the steady-state number of busy servers, is as follows.  We defer the proof to Section\ \ref{noprefactorsec15} of the supplemental appendix.
\begin{theorem}\label{sspdnoprefactora2}
Under the same assumptions as Theorem\ \ref{sspdnoprefactor}, for all $x \in (0, \frac{n - \frac{\mu_A}{\mu_S}}{\sqrt{\frac{\mu_A}{\mu_S}}}]$, $P\big( \textrm{Num}_{\textrm{service}}(\infty) \geq \frac{\mu_A}{\mu_S} + x \sqrt{\frac{\mu_A}{\mu_S}} \big) \leq \frac{1}{2 x}.$
\end{theorem}
\ \\The bound again has no large prefactors, but is restricted to Markovian arrivals, and unable to demonstrate faster decay rates in $x$ under the assumption that more moments of $S$ are finite.  
\\\indent Our next result is for $M/PH/n + M$ queues, i.e. multi-server queues with Markovian inter-arrival times, Markovian abandonments, and phase-type service times.  For this result we restrict to the class of phase-type service distributions, which are well-known to be dense within the space of all distributions (see \citet{AS96}).  This restriction enables us to apply certain technical results regarding e.g. existence and construction of the relevant stationary measures and processes.  However, our actual bounds will not depend on any parameters of the phase-type distribution beyond the mean-interarrival and service time and abandonment rate, and could likely be extended to general service times by simple continuity arguments (using e.g. the results of \citet{W74}).  Several past works have studied such queueing systems with abandonments using Lyapunov arguments (see e.g. \citet{GS12,DDG14,BD15}), but to our knowledge such a simple and explicit bound has not appeared previously in the literature.
\\\indent First, let us more formally describe the relevant system.  Let ${\mathcal Q}^n_a$ be a $M/PH/n + M$ multi-server queue with Markovian inter-arrival times distributed as (the exponentially distributed r.v.) $A$, service times distributed as the (phase-type) r.v. $S$, and patience times distributed as the (exponentially distributed) r.v. $B$ with $\theta = \frac{1}{\E[B]}$.  Suppose the queue is initially empty.  For a more formal description of such multi-server queues with abandonments, we refer the reader to e.g. \citet{GS12,DDG14,DH13,BD15,Mandelbaum2012}.  It follows from the results of \citet{DDG14} that such a system is positive Harris recurrent, and  the total number of jobs in ${\mathcal Q}^n_a$ (number in service + number waiting in queue) converges in distribution (as time goes to infinity, independent of the particular initial condition) to a steady-state r.v. $Q^n_a(\infty)$.  Also, let $L^n_a(\infty)$ denote a r.v. distributed as the steady-state number of jobs waiting in queue.  For such a system with abandonments, we again define $\rho = \frac{\mu_A}{n \mu_S}$.  We note that in general a system with abandonments will be stable even when $\rho > 1$, although here consider only the case $\rho < 1$ in analogy to our results without abandonments.  
\\\indent Then our result is as follows.  We defer the proof to Section\ \ref{noprefactorsec2} of the supplemental appendix.
\begin{theorem}\label{sspdnoprefactorb}
Consider an $M/PH/n + M$ queue ${\mathcal Q}^n_a$ with Markovian inter-arrival times having the same distribution as r.v. $A$, service times having the same distribution as phase-type r.v. $S$, and patience times having the same distribution as exponentially distributed r.v. $B$ with $\theta = \frac{1}{\E[B]}.$  Suppose also that $\mu_A < n \mu_S.$  Then $\E\big[ L^n_a(\infty) \big] \leq \sqrt{\frac{2}{\pi}} \frac{\mu_S}{\theta} \rho \sqrt{n}$.  If in addition $\rho \in [\frac{3}{4}, 1 - \frac{4}{n}]$, then $\E\big[ L^n_a(\infty) \big] \leq 2 \frac{\mu_S}{\theta} \sqrt{n} \exp\big( - \frac{1}{4} n (1 - \rho)^2 \big)$.
\end{theorem}
\ \\Note that the first result above is most applicable for values of $\rho$ in the Halfin-Whitt scaling (when $\frac{1}{1-\rho}$ scales as $\Theta(\sqrt{n})$), while the second result is applicable in heavy-traffic more broadly.  Our second result recaptures a scaling qualitatively similar to that previously shown for other related metrics in different particular asymptotic scaling regimes (see e.g. \citet{wang2021zero} and \citet{G16}).  The second result is also meaningful (i.e. yields bounds strictly less than unity) in the so-called quality-driven regime, with $\rho \geq \frac{3}{4}$ fixed and $n \rightarrow \infty$.  We also note the restriction $\rho \in [\frac{3}{4}, 1 - \frac{4}{n}]$ was chosen as a technical convenience, and different bounds (with a different exponent) could be derived under different assumptions.
\\\indent Unfortunately, as the system with abandonments approaches a system without abandonments (i.e. $\theta \downarrow 0$), the bound of Theorem\ \ref{sspdnoprefactorb} becomes meaningless.  For $M/GI/n$ queues without abandonments, we now present a conjecture which would imply an analogous result.  To state the conjecture, we first state a particular equation for $\E\big[L(\infty)\big]$.  This equation follows directly (after some straightforward algebra) from a known equation for the steady-state expected work-in-system studied via drift arguments in several past works (see e.g. \citet{hokstad85, scully2020gittins, grosof2021finite, wang2021zero, hong2021sharp}), and for completeness we provide a proof in Section\ \ref{noprefactorsec3} of the supplemental appendix.

\begin{theorem}\label{sspdnoprefactorc} Consider an $M/GI/n$ queue ${\mathcal Q}^n$ with with Markovian inter-arrival times having the same distribution as r.v. $A,$ service times having the same distribution as r.v. $S$ satisfying $\E[S^2] < \infty$, s.t. the c.d.f. of $S$ is absolutely continuous.  Suppose also that $\mu_A < n \mu_S$, and that $Q(\infty), \textrm{Work}(\infty),$ and $\overline{W}_{\textrm{service}}(\infty)$ exist.  Then $\E\big[L(\infty)\big]$ equals
$$\frac{1}{2} \E[(S \mu_S)^2] \frac{\rho}{1 - \rho} + \frac{
\E\big[\textrm{Num}_{\textrm{service}}(\infty)\big] \times \E\big[\textrm{Work}_{\textrm{service}}(\infty)\big] - 
\E\big[\textrm{Num}_{\textrm{service}}(\infty) \times \textrm{Work}_{\textrm{service}}(\infty)\big] }{n (1 - \rho) \E[S]}.$$
\end{theorem}

Theorem\ \ref{sspdnoprefactorb} would imply the simple bound $\E\big[L(\infty)\big] \leq \frac{1}{2} \E[(S \mu_S)^2] \frac{\rho}{1 - \rho}$ if the term 
\begin{equation}\label{covar11}
\E\big[\textrm{Num}_{\textrm{service}}(\infty)\big] \times \E\big[\textrm{Work}_{\textrm{service}}(\infty)\big] - 
\E\big[\textrm{Num}_{\textrm{service}}(\infty) \times \textrm{Work}_{\textrm{service}}(\infty)\big]
\end{equation}
was negative.  We note that $\frac{1}{2} \E[(S \mu_S)^2] \frac{\rho}{1 - \rho}$ closely resembles the steady-state expected number in queue in an $M/GI/1$ queue with inter-arrival times having the same distribution as $A$, and service times having the same distribution as $\frac{S}{n}$.  Indeed, that expectation equals $\frac{1}{2} \E[(S \mu_S)^2] \frac{\rho^2}{1 - \rho}.$  The fact that such relationships generally elucidate connections to such a sped-up single-server queue is well-known, see e.g. \citet{hokstad85, scully2020gittins, grosof2021finite}.  However, (\ref{covar11}) is simply the negative of the covariance between the total work in service and the total number in service.  One would intuitively expect this covariance to be positive for non-pathological FCFS $M/GI/n$ systems, and we indeed conjecture that such a result holds for a broad class of $M/GI/n$ queues.  

\begin{conjecture}\label{conjecture1} In any $M/G/n$ queue in which $\E[S^2] < \infty$ and all relevant steady-state distributions exist, it holds that 
$$\E\big[\textrm{Num}_{\textrm{service}}(\infty)\big] \times \E\big[\textrm{Work}_{\textrm{service}}(\infty)\big] \leq \E\big[\textrm{Num}_{\textrm{service}}(\infty) \times \textrm{Work}_{\textrm{service}}(\infty)\big],$$
and hence $\E\big[L(\infty)\big] \leq \frac{1}{2} \E[(S \mu_S)^2] \frac{\rho}{1 - \rho}.$
\end{conjecture}

We leave a further study as an interesting open question, and note that tools from the theory of associated r.v.s (which have been used to prove related results in the literature) may be relevant here (see e.g. \citet{Bac13}).  We note that bounds for this covariance are implicit (or in some cases explicit) in past work (see e.g. \citet{hokstad85, scully2020gittins, grosof2021finite, wang2021zero, hong2021sharp}), but do not seem to have the desired scaling when $n \rightarrow \infty$ and $\rho \uparrow 1$ simultaneously (at least in certain parameter regimes).  We also note that closely related conjectures for multi-server queues involving similar (albeit perhaps less interpretable) covariance terms appear throughout the literature, and we refer the reader to \citet{Mori.75,Daley.97} for additional discussion.

\section{Proof of our main results : bounds for $L(\infty)$ in Theorems\ \ref{mastertheoremsmallr} - \ref{mastertheoremlarger}.}\label{mainproofsec}
In this section we prove our central main results, i.e. the first part of Theorems\ \ref{mastertheoremsmallr} - \ref{mastertheoremlarger} in which $\pr\big( L(\infty) \geq \frac{x}{1-\rho} \big)$ is bounded.  To maximize readability, we proceed as follows.  First, we sketch a high-level outline of the proof in Section\ \ref{highleveloutlinesec}.  Second, we provide a more detailed proof (but still without most technical details), containing all of the most important auxiliary results and main flow of logic (albeit in many cases without their proofs), in Section\ \ref{abitmoredetailsec}.  Third, we provide many of the most important technical details of the proofs (but still with many of the finer subarguments omitted) in the technical appendix Section\ \ref{moredetailsec}.  Finally, we defer many of the finer subarguments of these proofs to the supplemental appendix Section\ \ref{appsec}. 

\subsection{High-level outline.}\label{highleveloutlinesec}
We begin by sketching the high-level outline of our proof of the bounds for $\pr\big( L(\infty) \geq \frac{x}{1-\rho} \big)$ appearing in Theorems\ \ref{mastertheoremsmallr} - \ref{mastertheoremlarger}.
\ \\\begin{enumerate}
\item  \label{step1outline} Use stochastic comparison results of \citet{GG13} to bound $\pr\big(  L(\infty) \geq x \big)$ by 
\begin{equation}\label{mainuboundgg}
\pr\Bigg( \sup_{t \geq 0} \bigg( A_e(t) - \sum_{i=1}^n N_{e,i}(t) \bigg) \geq x \Bigg).
\end{equation}
\item Use a union bound, and the connection between ${\mathcal A}_e$ and ${\mathcal A}_o$, to bound $(\ref{mainuboundgg})$ by
\begin{eqnarray}
\ &\ &\ \ \pr\Bigg( \sup_{t \geq 0} \bigg( A_o(t) - \mu_{A} t - \frac{1}{2}( n - \mu_{A} ) t \bigg) \geq \frac{1}{2} x - 1 \Bigg) \label{apart}
\\&\ &\ \ \ +\ \ \ \pr\Bigg( \sup_{t \geq 0} \bigg( n t - \sum_{i=1}^{n} N_{e,i}(t) - \frac{1}{2}( n - \mu_{A} ) t  \bigg) \geq \frac{1}{2} x \Bigg), \label{spart}
\end{eqnarray}
thus separating the proof into an analysis of a supremum arising from the arrival process (\ref{apart}) and a supremum arising from a centered pooled equilibrium renewal process with renewal intervals distributed as service times (\ref{spart}). \label{step2outline}
\item Bound (\ref{apart}) by relating the supremum to a simple single-server queue and using known bounds for that setting (specifically a martingale inequality proven in \citet{K64}). \label{step3outline}
\item Conditionally bound (\ref{spart}) by proving that IF one could suitably bound $\E[\big|n t - \sum_{i=1}^{n} N_{e,i}(t)|^r\big]$ for some $r > 2$ (and all $t$) THEN one could bound (\ref{spart}) as required (using modifications of known maximal inequalities). \label{step4outline}
\item Combine the bound for (\ref{apart}) and conditional bound for (\ref{spart}) to conditionally bound (\ref{mainuboundgg}).\label{step4halfoutline}
\item Prove that one can indeed bound $\E[\big|n t - \sum_{i=1}^{n} N_{e,i}(t)|^r\big]$ as needed in the conditional bounds (by making completely explicit, and enhancing, an approach to bounding centered renewal processes sketched in \citet{Gut09}). \label{step5outline}
\item Combine all of the above to yield the desired result . \label{step6outline}
\end{enumerate}
\subsubsection{High-level outline of proofs of other results.}  Our bounds for the s.s.p.d. appearing in Theorems\ \ref{mastertheoremsmallr} - \ref{mastertheoremlarger}, as well as our bound for the number of busy servers (Theorem\ \ref{numbusybbb} ), follow from a very similar logic, and those proofs can be found in Section\ \ref{sspdsec}.  We structure our proofs so the arguments used in proving our bounds for $\pr\big( L(\infty) \geq \frac{x}{1-\rho} \big)$ in Theorems\ \ref{mastertheoremsmallr} - \ref{mastertheoremlarger} can be easily ported over to these other settings.  Our bounds under additional assumptions on the moments of $S$ (Theorem\ \ref{tailboundedmoments}) follows by optimizing our bounds (by solving for the ``best $r$" for each $x$) from Theorem\ \ref{mastertheoremlarger}, and we defer the proof to the supplemental appendix Section\ \ref{optrsec}.  Our results from Section\ \ref{noprefacsec} (i.e. Theorems\ \ref{sspdnoprefactor} - \ref{sspdnoprefactorc}), whose proofs appear in the supplemental appendix Section\ \ref{noprefactorsec0}, are derived using very different (and simpler) drift arguments.

\subsection{Proof of our main results : bounds for $L(\infty)$ in Theorems\ \ref{mastertheoremsmallr} - \ref{mastertheoremlarger} .}\label{abitmoredetailsec}
We now prove our main results, the bounds for the tail of $L(\infty)$ in Theorems\ \ref{mastertheoremsmallr} - \ref{mastertheoremlarger}, by implementing the proof outlined in Section\ \ref{highleveloutlinesec} above.  
\subsubsection{\ref{step1outline}. : Use stochastic comparison results of \citet{GG13} to bound $\pr\big( L(\infty) \geq x \big)$ by (\ref{mainuboundgg}).}
In \citet{GG13}, the authors proved that $L(\infty)$ is stochastically dominated by the supremum of a certain one-dimensional random walk.  This random walk arises from analyzing a modified queueing system in which an artificial arrival is added to the system whenever a server would otherwise go idle.  To simplify notation the authors of \citet{GG13} imposed the restriction that $\pr(A = 0) = \pr(S = 0) = 0$ (to preclude having to deal with simultaneous events).  However, this restriction is unnecessary and the proofs of \citet{GG13} can be trivially modified to accomodate this setting.  As such, we state the relevant stochastic-comparison result of \citet{GG13} without that unnecessary assumption. 

\begin{lemma}[\citet{GG13}]\label{ggoldb}
Suppose that $\mu_A < n \mu_S$, and that $Q(\infty)$ exists.  Then for all $x > 0$,
$$\pr\big(  L(\infty) \geq x \big) \leq \pr\Bigg( \sup_{t \geq 0} \bigg( A_e(t) - \sum_{i=1}^{n} N_{e,i}(t) \bigg) \geq x \Bigg).
$$
\end{lemma}

\subsubsection{\ref{step2outline}. : Use a union bound, and the connection between ${\mathcal A}_e$ and ${\mathcal A}_o$, to bound (\ref{mainuboundgg}) by the sum of (\ref{apart}) and (\ref{spart}).} Note that we may construct ${\mathcal A}_e, {\mathcal A}_{o}$ on the same probability space s.t. w.p.1, 
\begin{equation}\label{chaaaanges}
A_e(t) \leq 1 + A_o(t)\ \textrm{for all}\ t \geq 0.
\end{equation}
The above inequality follows by observing that the set of event times in ${\mathcal A}_e$, after the first event, is an ordinary renewal process.  Next, we apply a straightforward union bound to reduce the problem of bounding (\ref{mainuboundgg}) to that of bounding (\ref{apart}) and (\ref{spart}).  We defer a formal proof to the supplemental appendix Section\ \ref{2partboundsec}.  
\begin{lemma}\label{2partbound}
Suppose that $\E[S] = 1$ and $\mu_A < n$.  Then for all $x > 2$,
\begin{eqnarray*}
\ &\ &\ \pr\bigg( \sup_{t \geq 0} \big(A_e(t) - \sum_{i=1}^{n} N_{e,i}(t) \big) \geq x \bigg)
\\&\leq&\ \pr\Bigg( \sup_{t \geq 0} \bigg( A_o(t) - \mu_{A} t - \frac{1}{2}( n - \mu_{A} ) t \bigg) \geq \frac{1}{2} x - 1 \Bigg) 
\\&\ &\ \ \ +\ \ \ \pr\Bigg( \sup_{t \geq 0} \bigg( n t - \sum_{i=1}^{n} N_{e,i}(t) - \frac{1}{2}( n - \mu_{A} ) t  \bigg) \geq \frac{1}{2} x \Bigg). 
\end{eqnarray*}
\end{lemma}

\subsubsection{ \ref{step3outline}. : Bound (\ref{apart}) by relating the supremum to a simple single-server queue.}
We now bound  (\ref{apart}).  Here we bound the corresponding supremum for general positive linear drift $\nu$, but will later plug in $\frac{1}{2}( n - \mu_A )$.  In particular, we prove the following.

\begin{lemma}\label{arrsuplem2}
Suppose that $\E[A^2] < \infty$.  Then for all $\nu > 0$ and $x > 0$,
$$\pr\bigg( \sup_{t \geq 0} \big( A_o(t) - \mu_A t - \nu t \big) \geq x \bigg) \leq 
\exp\bigg( - .09 \frac{\nu}{\nu + \mu_{A}} \big( \E[ (A \mu_A)^2 ] \big)^{-1} x \bigg).$$
\end{lemma}

Our proof proceeds by relating the supremum to the steady-state waiting time in a certain single-server queue, and then applying a result of \citet{K64} bounding the relevant tail probabilities.  We defer the proof to the technical appendix Section\ \ref{arrivaldetailsec}.

\subsubsection{\ref{step4outline}. : Conditionally bound (\ref{spart}).}
In this section, we prove that IF one could suitably bound $\E[\big|n t - \sum_{i=1}^{n} N_{e,i}(t)|^r\big]$ (for some $r > 2$ and all $t$), THEN one could bound (\ref{spart}) as required (using known maximal inequalities)
We defer all proofs to the technical appendix Section\ \ref{spartlatersec}, and instead simply state the most relevant result.  We restrict to the setting $\E[S] = 1$, which suffices since one can derive the general case by simply rescaling time (i.e. multiplying both the service and inter-arrival times by $\mu_S$).  This follows from the fact that such a rescaling does not change the distribution of $Q(\infty)$, and only impacts the proven bounds by replacing terms of the form $\E[S^k]$ by $\E[(S \mu_S)^k]$, leaving all other quantities unchanged (as $\E[(A \mu_A)^k]$ and $\rho$ are unchanged by such a rescaling).

\begin{lemma}\label{alltimepooled1}
Suppose that $\E[S] = 1$, and that for some fixed integer $n \geq 1$ and constants $C_1, C_2 > 0; r_1 > s > 1$; and $r_2 > 2$:
				\begin{enumerate}[(i)]
					\item For all $t \geq 1$,
$$\E\big[|\sum_{i=1}^{n} N_{e,i}(t) - n t|^{r_1}\big] \leq C_1 n^{\frac{r_1}{2}} t^{s}.$$
					\item For all $t \in [\frac{2}{n},1]$,
$$\E\big[|\sum_{i=1}^{n} N_{e,i}(t) - n t|^{r_2}\big] \leq C_2  (n t)^{\frac{r_2}{2}}.$$
\end{enumerate}
Then for all $\nu > 0$ and $x \geq 8$,
$$\pr\Bigg( \sup_{t \geq 0} \bigg( n t - \sum_{i=1}^{n} N_{e,i}(t) - \nu t  \bigg) \geq x \Bigg)$$ 
is at most
$$
3.6 \times (1 + \frac{1}{r_1 - s}) \times \bigg( (16 \frac{r_1+1}{s-1})^{r_1 + 1} \times C_1 n^{\frac{r_1}{2}} \nu^{-s} x^{- (r_1 - s)}
+
 (23 \frac{r_2+1}{\frac{r_2}{2}-1})^{r_2 + 1} \times  C_2 n^{\frac{r_2}{2}} (x \nu)^{- \frac{r_2}{2}} \bigg).
$$
\end{lemma}

\subsubsection{\ref{step4halfoutline}. : Conditionally bound (\ref{mainuboundgg}).}
In this section, we prove that IF one could suitably bound $\E[\big|n t - \sum_{i=1}^{n} N_{e,i}(t)|^r\big]$ (for some $r > 2$ and all $t$), THEN one could bound (\ref{mainuboundgg}) by combining our previous bound for (\ref{apart}) with our previous conditional bound for (\ref{spart}), as we have already bounded (\ref{mainuboundgg}) by the sum of (\ref{apart}) and (\ref{spart}).  The proof follows in a straightforward manner by using Lemma\ \ref{alltimepooled1} to bound (\ref{spart}), and Lemma\ \ref{arrsuplem2} to bound (\ref{apart}), combined with Lemma\ \ref{2partbound} and some straightforward algebra, and we omit the details.  $\Halmos$

\begin{theorem}\label{mastertheorem}
Suppose that $\E[S] = 1$, and that for some fixed integer $n \geq 1$ and constants $C_1, C_2 > 0; r_1 > s > 1$; and $r_2 > 2$, the following conditions hold:
				\begin{enumerate}[(i)]
					\item $\mu_A < n$.
					\item For all $t \geq 1$,
$$\E\big[|\sum_{i=1}^{n} N_{e,i}(t) - n t|^{r_1}\big] \leq C_1 n^{\frac{r_1}{2}} t^{s}.$$
					\item For all $t \in [\frac{2}{n},1]$,
$$\E\big[|\sum_{i=1}^{n} N_{e,i}(t) - n t|^{r_2}\big] \leq C_2  (n t)^{\frac{r_2}{2}}.$$
\end{enumerate}
Then for all $x \geq 18$, $\pr\bigg( \sup_{t \geq 0} \bigg( A_e(t) - \sum_{i=1}^{n} N_{e,i}(t) \bigg) \geq x \bigg)$ is at most
\begin{eqnarray*}
&\ &\ 1.8 \times (1 + \frac{1}{r_1 - s}) \times (32 \frac{r_1+1}{s-1})^{r_1 + 1} \times C_1 n^{\frac{r_1}{2}} (n - \mu_A)^{-s} x^{- (r_1 - s)} 
\\&\ \ \ &\ \ +\ \ 1.8 \times (1 + \frac{1}{r_1 - s}) \times (46 \frac{r_2+1}{\frac{r_2}{2}-1})^{r_2 + 1} \times C_2 n^{\frac{r_2}{2}} (n - \mu_A)^{- \frac{r_2}{2}} x^{- \frac{r_2}{2}}
\\&\ \ \ &\ \ +\ \ 1.1 \times \exp\bigg( - .045 \frac{n - \mu_A}{n + \mu_A} \big( \E[ (A \mu_A)^2 ] \big)^{-1} x \bigg).
\end{eqnarray*}
\end{theorem}

\subsubsection{\ref{step5outline}. : Prove that one can indeed bound $\E[\big|n t - \sum_{i=1}^{n} N_{e,i}(t)|^r\big]$ as needed to apply the conditional bounds of Lemma\ \ref{alltimepooled1} and Theorem\ \ref{mastertheorem}.}
We now show that one can indeed appropriately bound the central moments of $\sum_{i=1}^{n} N_{e,i}(t)$ s.t. the conditional bounds of Lemma\ \ref{alltimepooled1} and Theorem\ \ref{mastertheorem} can be applied.  We state two results, one for $t \geq 1$, and one for $t \in [\frac{2}{n},1]$.  In both cases, we defer the proofs to the supplemental appendix Sections\ \ref{poolerbound1sec} and \ref{binomial2sec}.

\begin{lemma}\label{poolerbound1}
Suppose that $\E[S] = 1$, and $\E[S^r] < \infty$ for some $r \geq 2$.  Then for all $n \geq 1$ and $t \geq 1$, $\E\big[ \big|\sum_{i=1}^n N_{e,i}(t) - n t \big|^r\big]$ is at most
$$\bigg( .76 \times \big( \E[S^2] + 1 \big)^r \times 1032^r \times r^{2 r} + .21 \times (\E[S^2] + 1) \times 516^r \times r^{1.5 r} \times (\E[S^r] + 1) \bigg) \times (n t)^{\frac{r}{2}}.$$
\end{lemma} 

We defer the proofs to the supplemental appendix Section\ \ref{poolerbound1sec}.

\begin{lemma}\label{binomial2}
Suppose that $\E[S] = 1$, and $\E[S^2] < \infty$.  Then for all $n \geq 1, r \geq 2, t \in [\frac{2}{n},1],$ 
\begin{equation}\label{focus1}
\E\bigg[ \big| \sum_{i=1}^n N_{e,i}(t) - n t \big|^r \bigg] \leq 
5.2 \times \big( 35 (1 + \E[S^2]) \big)^r \times r^{2.5 r} \times (nt)^{\frac{r}{2}}.
\end{equation}
\end{lemma}

We defer the proofs to the supplemental appendix Section\ \ref{binomial2sec}.

\subsubsection{\ref{step6outline}. : Proof of our main results, bounds for $L(\infty)$ in Theorems\ \ref{mastertheoremsmallr} - \ref{mastertheoremlarger}.}\label{togetherproofsubsec}
In this section, we prove our main results, the first part of Theorems\ \ref{mastertheoremsmallr} - \ref{mastertheoremlarger} in which $\pr\big( L(\infty) \geq \frac{x}{1-\rho} \big)$ is bounded.  We proceed by plugging in the bounds for $\E\bigg[ \big| \sum_{i=1}^n N_{e,i}(t) - n t \big|^r \bigg]$ shown in Lemmas\ \ref{poolerbound1} and \ref{binomial2} into the conditional bounds of Theorem\ \ref{mastertheorem}.
\proof{Proof of the bounds for $L(\infty)$ in Theorems\ \ref{mastertheoremsmallr} - \ref{mastertheoremlarger} : }
First, we again note that it suffices to prove the result for the case $\E[S] = 1$, by a simple rescaling argument (in which $S$ is replaced by $S \mu_S$).  Thus suppose $\E[S] = 1$.  Then combining Lemmas\ \ref{poolerbound1} and \ref{binomial2} with some straightforward algebra, we conclude the following.  For each integer $n \geq 1$ s.t. $n > \mu_A$, the conditions of Theorem\ \ref{mastertheorem} are met with the following parameters:
$$r_1 = r_2 = r\ \ \ ,\ \ \ s = \frac{r}{2},$$
$$C_1 =  .76 \times \big( \E[S^2] + 1 \big)^r \times 1032^r \times r^{2 r} + .21 \times (\E[S^2] + 1) \times 516^r \times r^{1.5 r} \times (\E[S^r] + 1),$$
$$C_2 =  5.2 \times \big( 35 (1 + \E[S^2]) \big)^r \times r^{2.5 r}.$$
Thus, applying Theorem\ \ref{mastertheorem} and some straightforward algebra (also using the fact that $\frac{n}{n-\mu_A} = \frac{1}{1-\rho}$ and $\frac{n - \mu_A}{n + \mu_A} = \frac{n}{n + \mu_A} (1 - \rho) \geq \frac{1}{2} (1 - \rho)$) , we find that for all $x \geq 18$ and $r \in (2,r^*)$, (\ref{mainuboundgg}) is at most
\begin{eqnarray*}
&\ &\ 906 \times (\frac{r+1}{\frac{r}{2}-1})^{r+1} \times \big( \E[S^2] + 1 \big) \times 33024^r \times \bigg( \big( \E[S^2] + 1 \big)^{r-1} \times r^{2.5 r} + r^{1.5 r} \times \E[S^r] \bigg) \times \big( x (1-\rho) \big)^{-\frac{r}{2}}
\\&\ &\ \ \ + 1.1 \times \exp\bigg( - .0225 \big( \E[ (A \mu_A)^2 ] \big)^{-1} (1-\rho) x \bigg).
\end{eqnarray*}
We now break into 2 cases, $r \leq 2.5$ and $r > 2.5$.  For $r \leq 2.5$, we find (after some straightforward algebra) that for $x \geq 18$ and $r \in \big(2, \min(2.5,r^*) \big)$, (\ref{mainuboundgg}) is at most
\begin{eqnarray*}
\ &\ &\ 5 \times 10^{18} \times \big( \E[S^2] + 1 \big) \times \bigg( \big( \E[S^2] + 1 \big)^{r-1} + \E[S^r] \bigg) \times 
(\frac{r}{2}-1)^{-(r+1)} \times \big( (1-\rho) x \big)^{-\frac{r}{2}} \nonumber
\\&\ \ \ &\ \ +\ \ 1.1 \times \exp\bigg( - .0225 \big( \E[ (A \mu_A)^2 ] \big)^{-1} (1 - \rho) x \bigg); 
\end{eqnarray*}
while for $r \in (2.5,r^*)$, using the easily verified fact that $\frac{r+1}{\frac{r}{2}-1} \leq 14$ for $r \geq 2.5$, 
we find that for $x \geq 18$ and $r \in (2.5,r^*)$, (\ref{mainuboundgg}) is at most
\begin{eqnarray*}
&\ &\ 1.3 \times 10^4 \times (\E[S^2] + 1) \times \big(4.6 \times 10^5\big)^r \times \big( (\E[S^2] + 1)^{r-1} \times r^{2.5 r} + r^{1.5 r} \times \E[S^r] \big) \times \big( (1-\rho) x \big)^{-\frac{r}{2}}
\\&\ \ \ &\ \ +\ \ 1.1 \times \exp\bigg( - .0225 \big( \E[ (A \mu_A)^2 ] \big)^{-1} (1-\rho) x \bigg).
\end{eqnarray*}
Noting that the bounds are anyways at least one for $x \in (0,18)$, and combining with the fact that $\E[S^2], \E[S^r] \geq 1$ and some straightforward algebra, completes the proof of our main results, the bounds for $\pr\big( L(\infty) \geq \frac{x}{1-\rho} \big)$ appearing in Theorems\ \ref{mastertheoremsmallr} - \ref{mastertheoremlarger}.  $\halmos$ \endproof

\section{Proofs of the bounds for the s.s.p.d. in Theorems\ \ref{mastertheoremsmallr} - \ref{mastertheoremlarger} and number of busy servers in Theorem\ \ref{numbusybbb}.}\label{sspdsec}
In this section we complete the bounds for the s.s.p.d. in the proofs of Theorems\ \ref{mastertheoremsmallr} - \ref{mastertheoremlarger}, and for the number of busy servers in Theorem\ \ref{numbusybbb}.  Note that the stochastic comparison result Lemma\ \ref{ggoldb}, upon which our entire analysis is premised, can only provide trivial bounds for the s.s.p.d. (i.e. plugging in $x = 0$ yields a trivial bound of 1), and similarly cannot be used to bound the number of busy servers.  Here we show that a certain enhancement to the bounds of \citet{GG13} can overcome this problem.  The enhancement is based on the intuition that the probability that an $n$-server queueing system has at least $n$ jobs in system is at most the probability that an $n'$-server queueing system has at least $n$ jobs in system for $n' < n$.  But the bounds of Lemma\ \ref{ggoldb} DO yield meaningful bounds for the latter quantity, as it is equivalent to the probability that an $n'$-server system has at least $n - n'$ jobs waiting in queue.  We defer the proof to the technical appendix Section\ \ref{ggoldbsec}.

\begin{lemma}\label{ggoldb222}
Suppose that for the $GI/GI/n$ queue with inter-arrival times distributed as the r.v. $A$ and service times distributed as the r.v. $S$, it holds that $\mu_A < n \mu_S$, and $Q(\infty)$ exists.  Then for all $n' \in \lbrace 1,\ldots,n \rbrace$ and $x \geq n' - n,$  

$$\pr\big(  Q(\infty) - n \geq x \big) \leq \pr\Bigg( \sup_{t \geq 0} \bigg( A_e(t) - \sum_{i=1}^{n'} N_{e,i}(t) \bigg) \geq x + (n - n') \Bigg).
$$
\end{lemma}

\subsection{Proofs of the bounds for the s.s.p.d. in Theorems\ \ref{mastertheoremsmallr} - \ref{mastertheoremlarger}.}
In this section, we again suppose without loss of generality (by rescaling) that $\E[S] = 1$.  As a direct corollary of Lemma\ \ref{ggoldb222} (by plugging in $x = 0, n' = n - \lfloor \frac{1}{2} (n - \mu_A) \rfloor$), we conclude the following bound for the s.s.p.d.

\begin{corollary}\label{ggold2}
Under the same assumptions as Lemma\ \ref{ggoldb222}, and supposing also $\E[S] = 1$, it holds that 
$$
\pr\big( Q(\infty) \geq n \big) \leq \pr\Bigg( \sup_{t \geq 0} \bigg( A_e(t) - \sum_{i=1}^{n - \lfloor \frac{1}{2} (n - \mu_A) \rfloor} N_{e,i}(t) \bigg) \geq \lfloor \frac{1}{2} (n - \mu_A) \rfloor \Bigg).
$$
\end{corollary}

Note that Corollary\ \ref{ggold2} reduces bounding the s.s.p.d. to bounding $\pr\Bigg( \sup_{t \geq 0} \bigg( A_e(t) - \sum_{i=1}^{n'} N_{e,i}(t) \bigg) \geq x' \Bigg)$ for $x' = \lfloor \frac{1}{2} (n - \mu_A) \rfloor$ and $n' = n - \lfloor \frac{1}{2} (n - \mu_A) \rfloor$ some integer different from the number of servers $n$ in the original system.  None-the-less, we can still apply those parts of Theorems\ \ref{mastertheoremsmallr} - \ref{mastertheoremlarger} which we have already proven, i.e. the tail bounds for the queue length, with these different parameters to complete the proofs of our bounds for the s.s.p.d.  Here we implicitly use the fact that the relevant bounds from Theorems\ \ref{mastertheoremsmallr} - \ref{mastertheoremlarger} are actually bounds for $\sup_{t \geq 0} \bigg( A_e(t) - \sum_{i=1}^n N_{e,i}(t) \bigg)$.
\proof{Proof of bounds for s.s.p.d. in Theorems\ \ref{mastertheoremsmallr} - \ref{mastertheoremlarger}.}
First, suppose $\rho \leq 1 - \frac{4}{n}$, which implies that
$$\frac{1}{4}(n - \mu_A) \leq \lfloor \frac{1}{2}(n - \mu_A) \rfloor \leq \frac{1}{2}(n - \mu_A).$$ 
Let $x' = \lfloor \frac{1}{2} (n - \mu_A) \rfloor$ and
$n' = n - \lfloor \frac{1}{2} (n - \mu_A) \rfloor.$  Then it follows from some straightforward algebra that $x' (1 - \rho_{n'})$ is at least
\begin{eqnarray*}
\ &\ &\ \frac{1}{4}(n - \mu_A) \times \Bigg(1 - \frac{\mu_A}{\bigg( n - \Big( \frac{1}{2}(n - \mu_A) \Big) \bigg) }\Bigg)
\\&\ &\ \ \ \ \ \ =\ \ \ \frac{1}{4} \frac{ (n - \mu_A)^2 }{n + \mu_A}\ \ \ \ \geq\ \ \ \frac{1}{8} \frac{ (n - \mu_A)^2 }{n}\ \ \ =\ \ \ \frac{n}{8} (1 - \rho)^2.
\end{eqnarray*}
Combining with some straightforward algebra, Corollary\ \ref{ggold2}, and Theorems\ \ref{mastertheoremsmallr} - \ref{mastertheoremlarger} proves the desired result for all $n \geq 5$ s.t. $\rho \leq 1 - \frac{4}{n}$.  Noting that if either $n \leq 4$, or $n \geq 5$ and $\rho > 1 - \frac{4}{n},$ then all relevant bounds are at least one (and hence hold for the s.s.p.d.) completes the proof of the bounds for the s.s.p.d. of Theorems\ \ref{mastertheoremsmallr} - \ref{mastertheoremlarger}.$\halmos$ \endproof
\subsection{Proof of bound for the number of busy servers in Theorem\ \ref{numbusybbb}.}
Again suppose $\E[S] =1$.  Similar to our logic in bounding the s.s.p.d, let $x' = \mu_A - n + x \sqrt{\mu_A} , n' = \lceil \mu_A  + \frac{x}{2} \sqrt{\mu_A} \rceil.$  Noting that $x' \geq n' - n$ by construction (as $x \sqrt{\mu_A} > \frac{x}{2} \sqrt{\mu_A}$), and supposing that $n' \leq n$ (which we will enforce by requiring $x \in [0, 2 \frac{n - \mu_A - 1}{\sqrt{\mu_A}}]$), we may thus apply Lemma\ \ref{ggoldb222} with $n, x' = \mu_A - n + x \sqrt{\mu_A} , n' = \lceil \mu_A  + \frac{x}{2} \sqrt{\mu_A} \rceil.$  As $Q(\infty) - n \geq x' \leftrightarrow Q(\infty) \geq \mu_A + x \sqrt{\mu_A},$ and $x' + (n - n') \geq  \frac{x}{2} \sqrt{\mu_A} - 1$, we conclude the following bound for the tail of the number of busy servers (by applying Lemma\ \ref{ggoldb222}).

\begin{corollary}\label{ggoldlightb}
Under the same assumptions as Lemma\ \ref{ggoldb222}, for all $x \in [0, 2 \frac{n - \mu_A - 1}{\sqrt{\mu_A}}]$, it holds that 
$$\pr\bigg( Q(\infty) \geq \mu_A + x \sqrt{\mu_A} \bigg) \leq \pr\Bigg( \sup_{t \geq 0} \bigg( A_e(t) - \sum_{i=1}^{\lceil \mu_A + \frac{x}{2} \sqrt{\mu_A} \rceil} N_{e,i}(t) \bigg) \geq \frac{x}{2} \sqrt{\mu_A} - 1 \Bigg).
$$
\end{corollary}

Note that Corollary\ \ref{ggoldlightb} reduces bounding the number of busy servers to bounding $\pr\Bigg( \sup_{t \geq 0} \bigg( A_e(t) - \sum_{i=1}^{n''} N_{e,i}(t) \bigg) \geq x'' \Bigg)$ for $x'' = \frac{x}{2} \sqrt{\mu_A} - 1$ and $n'' = \lceil \mu_A + \frac{x}{2} \sqrt{\mu_A} \rceil.$ As in the proof of our bounds for the s.s.p.d., we can still apply the tail bounds for the queue length of Theorem\ \ref{mastertheoremlarger} with these different parameters to complete the proof of our bound for the number of busy servers.

\proof{Proof of Theorem\ \ref{numbusybbb}:}
Let $x'' = \frac{x}{2} \sqrt{\mu_A} - 1$ and $n'' = \lceil \mu_A + \frac{x}{2} \sqrt{\mu_A} \rceil$.  First, let us point out that requiring $\rho \in [\frac{1}{n},  1 - \frac{2}{n}]$ may be seen to imply two properties (after some straightforward algebra) : (1)\ $\frac{ 2 (n - \mu_A - 1)}{\sqrt{\mu_A}} \geq \frac{1}{2} \frac{n - \mu_A}{\sqrt{\mu_A}}$, and hence Corollary\ \ref{ggoldlightb} is applicable to all $x$ in the range $\bigg[ 4, \frac{1}{2} \min\big(\sqrt{\mu_A}, \frac{n - \mu_A}{ \sqrt{\mu_A}}\big) \bigg]$; and (2)\ $\mu_A \geq 1$.  Next, note that 
\begin{eqnarray*}
x'' (1 - \rho_{n''}) &=& \big(\frac{x}{2} \sqrt{\mu_A} - 1\big) \times \big(1 - \frac{\mu_A}{\lceil \mu_A + \frac{x}{2} \sqrt{\mu_A}\rceil}\big)
\\&\geq& \frac{x}{4} \sqrt{\mu_A} \times \big(1 - \frac{ \mu_A }{\mu_A + \frac{x}{2} \sqrt{\mu_A}}\big)\ \ \textrm{by our assumptions that}\ x \geq 4\ \textrm{and}\ \mu_A \geq 1
\\&=& \frac{x}{4} \sqrt{\mu_A} \times \frac{ \frac{x}{2} \sqrt{\mu_A}}{\mu_A + \frac{x}{2} \sqrt{\mu_A}}
\\&\geq& \frac{x}{4} \sqrt{\mu_A} \times \frac{ \frac{x}{2} \sqrt{\mu_A} }{\mu_A + \frac{1}{2} \sqrt{\mu_A} \times \sqrt{\mu_A}}\ \ \textrm{by our assumption that}\ x \leq \frac{1}{2} \sqrt{\mu_A}
\\&\geq& \frac{x^2}{16}.
\end{eqnarray*}
The desired result then follows by using the tail bounds for the queue length of Theorem\ \ref{mastertheoremlarger} with parameters $n'',x''$ to bound $\pr\Bigg( \sup_{t \geq 0} \bigg( A_e(t) - \sum_{i=1}^{n''} N_{e,i}(t) \bigg) \geq x'' \Bigg)$, along with some straightforward algebra.  $\Halmos$ \endproof

\section{Conclusion and future research directions.}\label{concsec}
In this paper, we proved the first simple and explicit bounds for $GI/GI/n$ queues which scale as $\frac{1}{1 - \rho}$ (analogous to Kingman's bound for single-server queues), assuming only that $\E[A^2] < \infty$ and $\E[S^{2 + \epsilon}] < \infty$ for some $\epsilon > 0$.  Our main results are bounds for the tail of the steady-state queue length and the steady-state probability of delay.  The strength of our bounds (e.g. in the form of tail decay rate) is a function of how many moments of the service distribution are assumed finite, and we obtain even stronger results under the assumption that all moments of $S$ exist and satisfy a mild growth rate assumption.  Our bounds scale gracefully even when the number of servers grows large and the traffic intensity converges to unity simultaneously, as in the Halfin-Whitt scaling regime.  Some of our bounds scale better than $\frac{1}{1-\rho}$ in certain asymptotic regimes, for which we also prove bounds for the tail of the steady-state number in service.  We also prove several additional bounds using drift arguments (which have much smaller pre-factors), and point out a conjecture which would imply further related bounds and generalizations.
\\\indent Our results leave many interesting directions for future research.
\begin{itemize}
\item Our demonstrated tail decay rates are suboptimal, and a bound which is uniformly accurate (in $n,\rho$, and $x$) across all scaling regimes remains elusive.  This is also true regarding our bounds for the s.s.p.d., which are similarly suboptimal.\footnote{We refer the reader to \citet{G16} for further progress on asymptotic bounds for the s.s.p.d.  In the previous version of \citet{G16} and the present manuscript, a version of Lemma\ \ref{ggoldb222} and Corollary\ \ref{ggold2} instead appeared in \citet{G16} (with the present manuscript modifying and using those results).  In this version of the present manuscript these results instead appear in the present manuscript, and will be cited and used as appropriate by a new version of \citet{G16}, making the present manuscript self-contained.  We note that in contrast to the present work, \citet{G16} is restricted to asymptotic results in the Halfin-Whitt regime (and furthermore does not study the tail or mean of the queue length).}
\item The necessity of $r^{\Omega(r)}$ prefactors in our bounds (as well as what numerical constants are actually required) remains an open question, as does the question of whether tighter bounds (with smaller prefactors, possibly of a qualitatively different nature) can be derived (possibly under additional assumptions on $A,S$ and/or using a different proof technique).
\item Kingman's bound for single-server queues requires only that $\E[S^2] < \infty$, while our bounds require $\E[S^{2 + \epsilon}] < \infty$ for some $\epsilon > 0$.  Closing this gap remains an interesting open question.
\item Conjecture\ \ref{conjecture1}, which would imply a very simple bound for the expected queue length in $M/GI/n$ queues, remains a very interesting open question.
\item We conjecture that our approach can be modified to yield general and explicit bounds with an appropriate analogue of $\frac{1}{1-\rho}$ scaling for a broad range of queuing networks.  Although we are not aware of any simple and explicit analogues of Kingman's bound for queueing networks conjectured in the literature, we note that past work on heavy-traffic in queueing networks suggests that the number in queue at each station $i$ should scale as $\frac{1}{1-\rho_i}$ with $\rho_i$ the effective traffic intensity at that station (as dictated by the so-called traffic equations, see e.g. \citet{Reiman84,MMR.98,GZ06,Dai20}).  We leave a formal investigation along these lines as an interesting direction for future research, but do provide the sketch of a plausible approach in the supplemental appendix Section\ \ref{appnetworksec}.  It is also interesting to ask for what other queueing systems our approach can be implemented.  For example, in the parallel work \citet{G18}, the authors extend the stochastic comparison approach of \citet{GG13} to certain multi-server systems with abandonments and hyper-exponentially distributed service times.  The authors have also extended this approach to heavy-tailed systems (in the Halfin-Whitt regime) in the parallel work \citet{G17}.\footnote{We note that the work \citet{G18} is restricted to the Halfin-Whitt scaling regime, and does not yield explicit non-asymptotic bounds.  We also note that the analysis of \citet{G17} for multi-server queues with heavy tails (in the Halfin-Whitt regime) relies heavily on the bounds in the present manuscript, combined with a novel analysis of heavy-tailed renewal processes.}  Understanding the complete set of systems to which our stochastic comparison approach (or a suitable modification thereof) can be applied to derive simple and explicit bounds remains an interesting direction for future research.
\end{itemize}
\ \\\textbf{Broader connections to the applied probability and operations research communities.}  Taking a broader view of the literature not just on queueing, but on applied probability and operations research more broadly, let us reflect on some of the high-level take-aways of this work.  
\begin{itemize}
\item Our work provides an example of how ``quantifying" a simple coupling argument yields explicit bounds for a fundamental model.  Such an approach has been sucessful in several operations research models recently, and we point the reader to \citet{Xin16, vera2021}.
\item Our work provides an example of a setting where a powerful inuition/scaling for a very simple ``base model" (here $\frac{1}{1-\rho}$ scaling for single-server queues) carries over to natural generalizations more relevant in practice (here multi-server queues).  
\item Our work contextualizes and compares many approaches taken to multi-server queues, including stochastic comparison and Lyapunov drift, also surveying the relevant results.  As these different approaches appear in the study of many stochastic models, lessons learned from the queueing setting can inform the study of other stochastic models.  Our work also provides a useful reference for the vast literature on multi-server queues.
\item Our work touches on meta-questions about how to conceptualize the trade-off between simplicity/explicitness, and accuracy, in approximations for operations research models.  
\end{itemize}
\ \\\textbf{Final thoughts on the trade-off between simplicity and accuracy in operations research models.}  Several pressing big-picture questions along these lines remain unresolved in the study of stochastic models broadly.
\begin{itemize}
\item What is the right notion of ``complexity" in approximations for such models?  
\item How should one compare analytical bounds with results derived from simulation and numerical procedures? 
\item What is the formal algorithmic complexity of both numerical computation, and simulation, for the limiting processes which arise?   
\item And last, but by no means least, which types of approximations may be most useful in practice?  
\end{itemize}

\section{Technical Appendix.}\label{moredetailsec}
\subsection{Bound (\ref{apart}) by relating to a GI/GI/1 queue, and proof of Lemma\ \ref{arrsuplem2}.}\label{arrivaldetailsec}
In this section we fill in the details in our proof of Lemma\ \ref{arrsuplem2}, which we used to bound the supremum associated with the arrival process, (\ref{apart}).
Our proof proceeds as follows.  First, we relate the desired supremum to a discrete-time supremum associated with $k -  \sum_{i=1}^k (\mu_A A_i)$.  Second, we observe that this supremum is the steady-state waiting time in a certain single-server queue, and then apply a result of \citet{K64} bounding the relevant tail probabilities.  We also prove a novel result showing that a certain exponent appearing in Kingman's results can be bounded only in terms of the first two moments of $A$ and the associated drift parameter, which may be of independent interest. 
\\\\We begin with the following lemma relating the supremum appearing in (\ref{apart}) to the steady-state waiting time in a certain single-server queue, whose proof we defer to the supplemental appendix Section\ \ref{adiscboundsec}.  We note that the result follows in a straightforward manner by applying the basic definitions of renewal processes and some standard transformations.
\begin{lemma}\label{adiscbound}
For all $\nu > 0$ and $x > 0$, 
\begin{equation}\label{tobounda1}
\pr\bigg( \sup_{t \geq 0} \big( A_o(t) - \mu_A t - \nu t \big) \geq x \bigg)
\end{equation}
equals
\begin{equation}\label{apart22}
\pr\Bigg( \sup_{k \geq 0} \bigg(  \frac{\mu_A}{\mu_A + \nu} k - \sum_{i=1}^k (\mu_A A_i) \bigg) \geq x (1 + \frac{\nu}{\mu_{A}})^{-1} \Bigg).
\end{equation}
\end{lemma}
We note that the supremum term appearing in (\ref{apart22}) is exactly the steady-state waiting time in a single-server queue with inter-arrival times i.i.d. distributed as $A \mu_A$, and service times the constant $\frac{\mu_A}{\mu_{A} + \nu}.$  
\\\\Next, we recall the relevant result of Kingman for bounding the tails of such suprema.  We state the result in a specific form that we will need it, but note that the results of \citet{K64} hold in a more general setting.  We note that the results of \citet{K64} follow from standard martingale techniques, i.e. looking at an appropriate exponential martingale and applying a maximal inequality for martingales.  We state the relevant result in terms of a more general supremum $\sup_{k \geq 0} (c \times k - Z_k)$ for $c \in (0,1)$ and $\lbrace Z_k, k \geq 1 \rbrace$ an i.i.d. sequence of non-negative mean one r.v.s.  

\begin{lemma}[\citet{K64}] \label{kingmanmart1}
Suppose $\lbrace Z_i, i \geq 1 \rbrace$ is an i.i.d. sequence of non-negative r.v.s with $\E[Z_1] = 1$, and $c \in (0,1)$ is some constant.  Suppose $\E\big[ \exp\big( \theta(c - Z_1) \big) \big] \leq 1$ for some $\theta > 0$.  Then $P\big( \sup_{k \geq 0} (c \times k - Z_k) \geq x \big) \leq \exp(- \theta x)$ for all $x > 0$.
\end{lemma}

Next, we show that one can explicitly characterize a $\theta$ s.t. $\E\big[ \exp\big( \theta(c - Z_1) \big) \big] \leq 1$ in terms of $c$ and the first two moments of $Z_1$, using a Taylor expansion.  To our knowledge the result is novel, and may be of independent interest in the analysis of single-server queues.  We defer the proof to the supplemental appendix Section\ \ref{kingmanmart2sec}.
\begin{lemma}\label{kingmanmart2}
Under the same assumptions as Lemma\ \ref{kingmanmart1}, and supposing in addition that $\E[Z^2_1] < \infty$, it holds that $\E\big[ \exp\big( \frac{ 1 - c }{11 c \E[Z_1^2]} (c - Z_1) \big) \big] \leq 1$.  Namely, one can take $\theta = \frac{ 1 - c }{11 c \E[Z_1^2]}$ in Lemma\ \ref{kingmanmart1}.
\end{lemma}

Note that $\frac{ 1 - c }{11 c \E[Z_1^2]}$ scales as $\Omega(1-c)$ as $c \uparrow 1$, and as $\Omega(\frac{1}{c})$ as $c \downarrow 0$.  With Lemmas\ \ref{kingmanmart1} and \ref{kingmanmart2} in hand, we now complete the proof of the desired result Lemma\ \ref{arrsuplem2}.
\proof{Proof of Lemma\ \ref{arrsuplem2}:}
The result follows by applying Lemmas\ \ref{kingmanmart1} and \ \ref{kingmanmart2} with $c = \frac{\mu_A}{\mu_{A} + \nu}, \lbrace Z_i, i \geq 1 \rbrace = \lbrace \mu_{A} A_{i}, i \geq 1 \rbrace$ to Lemma\ \ref{adiscbound}, after some straightforward algebra.  $\Halmos$
\endproof

\subsection{Conditionally bound (\ref{spart}) and proof of Lemma\ \ref{alltimepooled1}.}\label{spartlatersec}
In this section we fill in the details in our proof of Lemma\ \ref{alltimepooled1}, which provides conditional bounds on the supremum associated with the pooled renewal process, (\ref{spart}).  We proceed as follows.
\begin{itemize}
\item First,we prove a conditional result, which asserts that if the supremum of a centered continuous-time stochastic process can be controlled over : (1) sets of consecutive integers; and (2) intervals of size at most 1, then one can bound the tail of the all-time supremum of the same centered stochastic process with any negative linear drift.  We will ultimately apply this to 
the centered process $nt - \sum_{i=1}^n N_{e,i}(t)$ with drift $-\frac{1}{2}(n - \mu_A)$.
\item Second, we prove a conditional result that if the moments of $|nt - \sum_{i=1}^n N_{e,i}(t)|$ can be bounded as in the conditions of Lemma\ \ref{alltimepooled1}, then the supremum of $nt - \sum_{i=1}^n N_{e,i}(t)$ over both sets of consecutive integers and intervals of size at most one can indeed by suitably controlled.
\item Third, we combine the above to complete the proof of Lemma\ \ref{alltimepooled1}.
\end{itemize}
\subsubsection{Proof that controlling the supremum of a centered stochastic process over sets of consecutive integers, and short intervals, implies control of its all-time supremum with linear drift.}\label{controlmoredetailsec}
In this section we prove a conditional result, which converts bounds for the supremum of a suitable stochastic process over sets of consecutive integers, and intervals of length at most one, to bounds for the general all-time supremum (with negative drift).  We will ultimately use this result to bound (\ref{spart}), the supremum associated with $n t - \sum_{i=1}^{n} N_{e,i}(t).$  We note that similar arguments have been used to bound all-time suprema of stochastic processes (\citet{Sz.99}), also in the heavy-tailed setting (\citet{Sz.04}).  We include a self-contained exposition and proof in the supplemental appendix Section\ \ref{TailBoundsec}.

\begin{lemma} \label{TailBound}
			Let $\lbrace \phi(t) , t \geq 0 \rbrace$ be a stochastic process with stationary increments such that $\phi(0) = 0$.  Here, stationary increments means that for all $s_0 \geq 0$, $\lbrace \phi(s + s_0) - \phi(s_0), s \geq 0 \rbrace$ has the same distribution (on the process level) as $\lbrace \phi(s), s \geq 0 \rbrace$.  Suppose there exist strictly positive finite constants $H_1,H_2,s,r_1,r_2$ and $Z \geq 0$ such that $r_1 > s > 1$ and $r_2 > 2$, and the following two conditions hold:
				\begin{enumerate}[(i)]
					\item For all integers $m\ge 1$ and real numbers $x \geq Z$, 
					$$
						\pr (\max_{j\in \{0, ..., m\}} \phi(j) \geq x) \le H_1 m^s x^{-r_1}.
					$$
					\item For all $t_0 \in (0,1]$ and $x \geq Z$, 
					$$
						\pr (\sup_{0\le t\le t_0}\phi(t) \geq x) \le H_2 t_0^{\frac{r_2}{2}} x^{-r_2}.
					$$
				\end{enumerate}
				Then for any drift parameter $\nu > 0$, and all $x \geq 4 Z$, 
$$
\pr\bigg(\sup_{t\ge 0} (\phi(t) - \nu t) \geq x\bigg) \leq 12 (1 + \frac{1}{r_1 - s})\big( H_1 4^{r_1} x^{-(r_1-s)} \nu^{-s} + H_2 4^{r_2} (x \nu)^{ - \frac{r_2}{2} } \big).
$$
\end{lemma}

\subsubsection{Proof that if the central moments of $nt - \sum_{i=1}^n N_{e,i}(t)$ can be suitably bounded then the conditions of Lemma\ \ref{TailBound} hold.}

We now prove that if the central moments of $nt - \sum_{i=1}^n N_{e,i}(t)$ can be suitably bounded then the conditions of Lemma\ \ref{TailBound} hold, i.e. one can control the supremum of the centered pooled renewal process s.t. one can plug in $n t - \sum_{i=1}^n N_{e,i}(t)$ for $\phi(t)$ in Lemma\ \ref{TailBound}.  Our proof proceeds as follows.
\\\\First, we prove a relevant conditional bound for the supremum of of centered pooled renewal process over consecutive integers, as required by the first condition of Lemma\ \ref{TailBound}.  We defer the proof to the supplemental appendix Section\ \ref{maximalssec}.
\begin{lemma}\label{maximal2}
Suppose that $\E[S] = 1$, and that for some fixed $n \geq 1, C_1 > 0, s > 1$, and $r_1 \geq s$, the following condition holds:
\begin{enumerate}[(i)]
\item For all $t \geq 1$,
$$\E\big[|\sum_{i=1}^{n} N_{e,i}(t) - n t|^{r_1}\big] \leq C_1 n^{\frac{r_1}{2}} t^{s}.$$
\end{enumerate}
Then it also holds that for all non-negative integers $k$ and $x > 0$, 
$$
\pr\big(\max_{j \in \lbrace 1,\ldots, k\rbrace} \big|n j - \sum_{i=1}^{n} N_{e,i}(j)\big| \geq x\big) \leq 2.25 \times 2^{s} \times (2 \frac{r_1 + 1}{s - 1})^{r_1 + 1} \times C_1 n^{\frac{r_1}{2}} k^{s} x^{-r_1}.
$$
\end{lemma}
\ \\\\Second, we prove a relevant conditional bound for the supremum of of centered pooled renewal process over small intervals, as required by the second condition of Lemma\ \ref{TailBound}.  We defer the proof to the supplemental appendix Section\ \ref{maximalssec}.
\begin{lemma}\label{maximal3}
Suppose that $\E[S] = 1$, and that for some fixed $n \geq 1, C_2 > 0$ and $r_2 > 2$, the following condition holds:
\begin{enumerate}[(i)]
\item For all $t \in [\frac{2}{n},1]$,
$$\E\big[|\sum_{i=1}^{n} N_{e,i}(t) - n t|^{r_2}\big] \leq C_2 (n t)^{\frac{r_2}{2}}.$$
\end{enumerate}
Then it also holds that for all $t_0 \in [0,1]$ and $x \geq 4$,
\begin{equation}\label{tobound3}
\pr\bigg(\sup_{t \in [0,t_0]} \big( n t - \sum_{i=1}^{n} N_{e,i}(t) \big) \geq x\bigg)
\end{equation}
is at most $.8 \times (5.7 \frac{r_2 + 1}{\frac{r_2}{2} - 1})^{r_2 + 1} \times C_2 \times (n t_0)^{\frac{r_2}{2}} x^{-r_2}.$
\end{lemma}
\ \\We note that both Lemmas\ \ref{maximal2} and \ref{maximal3} will follow from a general maximal inequality of \citet{LS77} which converts bounds on the moments/tail of the partial sums of a stochastic process to bounds on the supremum of that process (see Lemma\ \ref{maximal1} of the supplemental appendix Section\ \ref{maximalssec}).  Let us also note that although pooled renewal processes are a special family of stochastic processes, they can still exhibit complex behaviors, and it is not clear how to prove the necessary explicit bounds (at the level of precision required to prove the desired $\frac{1}{1-\rho}$ scaling) without such general tools from probability theory.  Let us also point out that intuitively, statements of the form $\E\big[|\sum_{i=1}^{n} N_{e,i}(t) - n t|^r \big] \leq C (n t)^{\frac{r}{2}}$ capture a notion that the correlations/fluctuations of the process $\sum_{i=1}^{n} N_{e,i}(t)$ can be sufficiently controlled, which holds here due to the nice properties of renewal processes.

\subsubsection{Proof of Lemma\ \ref{alltimepooled1}.}

\proof{Proof of Lemma\ \ref{alltimepooled1}:}
By our assumptions and Lemma\ \ref{maximal2}, for all non-negative integers $k$ and $x > 0$, 
$$\pr\big(\max_{j \in \lbrace 1,\ldots,k \rbrace} \big|n j - \sum_{i=1}^{n} N_{e,i}(j)\big| \geq x\big) \leq 2.25 \times 2^{s} \times (2 \frac{r_1 + 1}{s - 1})^{r_1 + 1} \times C_1 n^{\frac{r_1}{2}} k^{s} x^{-r_1}.$$
Next, by our assumptions and Lemma\ \ref{maximal3}, for all $t_0 \in [0,1]$ and $x \geq 4$,
$$\pr\bigg(\sup_{t \in [0,t_0]} \big(n t - \sum_{i=1}^{n} N_{e,i}(t)\big) \geq x\bigg) \leq 
.8 \times (5.7 \frac{r_2 + 1}{\frac{r_2}{2} - 1})^{r_2 + 1} \times C_2 \times (n t_0)^{\frac{r_2}{2}} x^{-r_2}.
$$
It then follows from our assumptions that the conditions of Lemma\ \ref{TailBound} are met with $\phi(t) = n t - \sum_{i=1}^{n} N_{e,i}(t)$, $s,r_1, r_2, \nu$ their given values, $Z = 4$,
$$H_1 = 2.25 \times 2^{s} \times (2 \frac{r_1 + 1}{s - 1})^{r_1 + 1} \times C_1 n^{\frac{r_1}{2}}\ \ \ ,\ \ \ H_2 = .8 \times (5.7 \frac{r_2 + 1}{\frac{r_2}{2} - 1})^{r_2 + 1} C_2 n^{\frac{r_2}{2}}.$$
Combining the above with the implications of Lemma\ \ref{TailBound} and some straightforward algebra completes the proof.  $\halmos$
\endproof

\subsection{Proof of Lemma\ \ref{ggoldb222}.}\label{ggoldbsec}
To prevent having to make additional unnecessary (albeit minor) assumptions about the existence of steady-state distributions for different number of servers (both $n$ and $n'$ as opposed to only $n$), we first state a small variant of Lemma\ \ref{ggoldb} which also follows directly from the results of \citet{GG13}.  Let ${\mathcal Q}_{\textrm{res}}^{n'}$ be the FCFS $GI/GI/n'$ queue with inter-arrival times having the same distribution as r.v. $A$, service times having the same distribution as r.v. $S$, and the following initial conditions.  The time until the first arrival is distributed as $R(A)$, and there are exactly $n'$ jobs in service, with initial residual service times drawn i.i.d. distributed as $R(S)$, independent from the arrival process.  Let $Q_{\textrm{res}}^{n'}(t)$ denote the number in system at time $t$ in ${\mathcal Q}_{\textrm{res}}^{n'}$.  

\begin{lemma}[\citet{GG13}]\label{ggoldbtrans}
Suppose that $0 < \E[A],\E[S] < \infty$.  Then for all $t, x > 0$,
$$\pr\big( Q_{\textrm{res}}^{n'}(t) - n' \geq x \big) \leq \pr\Bigg( \sup_{0 \leq s \leq t} \bigg( A_e(s) - \sum_{i=1}^{n'} N_{e,i}(s) \bigg) \geq x \Bigg).
$$
\end{lemma}

We will also need to define an additional queueing system.  Let ${\mathcal Q}^{n,n'}_{\textrm{res}}$ be the FCFS $GI/GI/n$ queue with inter-arrival distribution $A,$ service time distribution $S$, and the following initial conditions.  The time until the first arrival is distributed as $R(A).$  There are exactly $n'$ jobs in service (note here $n - n'$ servers are initially empty in ${\mathcal Q}_{\textrm{res}}^{n,n'}$), with initial residual service times drawn i.i.d. distributed as $R(S)$, independent from the arrival process.  Let $Q^{n,n'}_{\textrm{res}}(t)$ denote the number in system at time $t$ in ${\mathcal Q}_{\textrm{res}}^{n,n'}$.  

Finally, we will need a well-known stochastic comparison result for multi-server queues.  In particular, it follows from known results in the stochastic comparison of multi-server queues as one varies the number of servers, see e.g. \citet{Berger92} Theorem 1 and also \citet{W81}, that a pathwise stochastic comparison holds between ${\mathcal Q}_{\textrm{res}}^{n,n'}$ and ${\mathcal Q}_{\textrm{res}}^{n'}.$  Here we only use the weaker implied distributional comparison, i.e. that 
\begin{equation}\label{comparesspd1}
\pr\big( Q^{n,n'}_{\textrm{res}}(t) \geq z \big) \leq \pr\big( Q^{n'}_{\textrm{res}}(t) \geq z \big)\ \textrm{for all}\ t,z \geq 0.
\end{equation}

With Lemma\ \ref{ggoldbtrans} and (\ref{comparesspd1}) in hand, we now complete the proof of Lemma\ \ref{ggoldb222}.
\proof{Proof of Lemma\ \ref{ggoldb222}:} Notice that ${\mathcal Q}^{n'}_{\textrm{res}}$ satisfies the conditions of Lemma\ \ref{ggoldbtrans}, and thus for all $t \geq 0$ and $z \geq 0$, 
$$\pr\big( Q^{n'}_{\textrm{res}}(t) - n' \geq z \big) \leq \pr\Bigg( \sup_{0 \leq s \leq t} \bigg( A_e(s) - \sum_{i=1}^{n'} N_{e,i}(s) \bigg) \geq z \Bigg).$$
Combining the above, we conclude that for all $t \geq 0$ and $x \geq n' - n$,
\begin{eqnarray*}
\pr\big( Q^{n,n'}_{\textrm{res}}(t) - n \geq x \big) &\leq& \pr\big( Q^{n'}_{\textrm{res}}(t) - n \geq x \big)
\\&=& \pr\big( Q^{n'}_{\textrm{res}}(t) - n' \geq x + n - n' \big)
\\&\leq& \pr\Bigg( \sup_{0 \leq s \leq t} \bigg( A_e(s) - \sum_{i=1}^{n'} N_{e,i}(s) \bigg)  \geq x + n - n'\Bigg).
\end{eqnarray*}
As our assumption that $Q^n(\infty)$ exists (which by our definitions must be independent of initial conditions) implies that $\lbrace Q^{n,n'}_{\textrm{res}}(t), t \geq 0 \rbrace$ converges in distribution to $Q^n(\infty)$, and applying the monotonicity of the supremum operator and continuity of probability measures, completes the proof.  $\halmos$
\endproof

\section*{Acknowledgements.}
The authors gratefully acknowledge support from NSF grant no. 1333457, as well as the anonymous referees for very thoughtful feedback, and helpful conversations with Mor Harchol-Balter, Isaac Grosof, Jamol Pender, Alan Scheller-Wolf, Ziv Scully, and Weina Wang.

\renewcommand*{\bibfont}{\small}

\newpage
\section{Supplemental Appendix.}\label{appsec}

\subsection{Bounds for higher order moments.} \label{highermomentsec}
Theorems\ \ref{mastertheoremsmallr} and \ref{mastertheoremlarger} also imply bounds for $\E[L^{s}(\infty)]$ for all $s < \frac{r^*}{2}$, by integrating the corresponding bounds for the tail of the queue length (and using the tail integral form for higher moments, see e.g. \citet{Nad.22}).  To most accurately state the implied bounds for higher moments, the relevant expression should be the integral of an infimum over $r$ (to derive a bound for each $x$, which is then integrated).  The associated statement is somewhat cumbersome, and thus we do not state it here out of considerations of readability.  Instead, we present an illustrative implication of our main results.  The result follows essentially immediately from our main result Theorem\ \ref{mastertheoremlarger}, the fact that $\E[X^r] = r \int_0^{\infty} x^{r-1} P(S > x) dx$ for any non-negative r.v. $X$ (see e.g. \citet{Nad.22}), the fact that $\Gamma(1 + x) \leq 1.5 x^x$ for all $x \geq 0$ (which follows from the bounds of \citet{Batir17} Theorem 2.3), and some straightforward algebra and calculus, and we omit the details.

\begin{corollary}\label{Highermomentbound1}
Under the same assumptions as Theorem\ \ref{mastertheoremlarger}, and supposing in addition that $r^* < \infty$ and $\E[S^{r^*}] < \infty$, for all $\epsilon \in (0,\frac{1}{4})$, it holds that $\E\big[L^{\frac{r^*}{2} - \epsilon}(\infty)\big]$ is at most
\begin{eqnarray*}
\ &\ &\ \Bigg( 4 \times 10^4 \times \frac{r^*}{\epsilon} \times \E[(S \mu_S)^2] \times \big(10^6\big)^{r^*} \times \bigg( \big(\E[(S \mu_S)^2]\big)^{r^*-1} \times {r^*}^{2.5 r^*} + {r^*}^{1.5 r^*} \times \E[(S \mu_S)^{r^*}] \bigg)
\\&\ \ \ &\ \ \ \ \ +\ 1.5 \times \big( 15 \sqrt{\E[(A \mu_A)^2]} \big)^{r^*} \times {r^*}^{.5 r^* + 1} \Bigg) \times (\frac{1}{1-\rho})^{\frac{r^*}{2} - \epsilon}.
\end{eqnarray*}
\end{corollary}

Let us make some additional clarifications regarding the above results.
\begin{itemize}
\item Analogous results for $\E[L^s(\infty)]$ also follow from our results (for all $s$) when $r^* = \infty$, although we do not present those results here.  Indeed, to present those results most accurately would require integrating an infimum, leading to somewhat combersome expressions.
\item The $(\frac{1}{1-\rho})^s$ scaling of $\E[L^s(\infty)]$ is consistent with known results for the single-server queue, and those limited settings of the multi-server queue where results are available.
\item It is easily verified that all prefactors appearing in Corollary\ \ref{Highermomentbound1} are asymptotically dominated by terms of the form ${r^*}^{O(r^*)}$.  Let us note that such an $r^{\Omega(r)}$ scaling is in fact unavoidable, as even for the $M/M/1$ queue with $\rho = \frac{1}{e}$, one has that (for integer $r \geq 2$)
\begin{eqnarray*}
\E\big[L^r(\infty)\big] &=& (1 - \frac{1}{e}) \sum_{k=1}^{\infty} k^r \times e^{-k}
\\&\geq& (1 - \frac{1}{e}) \int_{1}^{\infty} x^r \times e^{-x} dx 
\\&\geq& (1 - \frac{1}{e}) \big( \int_0^{\infty} x^r \times e^{-x} dx - 1 \big)
\\&=& (1 - \frac{1}{e}) \big( \int_0^{\infty} x^r \times e^{-x} dx - 1 \big)
\\&=& (1 - \frac{1}{e}) \big( r! - 1 \big),
\end{eqnarray*}
where it follows from well-known bounds for the factorial function that $r! = r ^{\Omega(r)}$ (\citet{Beesack69}).
\item It is known that under the assumptions of Corollary\ \ref{Highermomentbound1}, $\E[L^s(\infty)] < \infty$ for a range of $s$ strictly greater than $\frac{r^*}{2}$, and in fact the interplay between which moments of $S$ are finite, the number of servers, and which moments of $L(\infty)$ are finite is quite subtle (\citet{Scheller03,Scheller06}).  However, it is not known how these higher moments scale with $\frac{1}{1-\rho}$.  It is possible that at the level of generality considered in our results (where e.g. only few moments of $S$ may be finite), a complex behavior could in principle arise where some of these higher moments (say the $s$th moment) would no longer scale as $(\frac{1}{1-\rho})^s$ due to a complex interplay between $s,r^*,$ and $n$.  Better understanding the scaling of higher moments, and the relation of our results to those of \citet{Scheller03,Scheller06}, remains an interesting open question.
\item Although for the single-server queue simple recursive schemes are known for expressing higher moments of the steady-state waiting time in terms of $A$ and $S$ (see e.g. \citet{Tak62,Gong92}), even in that setting it is not clear that explicit bounds (not given as recursive formulas) have appeared previously in the literature, especially in the setting where only few moments of $S$ are finite.
\end{itemize}

\subsection{Implications for the Halfin-Whitt scaling regime and connections to \citet{Chawla17}.} \label{HalfinChawlasec}
In this section, we state several implications for queues in the Halfin-Whitt regime, in which $\rho$ scales as $1 - B n^{-\frac{1}{2}}$ for some excess capacity parameter $B > 0$.  Although all of our main results can be customized to this setting, here we present the illustrative example of our Theorem\ \ref{tailboundedmoments}.

\begin{corollary}\label{HWcor1}
Under the same assumptions and definitions as Theorem\ \ref{tailboundedmoments}, and supposing in addition that $\rho \leq 1 - B n^{-\frac{1}{2}}$ for some $B > 0$, the following holds.
\\\\For all $x > 0$, $\pr\big(  L(\infty) \geq x n^{\frac{1}{2}} \big)$ is at most 
$$f \times \exp\big( - \delta (B x)^{\frac{1}{3 + 2 c}} \big) + 1.1 \times \exp\bigg( - .0225 \big( \E[ (A \mu_A)^2 ] \big)^{-1} B x \bigg);$$
and the s.s.p.d. is at most
$$f \times \exp\big( - \delta B^{\frac{2}{3 + 2 c}} \big) + 1.1 \times \exp\bigg( - .0028 \big( \E[ (A \mu_A)^2 ] \big)^{-1} B^2 \bigg).$$
\end{corollary}

These results regarding the s.s.p.d. (as well as the other results for the s.s.p.d. in the Halfin-Whitt regime implied by our main results) make progress on an open question posed in \citet{Chawla17} related to distribution-independent bounds for the s.s.p.d. of multi-server queues in the Halfin-Whitt regime.  Indeed, those authors pose the question of whether the s.s.p.d. scales as $\exp\big(- \Omega(B^2) \big)$ for general service time distributions in the Halfin-Whitt regime.  Our results imply a bound of $\exp\big(-\Omega(B^{\alpha})\big)$ for some $\alpha \in (0,1)$ for quite general service time distributions, thus representing partial progress, but falling short of the $\exp\big(- \Omega(B^2) \big)$ scaling.  The question in its original form remains an interesting open problem, and we refer the reader to \citet{G16} for further related progress on this problem.

\subsection{Proof of Observation\ \ref{limitedgrowth} and Theorem\ \ref{tailboundedmoments}.}\label{optrsec}

\proof{Proof of Observation\ \ref{limitedgrowth}}
For simplicity (and as noted without loss of generality by a simple rescaling argument), suppose $\E[S] = 1$.  Suppose $P(S > x) \leq a \exp( - b x^c)$ for all $x > 0$.  Then by the tail integral form for higher moments, see e.g. \cite{Nad.22}, it follows that 
\begin{eqnarray*}
\E[S^r] &=& r \int_0^{\infty} x^{r-1} P(S > x) dx
\\&\leq&  r \int_0^{\infty} x^{r-1} a \exp( - b x^c ) dx
\\&=& a r \int_0^{\infty} x^{r-1} \exp( - b x^c ) dx
\\&=& a r \int_0^{\infty} x^{r-1} \exp( - (\frac{x}{b^{-\frac{1}{c}}})^c ) dx.
\end{eqnarray*}
Let $Z$ denote a Weibull r.v. with scale parameter $\alpha = b^{-\frac{1}{c}}$, shape parameter $\beta = c$, and location parameter 0, see \citet{Lehman63}.  Then $P(Z > x) = \exp( - (\frac{x}{b^{-\frac{1}{c}}})^c )$ for all $x > 0$, and combining with the above we find that $\E[S^r] = a \E[Z^r] = a \times  (b^{-\frac{1}{c}})^r \times \Gamma(1 + \frac{r}{c})$ (\citet{Lehman63}).  Combining with the fact that $\Gamma(1 + x) \leq 1.5 x^x$ for all $x \geq 0$ (which follows from the bounds of \citet{Batir17} Theorem 2.3), we find that 
\begin{eqnarray*}
\E[S^r] &\leq& a \times (b^{-\frac{1}{c}})^r \times 1.5 \times (\frac{r}{c})^{\frac{r}{c}}
\\&=& 1.5 \times a \times \big( (b c)^{-\frac{1}{c}}\big)^r \times r^{\frac{1}{c} r},
\end{eqnarray*}
completing the proof of the first part of the observation.  For the part regarding the fact that $\sup_{t > 0} \E[S - t | S > t] = \infty$ for a r.v. $S$ satisfying $P(S > x) = \exp( - b x^c)$ for all $x > 0$ and some $c < 1$, this follows from the results of \citet{Downey.91}, which imply that any r.v. with uniformly bounded mean residual life must have an exponentially decaying tail (in contrast to the Weibull distribution with $c < 1$ whose tails are heavier than any exponential).
$\Halmos$
\endproof

\proof{Proof of Theorem\ \ref{tailboundedmoments}}
We first prove the bound for $\pr\big(  L(\infty) \geq \frac{x}{1 - \rho} \big).$
First, let us derive a slightly simpler bound than that of Theorem\ \ref{mastertheoremlarger}, to simplify our analysis.  Note that as $c \geq 1$, the dominant term $\big(\E[(S \mu_S)^2]\big)^{r-1} \times r^{2.5 r} + r^{1.5 r} \ E[(S \mu_S)^r]$ will scale as $r^{(1.5 + c) r}$.  It is easily verified (from Theorem\ \ref{mastertheoremlarger}) that our results also imply $\pr\big(  L(\infty) \geq \frac{x}{1 - \rho} \big)$ is at most 
\begin{eqnarray}
\ &\ &\ \inf_{r > 2.5} \bigg( 4 \times 10^4 \times a \times \big(10^6 \E[(S \mu_S)^2] b \big)^r r^{(1.5 + c) r} \times x^{-\frac{r}{2}} \bigg) \label{optimizeme1}
\\&\ &\ + 1.1 \times \exp\bigg( - .0225 \big( \E[ (A \mu_A)^2 ] \big)^{-1} x \bigg), \nonumber
\end{eqnarray}
where the terms optimized over in (\ref{optimizeme1}) have the same $r^{(1.5 + c) r}$ scaling.  This is the bound which we will use in our analysis.  Let $f \stackrel{\Delta}{=} 4 \times 10^4 \times a, g \stackrel{\Delta}{=} 10^6 \E[(S \mu_S)^2] b , h \stackrel{\Delta}{=} 1.5 + c$.  Then we may rewrite the above as the fact that for all $x > 0$, $\pr\big(  L(\infty) \geq \frac{x}{1 - \rho} \big)$ is at most 
\begin{eqnarray}
\ &\ &\ \inf_{r > 2.5} \big(f \times g^r \times r^{h r} \times x^{-\frac{r}{2}} \big) \label{optimizeme1b}
\\&\ &\ + 1.1 \times \exp\bigg( - .0225 \big( \E[ (A \mu_A)^2 ] \big)^{-1} x \bigg), \nonumber
\end{eqnarray}
For a given $x > 0$, let $\hat{r}_x \stackrel{\Delta}{=} e^{-1} g^{-\frac{1}{h}} x^{\frac{1}{2 h}}$.  It turns out that looking at $\hat{r}_x$ is motivated by the fact that it can be shown to be the optimal solution to (\ref{optimizeme1b}) when the infimum is taken over $r > 0$ (instead of $r > 2.5$).  We will not need this fact here, although for completeness we prove this at the end of this section.  Note that for all $x$ s.t. $\hat{r}_x > 2.5$, the above implies that
$\pr\big(  L(\infty) \geq \frac{x}{1 - \rho} \big)$ is at most 
\begin{eqnarray}
\ &\ &\ f \times g^{\hat{r}_x} \times \hat{r}_x^{h \hat{r}_x} \times x^{-\frac{\hat{r}_x}{2}} \nonumber
\\&\ &\ + 1.1 \times \exp\bigg( - .0225 \big( \E[ (A \mu_A)^2 ] \big)^{-1} x \bigg), \nonumber
\\&=& f \exp\big( - e^{-1} h g^{-\frac{1}{h}} x^{\frac{1}{2 h}} \big)\label{optimizeme1c}
\\&\ &\ + 1.1 \times \exp\bigg( - .0225 \big( \E[ (A \mu_A)^2 ] \big)^{-1} x \bigg), \nonumber
\end{eqnarray}
completing the proof of Theorem\ \ref{tailboundedmoments} in that case.  However, note that $\hat{r}_x \leq 2.5$ iff  $x \leq (2.5 e)^{2 h} g^2$.  However, for any $x \leq (2.5 e)^{2 h} g^2$, it holds that 
\begin{eqnarray*}
f \exp\big( - e^{-1} h g^{-\frac{1}{h}} x^{\frac{1}{2 h}} \big) &\geq& f \exp\big( - e^{-1} h g^{-\frac{1}{h}} \times (2. 5 e g^{\frac{1}{h}}) \big)
\\&=& f \exp( - 2.5 h )\ \ \ \geq\ \ \ 1,
\end{eqnarray*}
where the final inequality can be easily verified for our particular choice of $f,h$.  As $\pr\big(  L(\infty) \geq \frac{x}{1 - \rho} \big) \leq 1$ in all cases, it follows that the desired bound holds for all $x$, completing the proof.  We note that in our actual result we have presented a slightly weaker bound so we do not need to define different constants in the statement of our result for the queue length and the s.s.p.d., and plugged in the bound implied by our assumptions for $\E[(S \mu_S)^2]$.
\\\\The proof for the s.s.p.d. follows by applying a nearly identical argument to the slightly different bound (implied by Theorem\ \ref{mastertheoremlarger} for the s.s.p.d.)
\begin{eqnarray*}
\ &\ &\ \inf_{r > 2.5} \bigg( 4 \times 10^4 \times a \times \big(10^7 \E[(S \mu_S)^2] b \big)^r r^{(1.5 + c) r} \times \big( n (1 - \rho)^2 \big)^{-\frac{r}{2}} \bigg)
\\&\ &\ + 1.1 \times \exp\bigg( - .0028 \big( \E[ (A \mu_A)^2 ] \big)^{-1} n (1 - \rho)^2 \bigg),
\end{eqnarray*}
and we omit the details.  Combining the above completes the proof. $\halmos$ \endproof

\subsubsection{Proof that $\hat{r}_x$ is the claimed minimizer.}
Here, for completeness (and to motivate our use of $\hat{r}_x$), we prove that $\hat{r}_x$ is the the claimed minimizer.
\begin{lemma}
Suppose $f,g,h \geq 1$ are some real numbers.  Then for all $x > 1$, $\hat{r}_x \in \argmin_{r > 0}\big( f \times g^r \times r^{h r} \times x^{-\frac{r}{2}}\big)$.
\end{lemma}
\proof{Proof : }
As it will have the same minimizer, let us take logarithms, and consider finding the $r$ minimizing
$\log(f) + r \log(g) + h r \log(r) - \frac{r}{2} \log(x).$  As $$\frac{d}{dr} \big( \log(f) + r \log(g) + h r \log(r) - \frac{r}{2} \log(x) \big) = \log(g) + h \log(r) + h - \frac{1}{2} \log(x),$$ we find that $\log(f) + r \log(g) + h r \log(r) - \frac{r}{2} \log(x)$ is a convex function of $r$, which will attain its minimum where $\log(g) + h \log(r) + h - \frac{1}{2} \log(x) = 0$, i.e. at 
$$\hat{r} \stackrel{\Delta}{=}  \exp\big( \frac{1}{2 h} \log(x) - \frac{\log(g)}{h} - 1\big) = \hat{r}_x,$$ 
completing the proof.  $\halmos$ \endproof

\subsection{Proof of Corollary\ \ref{meanwaitcorbeyond}.}\label{meanwaitcorbeyondsec}
\proof{Proof of Corollary\ \ref{meanwaitcorbeyond} : }
Combining the two bounds of Corollary\ \ref{mastertheorem0ccc1}, along with some straightforward algebra and the tail integral formula for expected value, we conclude that $(1 - \rho) \E[L(\infty)]$ is at most
\begin{eqnarray*}
\ &\ &\ \int_0^{n (1 - \rho)^2} \Bigg(8 \times 10^{28} \times \E[(S \mu_S)^2] \times \bigg( \big(\E[(S \mu_S)^2]\big)^2 + \E[(S \mu_S)^3] \bigg) \times \big( n (1 - \rho)^2 \big)^{-1.5}
\\&\ \ \ \ \ \ \ \ &\ \ \ \ \ \ \ +\ \ 1.1 \times \exp\bigg( - .0028 \big( \E[ (A \mu_A)^2 ] \big)^{-1} n (1 - \rho)^2 \bigg) \Bigg) dx
\\&\ &\ +\ \int_{n (1 - \rho)^2}^{\infty} \Bigg( 8 \times 10^{25} \times \E[(S \mu_S)^2] \times \bigg( \big(\E[(S \mu_S)^2]\big)^2 + \E[(S \mu_S)^3] \bigg) \times x^{-1.5}
\\&\ \ \ \ \ \ \ \ &\ \ \ \ \ \ \ +\ \ 1.1 \times \exp\bigg( - .0225 \big( \E[ (A \mu_A)^2 ] \big)^{-1} x \bigg) \Bigg) dx
\\&\leq& 8 \times 10^{28} \times \E[(S \mu_S)^2] \times \bigg( \big(\E[(S \mu_S)^2]\big)^2 + \E[(S \mu_S)^3] \bigg) \times \big( n (1 - \rho)^2 \big)^{-.5}
\\&\ &\ \ +\ 1.1 \times n (1 - \rho)^2 \times \exp\bigg( - .0028 \big( \E[ (A \mu_A)^2 ] \big)^{-1} n (1 - \rho)^2 \bigg)
\\&\ &\ \ +\ 1.6 \times 10^{26} \times \E[(S \mu_S)^2] \times \bigg( \big(\E[(S \mu_S)^2]\big)^2 + \E[(S \mu_S)^3] \bigg) \times \big( n (1 - \rho)^2 \big)^{-.5}
\\&\ &\ \ +\ 49 \E[ (A \mu_A)^2] \exp\bigg( - .0225 \big( \E[ (A \mu_A)^2 ] \big)^{-1} \big(n (1 - \rho)^2\big) \bigg)
\\&\leq&\ 8.001 \times 10^{28} \times \E[(S \mu_S)^2] \times \bigg( \big(\E[(S \mu_S)^2]\big)^2 + \E[(S \mu_S)^3] \bigg) \times \big( n (1 - \rho)^2 \big)^{-.5}
\\&\ &\ \ \ +\ \big(1.1 \times n (1 - \rho)^2 +  49 \E[ (A \mu_A)^2] \big) \exp\bigg( - .0028 \big( \E[ (A \mu_A)^2 ] \big)^{-1} n (1 - \rho)^2 \bigg),
\end{eqnarray*}
and thus $\E[L(\infty)]$ is at most 
\begin{eqnarray*}
\ &\ &\ \Bigg(8.001 \times 10^{28} \times \E[(S \mu_S)^2] \times \bigg( \big(\E[(S \mu_S)^2]\big)^2 + \E[(S \mu_S)^3] \bigg) \times \big( n (1 - \rho)^2 \big)^{-.5}
\\&\ &\ \ \ + \big(1.1 \times n (1 - \rho)^2 +  49 \E[ (A \mu_A)^2] \big) \exp\bigg( - .0028 \big( \E[ (A \mu_A)^2 ] \big)^{-1} n (1 - \rho)^2 \bigg) \Bigg) \times \frac{1}{1-\rho}.
\end{eqnarray*}
Combining with the fact that $e^{-x} \leq \frac{2}{x^2}$ for all $x > 0$ (by a simple Taylor series expansion), which implies (after some straightforward algebra) that $$\big(1.1 \times n (1 - \rho)^2 +  49 \E[ (A \mu_A)^2] \big) \exp\bigg( - .0028 \big( \E[ (A \mu_A)^2 ] \big)^{-1} n (1 - \rho)^2 \bigg)$$
is at most
$2 \times 49 \E[(A \mu_A)^2] \times \big (n (1-\rho)^2 \big)^{-.5}$ if $n(1-\rho)^2 \geq 10^6 \times \big(\E[(A \mu_A)^2]\big)^2$ completes the proof.  $\Halmos$. \endproof

\subsection{Proofs of Corollary\ \ref{rnotsmalltheorema} and Theorem\ \ref{rnotsmalltheoremb}, and further discussion of prefactors and $r^{\Omega(r)}$ scaling.}\label{rnotsmallsec}
\subsubsection{Proof of Corollary\ \ref{rnotsmalltheorema}.}
First, let us complete the proof of Corollary\ \ref{rnotsmalltheorema}, a natural implication of our main results very similar to the types of bounds which appeared in an earlier version of this manuscript (\citet{G17c}).
\ \proof{Proof of Corollary\ \ref{rnotsmalltheorema} :}
Since $\E[S^r] < \infty$, Theorem\ \ref{mastertheoremlarger} implies that $P\big( L(\infty) \geq \frac{x}{1-\rho} \big)$ is at most
\begin{eqnarray*}
\ &\ &\ 2 \times 10^4 \times \E[(S \mu_S)^2] \times \big(10^6\big)^r \times \bigg( \big(\E[(S \mu_S)^2]\big)^{r-1} \times r^{2.5 r} + r^{1.5 r} \times \E[(S \mu_S)^r] \bigg) \times x^{-\frac{r}{2}} 
\\&\ \ \ &\ \ +\ \ 1.1 \times \exp\bigg( - .0225 \big( \E[ (A \mu_A)^2 ] \big)^{-1} x \bigg).
\end{eqnarray*}
A Taylor expansion implies that 
\begin{eqnarray*}
\ &\ &\ 1.1 \exp\bigg( - .0225 \big( \E[ (A \mu_A)^2 ] \big)^{-1} x \bigg) \leq 
\\&\leq& 1.1 \times \lceil \frac{r}{2} \rceil! \times (\frac{1}{.0225})^{\frac{r}{2} + 1} \times \big( \E[(A \mu_A)^2] \big)^{\lceil \frac{r}{2} \rceil} \times x^{-\frac{r}{2}}.
\end{eqnarray*}
Further applying the fact that $\lceil x \rceil! \leq 2 x^x$ for all $x \geq 1$ (which follows from known bounds for the factorial function and some straightforward algebra, see e.g. \citet{Beesack69,Batir17}), we conclude that 
$$
1.1 \times \exp\bigg( - .0225 \big( \E[ (A \mu_A)^2 ] \big)^{-1} x \bigg)
\leq
108 \times 5^r \times r^{\frac{r}{2}} \times \big( \E[ (A \mu_A)^2 ] \big)^{\frac{r}{2}}.$$
Further combining with some straightforward algebra completes the proof.$\Halmos$ \endproof
\subsubsection{Proof of Theorem\ \ref{rnotsmalltheoremb}.}
The intuition of our proof is actually quite simple.  Intuitively, we will show that for single-server queues in which $A$ and $S$ are uniformly bounded, $\E[A^r]$ and $\E[S^r]$ grow only as $\exp\big( O(r) \big)$, while $\E[L^r(\infty)]$ grows as $r^{\Omega(r)}$ (essentially inheriting the moment sequence of an exponential distribution, since even when $A,S$ are uniformly bounded the queue length itself is not and will have an exponential decay by known resuilts for single-server queues).  The desired result will immediately follow, since if $\underline{c_r}$ did not similarly scale as $r^{\Omega(r)}$ a contradiction would be reached.  Before proving Theorem\ \ref{rnotsmalltheoremb}, let us recall a lower bound on the tail of the waiting time in a single-server queue proven in \citet{K70}.
\begin{lemma}[\citet{K70}]\label{kinglbaa} Suppose that : (1) $n = 1$; (2)  there exist $\theta^* > 0$ s.t. $\E\big[\exp\big(\theta^*(S-A)\big)\big] = 1$; and (3) there exists a finite constant $B$ s.t. $P(S \leq B) = 1$.  Then for all $x > 0$, $P\bigg( \sup_{n \geq 0} \big( \sum_{i=1}^n (S_i - A_i) \big) \geq x \bigg) \geq \exp(- \theta^* B) \exp(- \theta^* x)$.
\end{lemma}
With Lemma\ \ref{kinglbaa} in hand, we now complete the proof of Theorem\ \ref{rnotsmalltheoremb}.
\proof{Proof of Theorem\ \ref{rnotsmalltheoremb} : }
Consider the single-server queue in which $A$ is distributed uniformly on $[0,2]$, and $S$ is distributed uniformly on $[0,1]$.  It is easily verified that there exists a unique strictly positive $\theta^* \sim 2.851$ s.t. $\E\big[\exp\big(\theta^*(S-A)\big)\big] = 1$.  It follows from standard queueing results appearing in e.g. \citet{AS08} that $\textrm{Work}(\infty), W(\infty),$ and $Q(\infty)$ exist; and for all $x > 0$, $P\big( \textrm{Work}(\infty) > x \big) = P\big( W(\infty) + S - R(A) > x \big)$, with $W(\infty), S, R(A)$ independent r.v.s.  Using the fact that $P(S \leq 1) = 1$ and $P(A \leq 2) = 1$, it easily follows that for all $x > 0$ : (1)\ $P\big( \textrm{Work}(\infty) \geq x \big) \geq P\big( W(\infty) \geq x + 2 \big)$; (2)\ $P\big( L(\infty) \geq x \big) \geq P\big( \lfloor \textrm{Work}(\infty) \rfloor - 1 \geq x \big) \geq P\big( \textrm{Work}(\infty) \geq x + 2\big)$.  It follows that for all $x > 0$, $P\big( L(\infty) \geq x \big) \geq  P\big( W(\infty) \geq x + 4 \big).$  Combining with Lemma\ \ref{kinglbaa}, we conclude that for all $x > 0$, $P\big( L(\infty) \geq x \big) \geq e^{- \theta^* (x + 5)} \geq \exp(- 14.3) \times \exp(-2.86 x)$.  It then follows from the tail integral form for higher moments (\cite{Nad.22}) and known bounds for the factorial function (\citet{Beesack69,Batir17}) that for all $r > 16$,
\begin{eqnarray*}
\E[L^{\frac{r}{4}}(\infty)] &\geq& \frac{r}{4} \int_0^{\infty} x^{\frac{r}{4} - 1} \exp(- 14.3) \times \exp(-2.86 x) dx
\\&\geq& \exp(-14.3) \times \frac{r}{4} \times 2.86^{-\frac{r}{4}} \times \lfloor \frac{r}{4} - 1 \rfloor!
\\&\geq& \exp(-14.3) \times \frac{r}{4} \times 2.86^{-\frac{r}{4}} \times \lfloor \frac{r}{4} - 1 \rfloor^{\lfloor \frac{r}{4} - 1 \rfloor} \times \exp(- \lfloor \frac{r}{4} - 1 \rfloor)
\\&\geq& \exp(-14.3) \times \frac{r}{4} \times (2.86 e)^{-\frac{r}{4}} \times (\frac{r}{8})^{\frac{r}{8}} 
\\&\geq& \exp(-14.3) \times 2.2^{-r} \times r^{\frac{1}{8} r}.
\end{eqnarray*}
Now, suppose that for all $r > 16$ and $x > 0$, it holds that $P\big( L(\infty) \geq \frac{x}{1 - \rho} \big) \leq c_r \times \bigg(\big( \E[(A \mu_A)^2] \big)^r \times \E[(A \mu_A)^r] + \big(\E[(S \mu_S)^2] \big)^r \times \E[ (S \mu_S)^r ] \bigg) \times x^{-\frac{r}{2}}$ for $\rho,A,S$ as in the above single-server queue.  Note that $\rho = \frac{1}{2}$, and using standard results for the moments of a uniform r.v. we have $\E[(A \mu_A)^r] = \E[(S \mu_S)^r] = \frac{2^r}{r+1} \leq 2^r$, and $\big(\E[(A \mu_A)^2]\big)^r = \big(\E[(S \mu_S)^2]\big)^r = (\frac{4}{3})^r$.  It then follows from some straightforward algebra that for all $x > 0$, $P\big( L(\infty) \geq x \big) \leq 2 \times c_r \times (\frac{8}{3})^r \times x^{-\frac{r}{2}}.$  Applying the tail integral form for higher moments, see e.g. \cite{Nad.22}, we thus have that 
$E\big[ L^\frac{r}{4}(\infty)] \leq \frac{r}{4} + 2 \times c_r \times \frac{r}{4} \times (\frac{8}{3})^r
\times \int_1^{\infty}  x^{\frac{r}{4} - 1} \times x^{-\frac{r}{2}} dx \leq  \frac{r}{4} + c_r \times 2 \times (\frac{8}{3})^r$.
\\\indent Combining the above lower and upper bounds for $\E[L^{\frac{r}{4}}(\infty)]$, it follows that $\frac{r}{4} + c_r \times 2 \times (\frac{8}{3})^r
\geq \exp(-14.3) \times 2.2^{-r} \times r^{\frac{1}{8} r}.$  A straightforward contradiction argument then implies that there exists $\epsilon > 0$ s.t. $c_r \geq \epsilon \times r^{\frac{1}{9} r}$ for all $r > 16$, which completes the proof.  $\Halmos$ \endproof
\subsubsection{Where in our proofs does the $r^{\Omega(r)}$ scaling arise?}
We now comment explicitly on where in our proofs the $r^{\Omega(r)}$ scaling arises.  Although our proofs involve several moving parts and multiple bounds which are composed together, the vast majority of those bounds do not actually contribute terms scaling as $r^{\Omega(r)}.$  Essentially all aspects of our proof which contribute to the $r^{\Omega(r)}$ scaling can ultimately be attributed to our having to explicitly bound $\E\big[|\sum_{i=1}^n N_{e,i}(t) - n \mu_S t|^r\big]$, which we do in 
our Lemma\ \ref{poolerbound1} for the case $t \geq 1$.  It follows from some straightforward algebra that our Lemma\ \ref{poolerbound1} can be loosely interpreted as asserting the following (which our Lemma\ \ref{poolerbound1} can indeed be formally shown to imply using some straightforward algebra, the details of which we omit).
\begin{corollary}\label{pooledlba}
For each $r > 2$, there exists a finite constant $c_r$ (depending only on $r$) s.t. for all integers $n \geq 1$, all $t \geq 1$, and all $S$ s.t. $\E[S^r] < \infty$, it holds hat  $\E\big[|\sum_{i=1}^n N_{e,i}(t) - n \mu_S t|^r\big] \leq c_r \times \big(\E[(S \mu_S)^2]\big)^r \times \E[(S \mu_S)^r] \times (n t)^{\frac{r}{2}}.$  Furthermore, one can take $c_r = r^{O(r)}$.
\end{corollary}  
We now show that the $r^{\Omega(r)}$ scaling in Corollary\ \ref{pooledlba} is in fact unavoidable, even for the case $n = 1, t= 1$.  Intuitively, this will follow from the simple fact that athough the $r$th moment of a r.v. which puts probability on 0 and 2 scales at most exponentially in $r$, the $r$th central moment of the associated renewal process scales as $r^{\Omega(r)}$ (inherited from a related geometrically distributed r.v.).  As such bounds for pooled renewal processes (with the $(nt)^{\frac{r}{2}}$ scaling) are essential to implementing our overall approach, this is further suggestive of the fact that avoiding constants scaling in this way may require fundamentally different approaches, and/or imposing additional assumptions on $S$.
\begin{lemma}\label{pooledlbb}
let $\underline{c_r}$ denote the infimum of all constants $c_{r}$ for which the bound of Corollary\ \ref{pooledlba} holds.  Then there exists an absolute constant $\epsilon > 0$ s.t. $\underline{c_{r}} \geq \epsilon \times r^{\frac{1}{2} r}$ for all $r > 16$.
\end{lemma}
\proof{Proof :} Let $S$ be the r.v. s.t. $P(S = 0) = P(S = 2) = \frac{1}{2}$.  Note that $R(S)$ is uniformly distributed on the interval [0,2].  Note also that $P\big( N_o(1) \geq k \big) \geq 2^{-k}$ for all $k \geq 1$, and it follows that $P\big( N_e(1) - 1 \geq k \big) \geq 2^{-(k+3)}$ for all $k \geq 1$.  Thus applying the tail integral form for higher moments (\citet{Nad.22}), it follows that for all $r > 2$,
\begin{eqnarray*}
\E\big[|N_{e}(1) - 1|^r\big] &\geq& \int_0^{\infty} r x^{r-1} 2^{-\lceil x + 3 \rceil}
\\&\geq& \frac{r}{16} \int_0^{\infty} x^{r-1} 2^{-x}
\\&=& \big( \log_e(2) \big)^{-r} \Gamma(r).
\end{eqnarray*}
The desired result then follows from some straightforward asymptotics and standard bounds for the Gamma function (\citet{Beesack69,Batir17}), the details of which we omit.  $\Halmos$ \endproof
\ \\We next comment briefly on where precisely in our proofs associated with the centered moments of pooled renewal processes the $r^{\Omega(r)}$ scaling arises.  First, in 
bounding $\E[|N_o(t) - \mu_S t|^r]$, we apply the Burkholder-Rosenthal ineuqality in our Lemma\ \ref{rosenthal1}, to bound the higher moments of certain martingale-related terms.  Our application of this inequality is consistent with the approach sketched in \citet{Gut09}, which our analysis makes completely explicit.  Second, in bounding $\E\big[|\sum_{i=1}^n N_{e,i}(t) - n \mu_S t|^r\big]$ for $t \geq 1$, we apply the Marcinkiewicz–Zygmund inequality in our Lemma\ \ref{csumbound} to convert out bounds on $\E\big[|N_e(t) - \mu_s t|^r\big]$ into bounds for the corresponding pooled process.  This inequality is a powerful tool for converting bounds for the higher moments of individual mean-zero r.v.s into bounds for the higher moments of the sums of those r.v.s.  Although the family of renewal processes have special structure, and much is known about their general asymptotic scaling, the associated (pooled) counting processes can still exhibit complex behaviors, where the analyses of these behaviors has been at the core of several recent analyses of multi-server queues (see e.g. \citet{Bazhba19}).  As the $\frac{1}{1-\rho}$ bounds we prove are very sensitive to how all aspects of our proof scale in $n$ and $t$, it is not clear whether it is possible to explicitly and non-asymptotically bound the central moments of pooled renewal processes at the level of generality we consider without applying such general inequalities from probability theory.  Let us point out that additional $r^{\Omega(r)}$ scaling arises in the analysis for $t \leq 1$, e.g. in our Lemma\ \ref{csumboundbb2}, again due to the application of general inequalities from probability theory.
\subsubsection{Summary of discussion of prefactors arising in our bounds.}
In summary, it remains an interesting open question whether the $r^{\Omega(r)}$ prefactors arising in our main results are fundamental, or an artifact of our analysis.  That said, we have proven that for very closely related bounds implied by our main results, the $r^{\Omega(r)}$ results are indeed fundamental.  This arises at least in part due to the level of generality of our main results, e.g. that they hold for both the setting that $A,S$ are uniformly bounded, as well as the setting that $A,S$ have quite heavy tails with few finite moments.   It would be very interesting to derive other qualitatively different bounds for multi-server queues with $\frac{1}{1-\rho}$ scaling, possibly under different and/or stronger assumptions and using fundamentally different methods of analysis.  We note that our Theorem\ \ref{tailboundedmoments}, in which we prove stronger tail bounds with fundamentally different behavior and no dependence on any $r$ parameter, as well as our Theorems\ \ref{sspdnoprefactor},\ \ref{sspdnoprefactora2}, \ref{sspdnoprefactorb}, and\ \ref{sspdnoprefactorc} based on drift arguments, represent a step in this direction.

\subsection{Proofs of Theorems\ \ref{sspdnoprefactor},\ \ref{sspdnoprefactora2}, \ref{sspdnoprefactorb}, and\ \ref{sspdnoprefactorc}.}\label{noprefactorsec0}
\subsubsection{Proof of Theorem\ \ref{sspdnoprefactor}.}\label{noprefactorsec1}
Our proof proceeds by noting that $\E[\textrm{Num}_{\textrm{service}}(\infty)] = \frac{\mu_A}{\mu_S}$ implies $\E\big[ \max\big(0, \textrm{Num}_{\textrm{service}}(\infty) - \frac{\mu_A}{\mu_S}\big)\big] = \E\big[ \max\big(0, \frac{\mu_A}{\mu_S} - \textrm{Num}_{\textrm{service}}(\infty)\big)\big]$, and then relating the left hand side of this equality to the s.s.p.d. and the right hand side to the well-understood $M/GI/\infty$ queue.  We again note that the proof is very similar to those of \citet{wang2021zero,hong2021sharp}, which were primarily focused on the study of more general models.  
\proof{Proof of Theorem\ \ref{sspdnoprefactor}:}
Let $N$ denote $\textrm{Num}^n_{\textrm{service}}(\infty)$, the steady-state number of busy servers.  It follows from standard Little's Law type conservation arguments (\citet{Wolff2011,Heyman82,Whitt.83c}) that $\E[N] = \frac{\mu_A}{\mu_S}.$  Since for any integrable r.v. $X$, it holds that $\E[X] = \E\big[ \max(0,X) \big] - \E\big[ \max(0,-X) \big]$, we conclude that
\begin{equation}\label{noprefactoreq1}
\E\big[ \max(0, N - \frac{\mu_A}{\mu_S}) \big] = \E\big[ \max(0, \frac{\mu_A}{\mu_S} - N) \big].
\end{equation}
Recall that $Q^n(\infty)$, which we will denote simply by $Q$, is a r.v. distributed as the steady-state total number in system.  Since $\frac{\mu_A}{\mu_S} < n$, and the basic dynamics of the FCFS GI/GI/n queue imply that one can construct $N$ and $Q$ on a common probability space s.t. $I(N = k) = I(Q = k)$ for $k \in \lbrace 0,\ldots,n-1 \rbrace$, it follows that 
\begin{equation}\label{noprefactoreq2}
\E\big[ \max(0, \frac{\mu_A}{\mu_S} - N) \big] = \E\big[ \max(0, \frac{\mu_A}{\mu_S} - Q) \big].
\end{equation}
Let $Q_{\infty}$ denote a r.v. distributed as the steady-state total number in system in an $M/GI/\infty$ queue with the same inter-arrival and service time distribution as ${\mathcal Q}^n$.  It follows from standard and well-known stochastic comparison results between multi-server and infinite-server queues, see e.g. \citet{Whitt.00,GG13,wang2021zero,hong2021sharp}, that $P\big( Q \geq x \big) \geq P\big( Q_{\infty} \geq x \big)$ for all $x \in {\mathcal R}$.  As the function $f(x) \stackrel{\Delta}{=} 
\max(0, \frac{\mu_A}{\mu_S} - x)$ is non-increasing in $x$, we obseve (as in \citet{wang2021zero, hong2021sharp}) that the basic properties of stochastic dominance (see e.g. \citet{Brumelle75}) thus imply
$$\E\big[ \max(0, \frac{\mu_A}{\mu_S} - Q) \big] \leq \E\big[ \max(0, \frac{\mu_A}{\mu_S} - Q_{\infty}) \big].$$
Combining with (\ref{noprefactoreq1}) - (\ref{noprefactoreq2}), we conclude that
\begin{equation}\label{noprefactoreq3}
\E\big[ \max(0, N - \frac{\mu_A}{\mu_S}) \big] \leq \E\big[ \max(0, \frac{\mu_A}{\mu_S} - Q_{\infty}) \big].
\end{equation}
Next, very similar to \citet{wang2021zero, hong2021sharp}, we observe that since by assumption $\frac{\mu_A}{\mu_S} < n$, non-negativity implies
\begin{eqnarray}
\E\big[ \max(0, N - \frac{\mu_A}{\mu_S}) \big] &=& \E\big[ (N - \frac{\mu_A}{\mu_S}) I(N \geq \frac{\mu_A}{\mu_S}) \big] \nonumber
\\&\geq& \E\big[ (N - \frac{\mu_A}{\mu_S}) I(N \geq n) \big] \nonumber
\\&=& \E\big[ (n - \frac{\mu_A}{\mu_S}) I(N = n) \big]\ \ \ =\ \ \ (n - \frac{\mu_A}{\mu_S}) P(N \geq n) \label{noprefactoreq4},
\end{eqnarray}
where we have used the fact that the dynamics of a FCFS $GI/GI/n$ queue imply $\lbrace N \geq n \rbrace$ iff $\lbrace N = n \rbrace$.  Combining (\ref{noprefactoreq3}) - (\ref{noprefactoreq4}), we conclude that
\begin{equation}\label{noprefactoreq5}
P(N \geq n) \leq (n - \frac{\mu_A}{\mu_S})^{-1} \E\big[ \max(0, \frac{\mu_A}{\mu_S} - Q_{\infty}) \big].
\end{equation}
Using the well-known fact that $Q_{\infty}$ has a Poisson distribution with mean $\frac{\mu_A}{\mu_S},$ and letting $\textrm{Poi}$ denote a r.v. with this distribution, we have that
\begin{equation}\label{noprefactoreq6}
P(N \geq n) \leq (n - \frac{\mu_A}{\mu_S})^{-1} \E\big[ \max(0, \frac{\mu_A}{\mu_S} - \textrm{Poi}) \big].
\end{equation}
Next, as in \citet{wang2021zero, hong2021sharp}, we use a simple application of Jensen's inequality to bound $\E\big[ \max(0, \frac{\mu_A}{\mu_S} - \textrm{Poi}) \big]$ in terms of the standard deviation of $\textrm{Poi}$, which equals $\sqrt{\frac{\mu_A}{\mu_S}}$.  More formally, we proceed as follows.  Recall that $\E[X] = \E\big[ \max(0,X) \big] - \E\big[ \max(0,-X) \big]$ for a general integrable r.v. $X$, and for essentially identical reasons $\E[|X|] = \E\big[ \max(0,X) \big] + \E\big[ \max(0,-X) \big]$ for a general integrable r.v. $X$.  Applying with $X = \textrm{Poi} - \frac{\mu_A}{\mu_S}$, we conclude that $\E\big[ \max(0, \frac{\mu_A}{\mu_S} - \textrm{Poi}) \big] = \frac{1}{2} \E\big[ \big| \textrm{Poi} - \frac{\mu_A}{\mu_S} \big| \big]$.  As Jensen's inequality implies $\E\big[ \big| \textrm{Poi} - \frac{\mu_A}{\mu_S} \big| \big] \leq \sqrt{ \E\big[ \big( \textrm{Poi} - \frac{\mu_A}{\mu_S} \big)^2 \big] },$ and the variance of $\textrm{Poi}$ equals  $\frac{\mu_A}{\mu_S}$, we may combine with (\ref{noprefactoreq6}) to conclude that
\begin{equation}\label{noprefactoreq70}
\E\big[ \max(0, \frac{\mu_A}{\mu_S} - \textrm{Poi}) \big] \leq \frac{1}{2} \sqrt{\frac{\mu_A}{\mu_S}},
\end{equation}
and
\begin{equation}\label{noprefactoreq7}
P(N \geq n) \leq \frac{1}{2} \frac{ \sqrt{\frac{\mu_A}{\mu_S}} }{ n - \frac{\mu_A}{\mu_S} }.
\end{equation}
Combining with the fact that the s.s.p.d. equals $P(N \geq n)$, along with some straightforward algebra which yields 
\begin{eqnarray*}
\frac{ \sqrt{\frac{\mu_A}{\mu_S}} }{ n - \frac{\mu_A}{\mu_S} } &=& \frac{ \sqrt{\frac{\mu_A}{ n \mu_S}} }{ \sqrt{n} - \sqrt{n} \frac{\mu_A}{ n \mu_S} }
\\&=& \frac{ \sqrt{\rho} }{ \sqrt{n} (1 - \rho) },
\end{eqnarray*}
completes the proof.  $\halmos$
\endproof
\subsubsection{Proof of Theorem\ \ref{sspdnoprefactora2}.}\label{noprefactorsec15}
It follows from our proof of Theorem\ \ref{sspdnoprefactor}, specifically Equations (\ref{noprefactoreq3}) and (\ref{noprefactoreq70}), that $\E\big[ \max(0, N - \frac{\mu_A}{\mu_S}) \big] \leq \frac{1}{2} \sqrt{\frac{\mu_A}{\mu_S}}$.  The desired result then follows from Markov's inequality.  $\Halmos$ \endproof

\subsubsection{Proof of Theorem\ \ref{sspdnoprefactorb}.}\label{noprefactorsec2}
Our proof proceeds by relating the long-run fraction of jobs that abandon the system to $\theta \times \E\big[ L^n_a(\infty) \big]$, and then using known stochastic comparison results for multi-server systems with abandonments to further bound this quantity in terms of the well-understood Erlang loss model.  Let us also point out that although one might think that the case of abandonments would be ``more challenging", in this case it is actually simpler for such drift arguments, as the queue length manifests more directly in the rate at which jobs depart.  

\proof{Proof of Theorem\ \ref{sspdnoprefactorb}:}
Let $L^n_a(t)$ denote the number of jobs waiting in queue in ${\mathcal Q}^n_a$ at time $t$, and $\textrm{Aban}(t)$ denote the number of jobs that abandon from ${\mathcal Q}^n_a$ on $[0,t].$  It follows from standard Poisson constructions for $M/PH/n + M$ queues, see e.g. \citet{Dai10}, that $\E[\textrm{Aban}(t)] = \theta \E[\int_0^t L^n_a(s) ds]$ for all $t \geq 0$.  The positive Harris recurrence proven in \citet{DDG14}, along with standard implications of positive Harris recurrence for countable-state continuous time Markov chains, imply that $\lim_{t \rightarrow \infty} t^{-1} \E[\int_0^t L^n_a(s) ds] = \E[L^n_a(\infty)]$ and 
\begin{equation}\label{abandoneq0}
\lim_{t \rightarrow \infty} t^{-1}\E[\textrm{Aban}(t)] = \theta \E[L^n_a(\infty)].
\end{equation}
Let ${\mathcal Q}^n_{\textrm{loss}}$ denote an $n$-server Erlang loss model with the same inter-arrival and service distribution as ${\mathcal Q}^n_a$, also initially empty.  Note that ${\mathcal Q}^n_{\textrm{loss}}$ is equivalent to a $M/PH/n + GI$ queue in which patience times are w.p.1 equal to zero.  For a more formal review of this family of well-studied systems, we refer the reader to \citet{Davis.95,Sev57,Franken82}.  Let $\textrm{Loss}(t)$ denote the number of jobs that abandon from ${\mathcal Q}^n_{\textrm{loss}}$ on $[0,t]$.  It follows from the stochastic comparison results of \citet{Bhat91}, specifically Theorem 3.1 of that work, that $\E[\textrm{Aban}(t)] \leq \E[\textrm{Loss}(t)]$ for all $t \geq 0$, and thus by (\ref{abandoneq0})
\begin{equation}\label{abandoneq0b}
\E[L^n_a(\infty)] \leq \theta^{-1} \limsup_{t \rightarrow \infty} \frac{\E[\textrm{Loss}(t)]}{t}.
\end{equation}
Let $\textrm{Poi}$ denote a r.v. with a Poisson distribution, with mean $\frac{\mu_A}{\mu_S}$.  Then it follows from well-known insensitivity results for the Erlang loss model, see e.g. \citet{Davis.95,Sev57,Franken82}, that 
\begin{eqnarray}
\lim_{t \rightarrow \infty} \frac{\E[\textrm{Loss}(t)]}{t} &=& \mu_A \pr\big( \textrm{Poi} = n | \textrm{Poi} \leq n \big) \nonumber
\\&=& \mu_A \frac{ \frac{\exp(- \frac{\mu_A}{\mu_S} ) (\frac{\mu_A}{\mu_S})^n}{n!} }{ \sum_{k=0}^n \frac{\exp(- \frac{\mu_A}{\mu_S} ) (\frac{\mu_A}{\mu_S})^k}{k!}}. \label{abanboundeq1}
\end{eqnarray}
Although many bounds exist for this blocking probability, see e.g. \citet{Hariel.88,Hariel.10,Janssen08c}, to prove our first bound we provide a self-contained and very simple bound.  As our assumptions imply $\frac{\mu_A}{\mu_S} < n$, and it follows from \citet{Chen86} that the median of $\textrm{Poi}$ is at most $\lceil \frac{\mu_A}{\mu_S} \rceil \leq n$, we conclude that (\ref{abanboundeq1}) is at most
\begin{equation}\label{abanboundeq2}
2 \mu_A \frac{\exp(- \frac{\mu_A}{\mu_S} ) (\frac{\mu_A}{\mu_S})^n}{n!}.
\end{equation}
Applying Stirling's inequality, which implies $n! \geq  \sqrt{2 \pi n} (\frac{n}{e})^n$, we conclude that (\ref{abanboundeq1}) is at most
\begin{equation}\label{abanboundeq2}
\sqrt{\frac{2}{\pi}} \mu_A \frac{\exp(n - \frac{\mu_A}{\mu_S}) (1 - \frac{n - \frac{\mu_A}{\mu_S}}{n})^n}{\sqrt{n}}.
\end{equation}
Noting that $1 - x \leq e^{-x}$ for all $x \geq 0$, we conclude that (\ref{abanboundeq1}) is at most $\sqrt{\frac{2}{\pi}} \frac{\mu_A}{\sqrt{n}}.$  Combining with (\ref{abandoneq0b}), we conclude that 
\begin{eqnarray*}
\E[L^n_a(\infty)] &\leq& \sqrt{\frac{2}{\pi}} \theta^{-1} \frac{\mu_A}{\sqrt{n}}
\\&=& \sqrt{\frac{2}{\pi}} \theta^{-1} \frac{\mu_A}{n \mu_S} \frac{n \mu_S}{\sqrt{n}}
\\&=& \sqrt{\frac{2}{\pi}} \rho \frac{\mu_S}{\theta} \sqrt{n},
\end{eqnarray*}
completing the proof of the first bound.  For the second part, we note that since it is easily verified that $e^{-x} x^n$ is increasing on $[0,n]$, we may upper bound $\frac{\exp(- \frac{\mu_A}{\mu_S} ) (\frac{\mu_A}{\mu_S})^n}{n!}$ by $\frac{\exp(- \lceil \frac{\mu_A}{\mu_S} \rceil ) (\lceil \frac{\mu_A}{\mu_S} \rceil)^n}{n!}$.  Further supposing $\rho \in [\frac{3}{4}, 1 - \frac{4}{n}]$ (also implying $n \geq \lceil \frac{\mu_A}{\mu_S} \rceil + 2$) , then it follows from \citet{Glynn87} Proposition 2 that $\frac{\exp(- \lceil \frac{\mu_A}{\mu_S} \rceil ) (\lceil \frac{\mu_A}{\mu_S} \rceil)^n}{n!}$ is at most
\begin{eqnarray*}
\ &\ &\ (2 \pi \lceil \frac{\mu_A}{\mu_S} \rceil)^{-\frac{1}{2}} \exp\bigg( - \frac{(n - \lceil \frac{\mu_A}{\mu_S} \rceil)(n - \lceil \frac{\mu_A}{\mu_S} \rceil - 1)}{2 \lceil \frac{\mu_A}{\mu_S} \rceil} + \frac{(n - \lceil \frac{\mu_A}{\mu_S} \rceil)(n - \lceil \frac{\mu_A}{\mu_S} \rceil - 1)\big(2 (n - \lceil \frac{\mu_A}{\mu_S} \rceil) - 1 \big)}{12 (\lceil \frac{\mu_A}{\mu_S} \rceil)^2} \bigg)
\\&\leq& (2 \pi \lceil \frac{\mu_A}{\mu_S} \rceil)^{-\frac{1}{2}} \exp\Bigg( - \frac{(n - \lceil \frac{\mu_A}{\mu_S} \rceil)(n - \lceil \frac{\mu_A}{\mu_S} \rceil - 1)}{2 \lceil \frac{\mu_A}{\mu_S} \rceil} \bigg(1 -  \frac{n - \lceil \frac{\mu_A}{\mu_S} \rceil}{3 \lceil \frac{\mu_A}{\mu_S} \rceil} \bigg) \Bigg)
\\&\leq& (2 \pi \lceil \frac{\mu_A}{\mu_S} \rceil)^{-\frac{1}{2}} \exp\big( - \frac{8}{9} \frac{(n - \lceil \frac{\mu_A}{\mu_S} \rceil)(n - \lceil \frac{\mu_A}{\mu_S} \rceil - 1)}{2 \lceil \frac{\mu_A}{\mu_S} \rceil}  \big)
\\&\ &\ \ \ \textrm{since our assumptions imply}\ 1 -  \frac{n - \lceil \frac{\mu_A}{\mu_S} \rceil}{3 \lceil \frac{\mu_A}{\mu_S} \rceil} \geq \frac{8}{9}
\\&=& (2 \pi \lceil \frac{\mu_A}{\mu_S} \rceil)^{-\frac{1}{2}} \exp\big( - \frac{4}{9} \frac{(n - \lceil \frac{\mu_A}{\mu_S} \rceil)^2 - (n - \lceil \frac{\mu_A}{\mu_S} \rceil)}{\lceil \frac{\mu_A}{\mu_S} \rceil}  \big)
\\&\leq& (2 \pi \lceil \frac{\mu_A}{\mu_S} \rceil)^{-\frac{1}{2}} \exp\big( - \frac{4}{9} \frac{(n-\lceil \frac{\mu_A}{\mu_S} \rceil)^2}{n} + \frac{4}{9} \times \frac{4}{3} \times \frac{n - \lceil \frac{\mu_A}{\mu_S} \rceil}{n} \big)
\\&\leq& (2 \pi \lceil \frac{\mu_A}{\mu_S} \rceil)^{-\frac{1}{2}} \exp\bigg( - \frac{1}{4} \big( \sqrt{n}(1 - \rho) \big)^2 + (1 - \rho) \bigg),
\end{eqnarray*}
the final inequality following from the fact that $\frac{(n-\lceil \frac{\mu_A}{\mu_S} \rceil)^2}{n} = \big( \sqrt{n}(1 - \rho) \big)^2 \times (\frac{n - \lceil \frac{\mu_A}{\mu_S} \rceil}{n - \frac{\mu_A}{\mu_S}})^2$, and our assumption that $\rho \leq 1 - \frac{4}{n}$ implies that $(\frac{n - \lceil \frac{\mu_A}{\mu_S} \rceil}{n - \frac{\mu_A}{\mu_S}})^2 \geq \frac{9}{16}$.  We conclude that (\ref{abanboundeq1}) is at most $\sqrt{\frac{2}{\pi}}\sqrt{\mu_A \mu_S} \exp\bigg( - \frac{1}{4} \big( \sqrt{n}(1 - \rho) \big)^2 + (1 - \rho) \bigg).$
Combining with (\ref{abandoneq0b}), we conclude that 
\begin{eqnarray*}
\E[L^n_a(\infty)] &\leq& \sqrt{\frac{2}{\pi}} \theta^{-1} \sqrt{\mu_A \mu_S} \exp\bigg( - \frac{1}{4} \big( \sqrt{n}(1 - \rho) \big)^2 + (1 - \rho) \bigg)
\\&=& \sqrt{\frac{2}{\pi}} \sqrt{\rho} \frac{\mu_S}{\theta} \sqrt{n} \exp\bigg( - \frac{1}{4} \big( \sqrt{n}(1 - \rho) \big)^2 + (1 - \rho) \bigg)
\\&\leq& \sqrt{\frac{2}{\pi}} \frac{\mu_S}{\theta} \sqrt{n} \exp\bigg( - \frac{1}{4} \big( \sqrt{n}(1 - \rho) \big)^2 + (1 - \rho) \bigg).
\end{eqnarray*}
Combining the above with the fact that $\rho \geq \frac{3}{4}$ implies $\sqrt{\frac{2}{\pi}} \exp(1 - \rho) \leq 2$ completes the proof.  $\halmos$
\endproof
\subsubsection{Proof of Theorem\ \ref{sspdnoprefactorc}.}\label{noprefactorsec3}
We begin by stating the aforementioned equation for the steady-state expected work-in-system studied in several past works (see e.g. \citet{hokstad85, scully2020gittins, grosof2021finite, wang2021zero, hong2021sharp}), which arises when applying the Lyapunov drift method to the square of the total work in system.  We impose the technical condition that $S$ has absolutely continuous distribution function to be consistent with the assumptions and arguments of \citet{hokstad85}, although as noted in \citet{hokstad85} the result holds in greater generality.
\begin{theorem}[\citet{hokstad85, scully2020gittins, grosof2021finite}]\label{lyapquad}
Consider an $M/GI/n$ queue ${\mathcal Q}^n$ with Markovian inter-arrival times having the same distribution as r.v. $A$, service times having the same distribution as r.v. $S$ satisfying $\E[S^2] < \infty$, such that the c.d.f. of $S$ is absolutely continuous.  Suppose also that $\mu_A < n \mu_S,$ and that $Q(\infty), \textrm{Work}(\infty),$ and $\overline{W}_{\textrm{service}}(\infty)$ exist.  Then
$$
\E\big[\textrm{Work}(\infty)\big] = 
\bigg(\frac{1}{2} \mu_A \E[S^2] + \E\big[ \big(n - \textrm{Num}_{\textrm{service}}(\infty) \big) \textrm{Work}_{\textrm{service}}(\infty) \big] \bigg) \times \frac{1}{n(1-\rho)}.$$
\end{theorem}
With Theorem\ \ref{lyapquad} in hand, our result will follow from a straightforward rearranging of terms.
\proof{Proof of Theorem\ \ref{sspdnoprefactorc}:}
By the basic properties of a $GI/GI/n$ queue, and decomposing the work in system into that in service and that waiting in queue, note that 
$$\E\big[\textrm{Work}(\infty)\big] = \E\big[ \textrm{Work}_{\textrm{service}}(\infty) \big] + \E[S] \E\big[L(\infty)\big].$$
Combining with Theorem\ \ref{lyapquad} and some straightforward algebra then yields that $\E\big[L(\infty)\big]$ equals
\begin{eqnarray*}
\ &\ &\ \frac{ \frac{1}{2} \mu_A \E[S^2] }{n (1 - \rho) \E[S]} + \frac{\E\big[ \big(n - \textrm{Num}_{\textrm{service}}(\infty) \big) \textrm{Work}_{\textrm{service}}(\infty) \big] - n(1-\rho) \E\big[\textrm{Work}_{\textrm{service}}(\infty) \big]}{n (1 - \rho) \E[S]}
\\&=& \frac{1}{2}\E[(S \mu_S)^2] \frac{\rho}{1-\rho} + \frac{\E\big[ - \textrm{Num}_{\textrm{service}}(\infty) \times \textrm{Work}_{\textrm{service}}(\infty) \big] + n \times \rho \times \E\big[\textrm{Work}_{\textrm{service}}(\infty) \big]}{n (1 - \rho) \E[S]}
\\&=& \frac{1}{2}\E[(S \mu_S)^2] \frac{\rho}{1-\rho} + \frac{\E\big[\textrm{Num}_{\textrm{service}}(\infty)\big] \times \E\big[\textrm{Work}_{\textrm{service}}(\infty) \big]
- \E\big[\textrm{Num}_{\textrm{service}}(\infty) \times \textrm{Work}_{\textrm{service}}(\infty) \big]}{n (1 - \rho) \E[S]},
\end{eqnarray*}
the final equality using the fact that $\E\big[\textrm{Num}_{\textrm{service}}(\infty)\big] = n \rho$ (a fact also used in our proof of Theorem\ \ref{sspdnoprefactor}).  Combining the above completes the proof.  $\halmos$
\endproof

\subsection{Proof of Lemma\ \ref{2partbound}.}\label{2partboundsec}
\proof{Proof of Lemma\ \ref{2partbound}:}
\begin{eqnarray*}
\ &\ &\ \pr\bigg( \sup_{t \geq 0} \big( A_e(t) - \sum_{i=1}^{n} N_{e,i}(t) \big) \geq x \bigg)\nonumber
\\&=& \pr\Bigg( \sup_{t \geq 0} \bigg( \big( A_e(t) - \frac{1}{2}(n + \mu_A) t \big)
+ \big( \frac{1}{2}(n + \mu_A) t - \sum_{i=1}^{n} N_{e,i}(t) \big) \bigg) \geq x \Bigg)
\\&=& \pr\Bigg( \sup_{t \geq 0} \Bigg( \bigg( A_e(t) - \big( \mu_A +\frac{1}{2}(n - \mu_A) \big) t \bigg)
+ \bigg( \big(n - \frac{1}{2}(n - \mu_A) \big) t - \sum_{i=1}^{n} N_{e,i}(t) \bigg) \Bigg) \geq x \Bigg)
\\&=& \pr\Bigg( \sup_{t \geq 0} \Bigg( \bigg( A_e(t) - \mu_A t - \frac{1}{2}(n - \mu_A) t \bigg)
+ \bigg( n t - \sum_{i=1}^{n} N_{e,i}(t)  - \frac{1}{2}(n - \mu_A) t \bigg) \Bigg) \geq x \Bigg)
\\&\leq& \pr\bigg( \sup_{t \geq 0} \big( A_e(t) - \mu_A t - \frac{1}{2}(n - \mu_A) t \big) \geq \frac{x}{2} \bigg)
\\&+& \pr\bigg( \sup_{t \geq 0} \big( n t - \sum_{i=1}^{n} N_{e,i}(t) - \frac{1}{2}(n - \mu_A) t \big) \geq \frac{x}{2} \Bigg).
\end{eqnarray*}
Using (\ref{chaaaanges}) to bound $A_e(t)$ by $A_o(t) +  1$ completes the proof.  $\halmos$
\endproof

\subsection{Proof of Lemma\ \ref{adiscbound}.}\label{adiscboundsec}

\proof{Proof of Lemma\ \ref{adiscbound}:}
As $\lbrace A_o(t) - \mu_{A} t - \nu t, t \geq 0 \rbrace$ jumps up only at times $\lbrace \sum_{i=1}^k A_{i}, k \geq 1 \rbrace$ and at all other times drifts downward at linear rate $-(\mu_{A} + \nu)$, we conclude that we may examine the relevant supremum only at times $\lbrace \sum_{i=1}^k A_{i}, k \geq 0 \rbrace$, from which it follows that (\ref{tobounda1}) equals
\begin{equation}\label{tobounda2bb}
\pr\Bigg( \sup_{k \geq 0} \bigg( k - \sum_{i=1}^k \big( (\mu_{A} + \nu) A_{i} \big) \bigg) \geq x \Bigg).
\end{equation}
Further observing that
\begin{eqnarray*}
 k - (\mu_{A} + \nu) \sum_{i=1}^k A_{i} &=& (1 + \frac{\nu}{\mu_{A}}) k - (\mu_{A} + \nu) \sum_{i=1}^k A_{i} - \frac{\nu}{\mu_{A}} k
\\&=& (1 + \frac{\nu}{\mu_{A}}) \big( k - \mu_{A} \sum_{i=1}^k A_{i} - \frac{\nu}{\mu_{A} + \nu} k \big)
\\&=& (1 + \frac{\nu}{\mu_{A}}) \big( \frac{\mu_A}{\mu_A + \nu} k - \sum_{i=1}^k (\mu_{A} A_{i})\big)
\end{eqnarray*}
completes the proof.  $\halmos$
\endproof

\subsection{Proof of Lemma\ \ref{kingmanmart2}.}\label{kingmanmart2sec}
\proof{Proof of Lemma\ \ref{kingmanmart2}:}
First, note that $\E[\exp\big(\theta(c - Z_1)\big)]$ exists for all $\theta > 0$ since $Z_1$ is non-negative and $c$ is a fixed constant.  We consider two cases, since the desired scaling is somewhat different as $c \downarrow 0$.  First, suppose $c \geq \frac{1}{2}$ (and thus $c \in [\frac{1}{2},1)$).  In this case, we argue that if $\E[\exp\big(\theta(c - Z_1)\big)] > 1$, then it must hold that $\theta > \frac{1 - c}{2.75 \E[Z_1^2]}.$  Thus suppose $c \geq \frac{1}{2}$, and $\E[\exp\big(\theta(c - Z_1)\big)] > 1$ for some $\theta \in (0,1)$ (here it suffices to consider $\theta \in (0,1)$ as $c \geq \frac{1}{2}$ implies $\frac{1 - c}{2.75 \E[Z_1^2]} < 1$).  It follows from the fact that $\theta c \in(0,1)$, and a straightforward calculus exercise (the details of which we omit), that (1) : $\exp(\theta c) < 1 + \theta c + .75 \theta^2 c^2$; and (2) : $\exp(-\theta Z_1) < \frac{1}{1 + \theta Z_1}.$  (2) implies   
\begin{eqnarray*}
\E[\exp(- \theta Z_1)] &<& \E[\frac{1}{1+\theta Z_1}]
\\&=& \E[1 - \theta Z_1 + \frac{(\theta Z_1)^2}{1 + \theta Z_1}]
\\&\leq& 1 - \theta \E[Z_1] + \theta^2 \E[Z_1^2],
\end{eqnarray*}
where we note that a related logic appears in \citet{G16b}.  Thus $\E[\exp\big(\theta(c - Z_1)\big)] > 1$, combined with our other assumptions, implies
$$\big(1 + \theta c + .75 \theta^2 c^2 \big) \times \big(1 - \theta \E[Z_1] + \theta^2 \E[Z_1^2]\big) > 1,$$
which by some straightforward algebra is equivalent to
$$(c - \E[Z_1]) \theta + \big(  \E[Z_1^2] - c \E[Z_1] + .75 c^2 \big) \theta^2 + \big( c \E[Z_1^2] - .75 c^2 \E[Z_1] \big) \theta^3 + .75 c^2 \E[Z_1^2] \theta^4 > 0.$$
Combining with the fact that $c,\theta \in (0,1)$ and $\E[Z_1] = 1$, some additional straightforward algebra further implies that
$$(c - 1) \theta + 2.75 \E[Z_1^2] \theta^2 > 0,$$
itself implying that $\theta > \frac{ 1 - c }{2.75 \E[Z_1^2]}$. 
\\\\Next, we argue that if $c \in (0,\frac{1}{2})$ and $\E[\exp\big(\theta(c - Z_1)\big)] > 1$, then it must hold that $\theta > \big( 11 c \E[Z_1^2]\big)^{-1}.$  Thus suppose $c \in (0,\frac{1}{2})$, and $\E[\exp\big(\theta(c - Z_1)\big)] > 1$ for some $\theta \in \big( 0, (2 c)^{-1} \big)$.  Here it suffices to consider $\theta \in \big(0, (2c)^{-1}\big)$ since $\big( 11 c \E[Z_1^2]\big)^{-1} < (2c)^{-1}$.  It follows from the non-negativity of $Z_1$ and fact that $c \in (0,\frac{1}{2})$ that for all $\theta > 0$, w.p.1 $\theta (c - Z_1) < 2 \theta c (\frac{1}{2} - Z_1),$ and thus $\E[\exp\big(\theta(c - Z_1)\big)] < \E[\exp\big(2 \theta c (\frac{1}{2} - Z_1)\big)]$.  Thus $\E[\exp\big(\theta(c - Z_1)\big)] > 1$ for some $\theta \in \big( 0, (2 c)^{-1} \big)$ implies $\E[\exp\big(2 \theta c(\frac{1}{2} - Z_1)\big)] > 1$ for some $\theta \in \big( 0, (2 c)^{-1} \big)$.  As $\theta \in \big(0, (2 c)^{-1} \big)$ implies $2 \theta c \in (0,1)$, it then follows from an argument nearly identical to that used in our previous analysis of the $c \geq \frac{1}{2}$ case (and the details of which we omit) that $2 \theta c > \frac{ 1 - \frac{1}{2} }{2.75 \E[Z_1^2]}$, and thus $\theta > \big( 11 c \E[Z_1^2] \big)^{-1}$, which completes the argument in this case.
\\\\Combining the above with the fact that (1) : $\frac{ 1 - c }{2.75 \E[Z_1^2]} > \frac{ 1 - c }{11 c \E[Z_1^2]}$ for $c \in [\frac{1}{2},1)$; and (2) : $\big( 11 c \E[Z_1^2] \big)^{-1} > \frac{ 1 - c }{11 c \E[Z_1^2]}$ for $c \in (0,\frac{1}{2}),$ completes the proof.
$\halmos$
\endproof

\subsection{Proof of Lemma \ref{TailBound}.}\label{TailBoundsec}
		\proof{Proof of Lemma\ \ref{TailBound}:}
			Note that for $\lambda > 0$, $\pr\bigg(\sup_{t\ge 0} \big(\phi (t) - \nu t \big) \geq \lambda \bigg)$ equals
\begin{eqnarray}
\ &\ &\ \pr \Bigg( \bigg( \bigcup_{k=0}^\infty \left\{\phi(t) - \nu t \geq \lambda  \textrm{  for some $t\in [2^k, 2^{k+1}]$}\right\} \bigg) \bigcup  \bigg\{\phi(t) - \nu t \geq \lambda  \textrm{  for some $t\in [0,1]$}\bigg\}  \Bigg) \nonumber
\\&\ &\ \ \ \leq\ \ \ \sum_{k=0}^\infty \pr \bigg(\sup_{t\in [2^k,2^{k+1}]} \big(\phi(t) - \nu t\big) \geq \lambda \bigg) \label{TailBoundeq1a}
\\&\ &\ \ \ \ \ \ \ \ \ \ \ +\ \ \ \pr \bigg(\sup_{t\in [0,1]} \big(\phi(t) - \nu t\big) \geq \lambda\bigg). \label{TailBoundeq1b}
\end{eqnarray}
We now bound (\ref{TailBoundeq1a}), and proceed by bounding (for each $k \geq 0$)
\begin{equation}\label{eachk}
\pr \bigg(\sup_{t\in [2^k,2^{k+1}]} \big(\phi(t) - \nu t\big) \geq \lambda \bigg).
\end{equation}
Since $t\in [2^k,2^{k+1}]$ implies $\nu t \geq \nu 2^k$, we conclude that (\ref{eachk}) is at most $
\pr \bigg(\sup_{t\in [2^k,2^{k+1}]} \phi(t) \geq \lambda + \nu 2^k \bigg),$
which by adding and subtracting $\phi(2^k)$, and applying stationary increments and a union bound, is at most
\begin{eqnarray}
\ &\ &\ \pr\bigg(\big(\sup_{t\in [2^k, 2^{k+1}]} \phi(t) - \phi(2^k) \big) + \phi(2^k) \geq \lambda + \nu 2^k\bigg) \nonumber
\\&\ &\ \ \ \leq\ \ \ \pr\bigg(\sup_{t\in [2^k, 2^{k+1}]} \phi(t) - \phi(2^k) \geq \frac{1}{2} (\lambda + \nu 2^k) \bigg)
+ \pr\bigg(\phi(2^k) \geq \frac{1}{2} (\lambda + \nu 2^k) \bigg) \nonumber
\\&\ &\ \ \ =\ \ \ \pr\bigg(\sup_{t\in [0, 2^{k}]} \phi(t) \geq \frac{1}{2} (\lambda + \nu 2^k) \bigg)
+ \pr\bigg(\phi(2^k) \geq \frac{1}{2} (\lambda + \nu 2^k) \bigg) \nonumber
\\&\ &\ \ \ \leq\ \ \ 2 \pr\bigg(\sup_{t\in [0, 2^{k}]} \phi(t) \geq \frac{1}{2} (\lambda + \nu 2^k) \bigg). \label{TailBoundeq2}
\end{eqnarray}
We proceed to bound (\ref{TailBoundeq2}) by breaking the supremum into two parts, one part taken over integer points, one part taken over intervals of length one corresponding to the regions between these integer points.  In particular, the assumptions of the lemma, combined with a union bound and stationary increments, ensure that 
\begin{eqnarray}
\ &\ &\ \pr \bigg(\sup_{t\in [0, 2^k]} \phi(t) \geq \frac{1}{2} (\lambda + \nu 2^k) \bigg) \nonumber
\\&\ &\ \ \ \leq\ \ \ \pr\bigg(\sup_{j \in \{0,..., 2^k\}} \phi(j) + \sup_{\substack{j\in \{0,...,2^k-1\}\\t \in [0,1]}} \big( \phi(j+t) - \phi(j) \big) \geq \frac{1}{2} (\lambda + \nu 2^k) \bigg) \nonumber
\\&\ &\ \ \ \leq\ \ \ \pr\bigg(\sup_{j \in \{0,..., 2^k\}} \phi(j) \geq \frac{1}{4} (\lambda + \nu 2^k) \bigg) 
+ 2^k \pr\bigg(\sup_{t\in [0,1]}\phi(t) \geq \frac{1}{4} (\lambda + \nu 2^k) \bigg) \nonumber
\\&\ &\ \ \ \leq\ \ \ \frac{H_14^{r_1}2^{ks}}{(\lambda + \nu 2^k)^{r_1}} + \frac{H_2 4^{r_2} 2^k }{(\lambda+\nu 2^k)^{r_2}}, \label{TailBoundeq3}
			\end{eqnarray}
where the final inequality is applicable since $\lambda \geq 4 Z$ implies $\frac{1}{4} (\lambda + \nu 2^k) \geq Z$, in which case the inequality follows from our assumptions.  Combining (\ref{TailBoundeq2}) and (\ref{TailBoundeq3}), we conclude that (\ref{TailBoundeq1a}) is at most
\begin{equation}\label{tttail1eq}
2 \sum_{k=0}^{\infty} \frac{H_14^{r_1}2^{ks}}{(\lambda + \nu 2^k)^{r_1}} + 2 \sum_{k=0}^{\infty} \frac{H_2 4^{r_2} 2^k }{(\lambda+\nu 2^k)^{r_2}}.
\end{equation}
We now treat two cases.  First, suppose $\lambda > \nu$.  Then (\ref{tttail1eq}) is at most 
\begin{eqnarray*}
\ &\ &\ 2 H_1 4^{r_1} \sum_{k=0}^{\lceil \log_2(\frac{\lambda}{\nu}) \rceil - 1} \frac{2^{k s}}{\lambda^{r_1}}
\\&\ &\ \ \ \ +\ \ \ 2 H_2 4^{r_2} \sum_{k=0}^{\lceil \log_2(\frac{\lambda}{\nu}) \rceil - 1} \frac{2^{k}}{\lambda^{r_2}}
\\&\ &\ \ \ \ +\ \ \ 2 H_1 4^{r_1} \sum_{k = \lceil \log_2(\frac{\lambda}{\nu}) \rceil}^{\infty} \frac{2^{-(r_1 - s) k }}{\nu^{r_1}}
\\&\ &\ \ \ \ +\ \ \ 2 H_2 4^{r_2} \sum_{k = \lceil \log_2(\frac{\lambda}{\nu}) \rceil}^{\infty} \frac{2^{-(r_2 - 1) k }}{\nu^{r_2}}
\\&=&\ 2 H_1 4^{r_1} \lambda^{- r_1} \frac{2^{\lceil \log_2(\frac{\lambda}{\nu}) \rceil s} - 1}{2^s - 1}
\\&\ &\ \ \ \ +\ \ \ 2 H_2 4^{r_2} \lambda^{-r_2} (2^{\lceil \log_2(\frac{\lambda}{\nu}) \rceil} - 1)
\\&\ &\ \ \ \ +\ \ \ 2 H_1 4^{r_1} \nu^{-r_1} \frac{2^{-(r_1 - s)\lceil \log_2(\frac{\lambda}{\nu}) \rceil}}{1 - 2^{-(r_1 - s)}}
\\&\ &\ \ \ \ +\ \ \ 2 H_2 4^{r_2} \nu^{-r_2} \frac{2^{-(r_2 - 1)\lceil \log_2(\frac{\lambda}{\nu}) \rceil}}{1 - 2^{-(r_2 - 1)}}
\\&\leq&\ 4 H_1 4^{r_1} \lambda^{- r_1} (\frac{\lambda}{\nu})^s
\\&\ &\ \ \ \ +\ \ \ 4 H_2 4^{r_2} \lambda^{-r_2} \frac{\lambda}{\nu}
\\&\ &\ \ \ \ +\ \ \ 2 H_1 \big(1 - 2^{-(r_1 - s)}\big)^{-1} 4^{r_1} \nu^{-r_1} (\frac{\lambda}{\nu})^{-(r_1 - s)} 
\\&\ &\ \ \ \ +\ \ \ 2 H_2 \big(1 - 2^{-(r_2 - 1)}\big)^{-1} 4^{r_2} \nu^{-r_2} (\frac{\lambda}{\nu})^{-(r_2 - 1)},
\end{eqnarray*}
with the first line of the final inequality following from the fact that $2^{\lceil \log_2(\frac{\lambda}{\nu}) \rceil s} - 1 \leq 2^s (\frac{\lambda}{\nu})^s$ and $2^s - 1 \geq 2^{s - 1}$.  Combining with the fact that $r_2 > 2$ implies $\big(1 - 2^{-(r_2 - 1)}\big)^{-1} \leq 2$, we conclude that if $\lambda > \nu$, then (\ref{tttail1eq}) is at most
\begin{eqnarray}
\ &\ &\ 6 H_1 \big(1 - 2^{-(r_1 - s)}\big)^{-1} 4^{r_1} \lambda^{- (r_1 - s)} \nu^{-s} \label{tttbig1aaa}
\\&\ &\ \ \ \ \ \ +\ \ \ 8 H_2 4^{r_2} \lambda^{- (r_2 - 1)} \nu^{-1}. \nonumber
\end{eqnarray}
Combining with the fact that $\lambda > \nu > 0$ and $r_2 > 2$ implies $\lambda^{- (r_2 - 1)} \nu^{-1} \leq ( \lambda \nu )^{-\frac{r_2}{2}}$,
we conclude that if $\lambda > \nu$, then (\ref{tttail1eq}) is at most
\begin{eqnarray}
\ &\ &\ 6 H_1 \big(1 - 2^{-(r_1 - s)}\big)^{-1} 4^{r_1} \lambda^{- (r_1 - s)} \nu^{-s} \label{tttbig1}
\\&\ &\ \ \ \ \ \ +\ \ \ 8 H_2 4^{r_2} (\lambda \nu)^{- \frac{r_2}{2}}. \nonumber
\end{eqnarray}
Alternatively, suppose $\lambda \leq \nu$.  Then (\ref{tttail1eq}) is at most 
\begin{eqnarray}
\ &\ & 2 H_1 4^{r_1} \sum_{k = 0}^{\infty} \frac{2^{-(r_1 - s) k }}{\nu^{r_1}} \nonumber
\\&\ &\ \ \ \ \ \ +\ \ \ 2 H_2 4^{r_2} \sum_{k = 0}^{\infty} \frac{2^{-(r_2 - 1) k }}{\nu^{r_2}} \nonumber
\\&\leq&\ 2 H_1 4^{r_1} \nu^{-r_1} \big(1 - 2^{-(r_1 - s)}\big)^{-1} \nonumber
\\&\ &\ \ \ \ \ \ +\ \ \ 4 H_2 4^{r_2} \nu^{-r_2} \nonumber
\\&\leq&\ 2 H_1 \big(1 - 2^{-(r_1 - s)}\big)^{-1} 4^{r_1} \lambda^{-(r_1 - s)} \nu^{-s} \label{tttbig2}
\\&\ &\ \ \ \ \ \ +\ \ \ 4 H_2 4^{r_2} (\lambda \nu)^{-\frac{r_2}{2}}, \nonumber
\end{eqnarray}
the final inequality following from the fact that $\nu \geq \lambda, r_1 > s, r_2 > 2$ implies $\nu^{-r_1} \leq \lambda^{-(r_1 - s)} \nu^{-s}$, and $\nu^{-r_2} \leq (\lambda \nu)^{- \frac{r_2}{2}}$.  Next, we claim that $\big(1 - 2^{-(r_1 - s)}\big)^{-1} \leq 2(1 + \frac{1}{r_1 - s})$.  Indeed, first, suppose $r_1 - s < 1$.  In this case, as it is easily verified that $1 - 2^{-z} \geq \frac{z}{2}$ for all $z \in (0,1)$, the result follows.  Alternatively, if $r_1 - s \geq 1$, then $\big(1 - 2^{-(r_1 - s)}\big)^{-1} \leq 2$, completing the proof.  Combining with (\ref{tttbig1}) and (\ref{tttbig2}), and our assumptions, it follows that in all cases (\ref{tttail1eq}), and hence (\ref{TailBoundeq1a}), is at most
\begin{equation}\label{tttbig3}
12 H_1 (1 + \frac{1}{r_1 - s})  4^{r_1} \lambda^{-(r_1-s)} \nu^{-s} + 8 H_2 4^{r_2} (\lambda \nu)^{ - \frac{r_2}{2} }.
\end{equation}
We next bound (\ref{TailBoundeq1b}).  First, suppose $\lambda \geq \nu$.  Then our assumptions (applied with $t_0 = 1$) imply that (\ref{TailBoundeq1b}) is at most 
\begin{equation}\label{aux333}
H_2 \lambda^{-r_2} \leq H_2 (\lambda \nu)^{-\frac{r_2}{2}}.
\end{equation}
Alternatively, suppose that $\lambda < \nu$.  Then applying our assumptions with $t_0 = \frac{\lambda}{\nu}$, along with a union bound, we conclude that
\begin{eqnarray}
\ &\ &\ \pr \bigg(\sup_{t\in [0,1]} \big(\phi(t) - \nu t\big) \geq \lambda\bigg)\nonumber
\\&\leq&\ \pr\bigg(\sup_{t\in [0, \frac{\lambda}{\nu}]} \big(\phi(t) - \nu t\big) \geq \lambda\bigg)\label{newideaeq1}
\\&\ &\ \ \ \ \ \ +\ \ \ \pr\bigg(\sup_{t\in [\frac{\lambda}{\nu},1]} \big(\phi(t) - \nu t\big) \geq \lambda\bigg). \label{newideaeq2}
\end{eqnarray}
It follows from our assumptions that (\ref{newideaeq1}) is at most
\begin{equation}\label{newideaeq3}
 \pr\Big(\sup_{t\in [0, \frac{\lambda}{\nu}]} \phi(t) \geq \lambda\Big)\ \ \ \leq\ \ \ H_2 (\frac{\lambda}{\nu})^{\frac{r_2}{2}} \lambda^{-r_2}\ \ \ =\ \ \ H_2 (\lambda \nu)^{-\frac{r_2}{2}}.
\end{equation}
We next bound (\ref{newideaeq2}), which by stationary increments, a union bound, and our assumptions is at most
\begin{eqnarray}
\ &\ &\ \pr\bigg(\sup_{t\in [\frac{\lambda}{\nu},1]} \big( \phi(\frac{\lambda}{\nu}) + \phi(t) - \phi(\frac{\lambda}{\nu}) - \nu(t - \frac{\lambda}{\nu}) \big) \geq 2 \lambda \bigg) \nonumber
\\&\ &\ \ \ \leq\ \ \ \pr\big( \phi(\frac{\lambda}{\nu}) \geq \frac{1}{2} \lambda \big) \nonumber
\\&\ &\ \ \ \ \ \ +\ \ \ \pr\bigg( \sup_{t\in [\frac{\lambda}{\nu},1]} \big( \phi(t) - \phi(\frac{\lambda}{\nu}) - \nu(t - \frac{\lambda}{\nu}) \big) \geq \frac{3}{2} \lambda \bigg) \nonumber
\\&\ &\ \ \ =\ \ \ \pr\big( \phi(\frac{\lambda}{\nu}) \geq \frac{1}{2} \lambda \big) \nonumber
\\&\ &\ \ \ \ \ \ +\ \ \ \pr\Big( \sup_{s \in [0,1 - \frac{\lambda}{\nu}]} \big( \phi(s + \frac{\lambda}{\nu}) - \phi(\frac{\lambda}{\nu}) - \nu s \big) \geq \frac{3}{2} \lambda \Big) \nonumber
\\&\ &\ \ \ \leq\ \ \ \pr\Big(\sup_{t\in [0, \frac{\lambda}{\nu}]} \phi(t) \geq \frac{1}{2} \lambda\Big) \nonumber
\\&\ &\ \ \ \ \ \ +\ \ \ \pr\Big( \sup_{s \in [0,1 - \frac{\lambda}{\nu}]} \big( \phi(s + \frac{\lambda}{\nu}) - \phi(\frac{\lambda}{\nu}) - \nu s \big) \geq \frac{3}{2} \lambda \Big) \nonumber
\\&\ &\ \ \ \leq\ \ \ H_2 (\frac{\lambda}{\nu})^{\frac{r_2}{2}} (\frac{\lambda}{2})^{-r_2} \nonumber
\\&\ &\ \ \ \ \ \ +\ \ \ \pr\Big( \sup_{s \in [0,1 - \frac{\lambda}{\nu}]} \big( \phi(s) - \nu s \big) \geq \frac{3}{2} \lambda \Big) \nonumber
\\&\ &\ \ \ \leq\ \ \ 2^{r_2} H_2 (\lambda \nu)^{- \frac{r_2}{2}} \label{newideaeq4}
\\&\ &\ \ \ \ \ \ +\ \ \ \pr\Big( \sup_{t\in [0,1]} \big( \phi(t) - \nu t \big) \geq \frac{3}{2} \lambda \Big). \label{newideaeq5}
\end{eqnarray}
Let us define $f(z) \stackrel{\Delta}{=} \pr\bigg( \sup_{t\in [0,1]} \big( \phi(t) - \nu t \big) \geq z \bigg)$.  Then using (\ref{newideaeq3}) to bound (\ref{newideaeq1}), and (\ref{newideaeq4}) - (\ref{newideaeq5}) to bound (\ref{newideaeq2}), we conclude that for all $z \in [2 Z, \nu)$, 
\begin{equation}\label{neqideaeq6}
f(z) \leq 2^{r_2 + 1} H_2 (z \nu)^{-\frac{r_2}{2}} + f(\frac{3}{2} z).
\end{equation}
Let $j^* \stackrel{\Delta}{=} \sup\lbrace j \in Z^+\ :\ (\frac{3}{2})^j \lambda < \nu \rbrace$.  Then it follows from (\ref{neqideaeq6}) that for all $j \in [0,j^*]$, 
\begin{equation}\label{newideaeq7}
f\big( (\frac{3}{2})^j \lambda\big) \leq 2^{r_2 + 1} H_2 (\lambda \nu)^{-\frac{r_2}{2}} \big( (\frac{3}{2})^{\frac{r_2}{2}} \big)^{-j} + f\big( (\frac{3}{2})^{j+1} \lambda\big).
\end{equation}
Combining (\ref{newideaeq7}) with a straightforward induction and our assumptions, and noting that $f$ is a non-increasing function, we conclude that for all $\lambda \in [2 Z, \nu)$,
\ \\\\
\begin{eqnarray}
f(\lambda) &\leq& 2^{r_2 + 1} H_2 (\lambda \nu)^{-\frac{r_2}{2}} \sum_{j=0}^{j^*} \big( (\frac{3}{2})^{\frac{r_2}{2}} \big)^{-j} + f(\nu) \nonumber
\\&\leq& 2^{r_2 + 1} H_2 (\lambda \nu)^{-\frac{r_2}{2}} \sum_{j=0}^{\infty} (\frac{3}{2})^{-j} + H_2 \nu^{-r_2} \nonumber
\\&\leq& 2^{r_2 + 3} H_2 (\lambda \nu)^{-\frac{r_2}{2}} + H_2 \nu^{-r_2} \nonumber
\\&\leq& 2^{r_2 + 4} H_2 (\lambda \nu)^{-\frac{r_2}{2}}, \label{newideaeq8}
\end{eqnarray}
the final inequality following from the fact that by assumption $\nu > \lambda$ and thus $\nu^{-r_2} \leq (\lambda \nu)^{-\frac{r_2}{2}}$ by the same logic as in (\ref{aux333}).  Thus using (\ref{aux333}) to bound (\ref{TailBoundeq1b}) in the case $\lambda \geq \nu$, and (\ref{newideaeq8}) to bound (\ref{TailBoundeq1b}) in the case $\lambda < \nu$, we conclude that in all cases (\ref{TailBoundeq1b}) is at most 
\begin{equation}\label{newideaeq9}
2^{r_2 + 4} H_2 (\lambda \nu)^{-\frac{r_2}{2}}.
\end{equation}
Using (\ref{tttbig3}) to bound (\ref{TailBoundeq1a}), and (\ref{newideaeq9}) to bound (\ref{TailBoundeq1b}), combined with some straightforward algebra, demonstrates that for all $\lambda \geq 4 Z$,
$
\pr\bigg(\sup_{t\ge 0} \big(\phi (t) - \nu t \big) \geq \lambda \bigg)$ is at most
\begin{eqnarray*}
\ &\ &\ 12 H_1 (1 + \frac{1}{r_1 - s})  4^{r_1} \lambda^{-(r_1-s)} \nu^{-s}
\\&+& 8 H_2 4^{r_2} (\lambda \nu)^{ - \frac{r_2}{2} }
\\&+& 2^{r_2 + 4} H_2 (\lambda \nu)^{-\frac{r_2}{2}}
\\&\leq& 12 H_1 (1 + \frac{1}{r_1 - s})  4^{r_1} \lambda^{-(r_1-s)} \nu^{-s}
\\&+& 12 H_2 4^{r_2} (\lambda \nu)^{ - \frac{r_2}{2} }
\\&\leq& 12 (1 + \frac{1}{r_1 - s})\big( H_1 4^{r_1} \lambda^{-(r_1-s)} \nu^{-s} + H_2 4^{r_2} (\lambda \nu)^{ - \frac{r_2}{2} } \big),
\end{eqnarray*}
completing the proof.  $\halmos$
\endproof

\subsection{Proof of Lemmas\ \ref{maximal2} and \ref{maximal3}.}\label{maximalssec}
In this section we prove Lemmas\ \ref{maximal2} and \ref{maximal3}.  

\subsubsection{A maximal inequality we will use in the proof of both Lemmas\ \ref{maximal2} and \ref{maximal3}.}
As mentioned previously, our proofs of both Lemmas\ \ref{maximal2} and \ref{maximal3} will rely heavily on a maximal inequality of \citet{LS77}.  We begin by stating (a variant of) the relevant maximal inequality of \citet{LS77} which we will use in both proofs.

\begin{lemma}[\citet{LS77} Theorem 2]\label{maximal1}
Let $\lbrace X_l, 1 \leq l \leq L \rbrace$ be a completely general sequence of r.v.s.  Suppose that for some fixed $\gamma > 1$, $\nu \geq \gamma$, and $C > 0$ the following condition holds:
\begin{enumerate}[(i)]
\item For all $x > 0$ and non-negative integers $1 \leq i \leq j \leq L$, 
$$
\pr\big(|\sum_{k=i}^j X_k| \geq x\big) \leq \big( C (j - i + 1) \big)^{\gamma} x^{-\nu}.
$$
\end{enumerate}
Then it must also hold that
$$
\pr\big(\max_{i \in \lbrace 1,\ldots,L \rbrace} |\sum_{k=1}^i X_k| \geq x\big) \leq 2.25 \times 2^{\gamma} \times (2 \frac{\nu + 1}{\gamma - 1})^{\nu + 1} ( C L )^{\gamma} x^{-\nu}.
$$
\end{lemma}

For completeness we show how Lemma\ \ref{maximal1} follows from the results of \citet{LS77}.  First, we state the relevant result of \citet{LS77}.

\begin{lemma}[\citet{LS77} Theorem 2]\label{maximal000}
Let $\lbrace X_l, 1 \leq l \leq L \rbrace$ be a completely general sequence of r.v.s.  Suppose there exist $\nu > 0$, $\gamma > 1$, and $C > 0$ such that for all $x > 0$ and non-negative integers $1 \leq i \leq j \leq L$, it holds that
$$
\pr\big(|\sum_{k=i}^j X_k| \geq x\big) \leq \big( C (j - i + 1) \big)^{\gamma} x^{-\nu}.
$$

Then it must also hold that
$$
\pr\big(\max_{i \in \lbrace 1,\ldots,L \rbrace} |\sum_{k=1}^i X_k| \geq x\big) \leq 
2^{\gamma} \bigg(1 + \big(2^{-\frac{1}{\nu + 1}} - 2^{-\frac{\gamma}{\nu + 1}}\big)^{-(\nu + 1)} \bigg) ( C L )^{\gamma} x^{-\nu}.
$$
\end{lemma}

Here let us point out that although \citet{LS77} presents several related bounds, they all seem to ultimately yield bounds scaling similar to those of Lemma\ \ref{maximal1} in our setting.  We leave it as an interesting open question whether bounds which exhibit significantly tighter scaling (e.g. as $r \downarrow 2$ and $r \uparrow \infty$) are possible for such maximal inequalities in our setting.
\\\\With Lemma\ \ref{maximal000} in hand, we now complete the proof of Lemma\ \ref{maximal1}.
\proof{Proof of Lemma\ \ref{maximal1}:}
As the conditions of Lemma\ \ref{maximal000} and Lemma\ \ref{maximal1} are identical, it suffices to prove that $\nu \geq \gamma$ (along with the other assumptions of Lemma\ \ref{maximal000}) implies 
\begin{equation}\label{maxeq11}
2^{\gamma} \bigg(1 + \big(2^{-\frac{1}{\nu + 1}} - 2^{-\frac{\gamma}{\nu + 1}}\big)^{-(\nu + 1)} \bigg) \leq 2.4 \times 2^{\gamma} \times (2 \frac{\nu + 1}{\gamma - 1})^{\nu + 1}.
\end{equation}
Note that
\begin{equation}\label{maxeq22}
2^{-\frac{1}{\nu + 1}} - 2^{-\frac{\gamma}{\nu + 1}} = 2^{-\frac{1}{\nu + 1}} \big(1 - 2^{- \frac{\gamma - 1}{\nu + 1} }\big).
\end{equation}
As our assumptions imply $0 < \frac{\gamma - 1}{\nu + 1} < 1$, and it is easily verified that $1 - 2^{-z} \geq \frac{z}{2}$ for all $z \in [0,1]$, we conclude that
\begin{equation}\label{maxeq33}
\big(1 - 2^{- \frac{\gamma - 1}{\nu + 1} }\big)^{-(\nu + 1)} \leq (2 \frac{\nu + 1}{\gamma - 1})^{\nu + 1}.
\end{equation}
Combining (\ref{maxeq22}) and (\ref{maxeq33}) with the fact that (by our assumptions) $(\frac{\nu + 1}{\gamma - 1})^{\nu + 1} \geq 1,$
we conclude that 
\begin{eqnarray*}
\\ &\ &\ \ 2^{\gamma} \bigg(1 + \big(2^{-\frac{1}{\nu + 1}} - 2^{-\frac{\gamma}{\nu + 1}}\big)^{-(\nu + 1)} \bigg)
\\&=& 2^{\gamma} \Bigg(1 + \bigg(2^{-\frac{1}{\nu + 1}} \big(1 - 2^{- \frac{\gamma - 1}{\nu + 1} }\big)\bigg)^{-(\nu+1)}\Bigg)
\\&\leq& 2^{\gamma} \big(1 + 2 (2 \frac{\nu + 1}{\gamma - 1})^{\nu + 1}\big).
\end{eqnarray*}
As our assumptions imply that $(2 \frac{\nu + 1}{\gamma - 1})^{\nu + 1} \geq 4$, the desired result then follows from some straightforward algebra.  $\halmos$
\endproof

\subsubsection{Proof of Lemma\ \ref{maximal2}.}\label{maximal2sec}

\proof{Proof of Lemma\ \ref{maximal2}:}
We proceed by verifying that for each fixed $k \geq 1$, the conditions of Lemma\ \ref{maximal1} hold for $\bigg\lbrace n - \sum_{i=1}^{n} \big( N_{e,i}(j) - N_{e,i}(j-1) \big) , j = 1,\ldots,k \bigg\rbrace$.
Let us fix some $k \geq 1$, and non-negative integers $l \leq m \leq k$.  Then for any $x > 0$, it follows from the fact that the given sequence of r.v.s is centered and stationary, the independence of $\lbrace N_{e,i}(t), i \geq 1 \rbrace$, our assumptions, and Markov's inequality (after raising both sides to the $r_1$ power), that for all $1 \leq l \leq m \leq k$,
\begin{eqnarray*}
&\ &\ \pr\bigg( \bigg| \sum_{j=l}^m \bigg( n - \sum_{i=1}^{n} \big( N_{e,i}(j) - N_{e,i}(j-1) \big) \bigg) \bigg| \geq x \bigg) 
\\&\ \ \ &\ \ \ \leq\ \ \  \E\bigg[ \bigg|\sum_{j=l}^m \bigg( n - \sum_{i=1}^{n} \big( N_{e,i}(j) - N_{e,i}(j-1) \big) \bigg) \bigg|^{r_1}\bigg] x^{-r_1}
\\&\ \ \ &\ \ \ =\ \ \ \E\bigg[ \bigg| n (m - l + 1) - \sum_{i=1}^{n} \big(N_{e,i}(m) - N_{e,i}(l-1)\big) \bigg|^{r_1} \bigg] x^{-r_1}
\\&\ \ \ &\ \ \ =\ \ \ \E\bigg[ \bigg| \sum_{i=1}^{n} N_{e,i}(m - l + 1) - n (m - l + 1) \bigg|^{r_1} \bigg] x^{-r_1}
\\&\ \ \ &\ \ \ \leq\ \ \ C_1 n^{\frac{r_1}{2}} (m - l + 1)^{s} x^{-r_1}.
\end{eqnarray*}
Thus we find that the conditions of Lemma\ \ref{maximal1} are met with $L = k$, $\lbrace X_l, 1 \leq l \leq L \rbrace = \bigg\lbrace  n - \sum_{i=1}^{n} \big( N_{e,i}(l) - N_{e,i}(l-1) \big), l = 1,\ldots,k \bigg\rbrace$, $C = (C_1 n^{\frac{r_1}{2}})^{\frac{1}{s}}$, $\nu = r_1$, $\gamma = s$, 
and the desired result follows.  $\halmos$
\endproof

\subsubsection{Proof of Lemma\ \ref{maximal3}.}\label{maximal3sec}

\proof{Proof of Lemma\ \ref{maximal3}:}
We begin by noting that it suffices to bound the supremum of interest over a suitable mesh, which follows immediately from the fact that w.p.1 $n t - \sum_{i=1}^{n} N_{e,i}(t)$ can increase by at most $2$ over any interval of length at most $\frac{2}{n}$.  We thus take our mesh to be $\lbrace \frac{2 k}{n}, k = 0,\ldots, \lfloor \frac{n t_0}{2} \rfloor \rbrace$, and conclude that  (\ref{tobound3}) is at most
\begin{equation}\label{maximal3eq2}
\pr\bigg(2 + \max_{k \in \lbrace 0 ,\ldots, \lfloor \frac{n t_0}{2} \rfloor\rbrace} \big( 2 k - \sum_{i=1}^{n} N_{e,i}(\frac{2 k}{n}) \big) \geq x\bigg).
\end{equation}
We now verify that the conditions of Lemma\ \ref{maximal1} hold for $\bigg\lbrace 2 - \sum_{i=1}^{n} \big( N_{e,i}(\frac{2 k}{n}) - N_{e,i}(\frac{2(k-1)}{n}) \big) , k = 0,\ldots, \lfloor \frac{n t_0}{2} \rfloor \bigg\rbrace$.  Let us fix some non-negative integers $m \leq j \leq \lfloor \frac{n t_0}{2} \rfloor$.  Then for any $x > 0$, it follows from stationary increments, centeredness, and Markov's inequality (after raising both sides to the $r_2$ power) that
\begin{equation}\label{maximal3eqqq4}
\pr\bigg( \bigg| \sum_{l=m}^j  \bigg( 2 - \sum_{i=1}^{n} \big( N_{e,i}(\frac{2 l}{n}) - N_{e,i}(\frac{2(l-1)}{n}) \big) \bigg) \bigg| \geq x \bigg)
\end{equation}
is at most
\begin{eqnarray*}
\ &\ &\ \ \ \E\bigg[\bigg| \sum_{l=m}^j \bigg( 2 - \sum_{i=1}^{n} \big( N_{e,i}(\frac{2l}{n}) - N_{e,i}(\frac{2(l-1)}{n}) \big) \bigg) \bigg|^{r_2}\bigg] x^{-r_2}
\\&\ &\ \ \ \ \ \ \ =\ \ \ \E\bigg[\big| 2 (j - m + 1) - \sum_{i=1}^{n} N_{e,i}(\frac{2(j - m + 1)}{n})\big|^{r_2}\bigg] x^{-r_2},
\end{eqnarray*}
which by our assumptions (and noting that in this case the $n t$ appearing in our assumptions equals $2 (j - m + 1)$) is at most $C_2 \big(2 (j - m + 1) \big)^{\frac{r_2}{2}} x^{-r_2}$.  We thus find that the conditions of Lemma\ \ref{maximal1} are met with $L = \lfloor \frac{n t_0}{2} \rfloor$, $\lbrace X_l, 1 \leq l \leq L \rbrace = \bigg\lbrace 2 - \sum_{i=1}^{n} \big( N_{e,i}(\frac{2 l}{n}) - N_{e,i}(\frac{2(l-1)}{n}) \big), l = 1,\ldots, \lfloor \frac{n t_0}{2} \rfloor \bigg\rbrace$, 
$C = 2 (C_2)^{\frac{2}{r_2}}$, $\nu = r_2$, $\gamma = \frac{r_2}{2}$.  Thus for all $x > 0$,
$$
\pr\bigg(\max_{k \in \lbrace 0 ,\ldots, \lfloor \frac{n t_0}{2} \rfloor\rbrace} \big( 2 k - \sum_{i=1}^{n} N_{e,i}(\frac{2 k}{n}) \big) \geq x \bigg) \leq
2.25 \times 2^{\frac{r_2}{2}} \times (2 \frac{r_2 + 1}{\frac{r_2}{2} - 1})^{r_2 + 1} ( C L )^{\frac{r_2}{2}} x^{- r_2}.$$
It then follows from (\ref{maximal3eq2}), and the fact that $x \geq 4$ implies $(x - 2)^{-r_2} \leq 2^{r_2} x^{-r_2}$, that
(\ref{tobound3}) is at most 
\begin{eqnarray*}
\ &\ &\ 2^{r_2} \times 2.25 \times 2^{\frac{r_2}{2}} \times (2 \frac{r_2 + 1}{\frac{r_2}{2} - 1})^{r_2 + 1} \times C_2 \times 2^{\frac{r_2}{2}} \times (\frac{n t_0}{2})^{\frac{r_2}{2}} x^{-r_2}
\\&\ &\ \ \ \ \ \ \leq\ \ \ .8 \times (5.7 \frac{r_2 + 1}{\frac{r_2}{2} - 1})^{r_2 + 1} \times C_2 \times (n t_0)^{\frac{r_2}{2}} x^{-r_2},
\end{eqnarray*}
completing the proof.  $\halmos$
\endproof

\subsection{Proof of Lemma\ \ref{poolerbound1}.}\label{poolerbound1sec}
Our proof of Lemma\ \ref{poolerbound1} proceeds in several steps.  First, we bound the $r$th central moment of $N_o(t)$, showing that this moment scales (with $t$) like $t^{\frac{r}{2}}$ and providing a completely explicit bound along these lines.  Second, we modify these bounds to instead yield bounds for the $r$th central moment of $N_e(t).$  Third, we apply general results from the literature for converting bounds on the moments of zero-mean r.v.s to bounds on the moments of sums of those r.v.s to convert the above bounds for individual cenetered renewal processes into bounds for centered pooled renewal processes.
\\\indent We begin by bounding the $r$th central moment of $N_o(t)$.  Our approach can essentially be viewed as ``making completely explicit", e.g. all constants explicitly worked out, the approach to bounding the central moments of a renewal process sketched in \citet{Gut09}.  As noted in \citet{Gut09} (and used in \citet{GG13}), a non-explicit bound proving that the $r$th central moment indeed scales asymptotically (with $t$) like $t^{\frac{r}{2}}$ was first proven in \citet{Chao79}.  To our knowledge such a completely explicit bound is new, and may prove useful in other settings.  In particular, we begin by proving the following.
\begin{lemma}\label{explicitrenewaltheorem1}
Suppose that $\E[S] = 1$, and that $\E[S^r] < \infty$ for some $r \geq 2$.  Then for all $t \geq 1$, $\E\bigg[\bigg|N_o(t) - t\bigg|^r\bigg]$ is at most
$$
\bigg( 1.5 \times \big( \E[S^2] + 1 \big)^r \times 120^r \times r^{1.5 r} + .4 \times (\E[S^2] + 1) \times 60^r \times r^r \times ( \E[S^r] + 1) \bigg) \times t^{\frac{r}{2}}.$$
\end{lemma}
\subsubsection{Preliminary results for the proof of Lemma\ \ref{explicitrenewaltheorem1}.}
Before proving Lemma\ \ref{explicitrenewaltheorem1}, let us prove some preliminary technical results.
First, we recall the celebrated Burkholder-Rosenthal Inequality for bounding the moments of a martingale.  We state a particular variant given in \citet{hitczenko1990best}.
\begin{lemma}[Burkholder-Rosenthal Inequality, \citet{hitczenko1990best}]\label{rosenthal1}  Let $\lbrace X_i, i \geq 1 \rbrace$ be a martingale difference sequence 
w.r.t. the filtration $\lbrace {\mathcal F}_i, i \geq 0 \rbrace$.  Namely, we have that $\lbrace X_i, i \geq 1 \rbrace$ is adapted to $\lbrace {\mathcal F}_i, i \geq 0 \rbrace$; $\E[|X_i|] < \infty$ for all $i \geq 1$; and $\E[X_i | {\mathcal F}_{i-1}] = 0$ for all $i \geq 1$.  Suppose also that $\lbrace \sum_{i=1}^n X_i, n \geq 1 \rbrace$ converges a.s. to a limiting r.v. which we denote $\sum_{i=1}^{\infty} X_i$.  Then for all $r \geq 2$,
$$
\bigg(\E\big[\big|\sum_{i=1}^{\infty} X_i\big|^r\big]\bigg)^{\frac{1}{r}} \leq 10 r \Bigg( \E\bigg[\bigg( \sum_{i=1}^{\infty} \E[X^2_i | {\mathcal F}_{i-1}] \bigg)^{\frac{r}{2}}\bigg] \Bigg)^{\frac{1}{r}}
+ 10 r \bigg( \E\big[\sup_{i \geq 1} |X^r_i| \big] \bigg)^{\frac{1}{r}}.$$
\end{lemma}
Since for any sequence of r.v.s $\lbrace Z_i, i = 1,\ldots,n \rbrace$ and $r \geq 1$ it follows from convexity that w.p.1
\begin{equation}\label{convexxx1}
\big| \sum_{i=1}^n Z_i \big|^r \leq n^{r-1} \sum_{i=1}^n |Z_i|^r,
\end{equation}
we deduce the following corollary.
\begin{corollary}\label{rosenthalcor1}
Under the same definitions and assumptions as Lemma\ \ref{rosenthal1}, for all $r \geq 2$,
$$
\E\big[\big|\sum_{i=1}^{\infty} X_i\big|^r\big] \leq (20 r)^r \E\bigg[\bigg( \sum_{i=1}^{\infty} \E[X^2_i | {\mathcal F}_{i-1}] \bigg)^{\frac{r}{2}}\bigg]+ (20 r)^r \E\big[\sup_{i \geq 1} |X^r_i| \big].$$
\end{corollary}
We next recall a certain inequality for the non-central moments of $N_o(t)$, proven in \citet{Gut09} Equation 5.11.
\begin{lemma}[\citet{Gut09} Equation 5.11]\label{yeq10}
For all $r \geq 1$ and $t \geq 1$,
$$\E\Big[ \big(  N_o(t) + 1 \big)^r \Big] \leq (2  t)^r \E\Big[ \big(  N_o(1) + 1 \big)^r \Big].
$$
\end{lemma}
Next, we prove several explicit bounds for the higher moments of a renewal process.  Although there is a large literature on bounds for renewal processes, including explicit bounds for the mean and variance of the number of renewals and non-explicit/asymptotic bounds for the higher moments of these processes (see e.g. \citet{Daley.78, Daley.80, Smith59,hunter69,Leadbettter63,Grubel.81,Taga63}), we believe our bounds to be novel, and potentially of independent interest.
\begin{lemma}\label{renewalbound00a}
Suppose that $\E[S] = 1$.  Then for all $r \geq 1$ and $\theta > 0$,
$$\E\big[\big(N_o(1) + 1\big)^r\big] \leq 1 + 1.5 \times r \times (r-1)^{r-1} \times \exp( 2 \theta ) \times \big(1 - \E[\exp(- \theta S)]\big)^{-r}.$$
\end{lemma}
\proof{Proof:}
Note that for all $j \geq 1$ and $\theta > 0$,
\begin{eqnarray}
\pr\big(N_o(1) + 1 \geq j\big) &=& \pr(\sum_{i=1}^{j-1} S_i \leq 1) \nonumber
\\&=& \pr\bigg( \exp\big( - \theta \sum_{i=1}^{j-1} S_i \big) \geq \exp(-\theta) \bigg) \nonumber
\\&\leq& \exp(\theta) \times \E^{j-1}[\exp(-\theta S)]\ \ \ \textrm{by Markov's inequality.} \label{showmarkov1111}
\end{eqnarray}
Applying the tail integral form for higher moments (\citet{Nad.22}), it follows that
\begin{eqnarray*}
\E[(N_o(1) + 1)^r] &=& r \int_0^{\infty} x^{r-1} P\big( N_o(1) + 1 \geq x \big) dx
\\&\leq& 1 + r \exp(\theta) \int_1^{\infty} x^{r-1} \E^{x-1}[\exp(-\theta S)] dx
\\&\leq& 1 + r \exp(\theta) \big(\E[\exp(- \theta S)]\big)^{-1} \int_0^{\infty} x^{r-1} \big( \E[\exp(- \theta S)] \big)^x dx
\\&=& 1 + r \exp(\theta) \big(\E[\exp(- \theta S)]\big)^{-1} \Gamma(r) \times \log^{-r}\big( \frac{1}{\E[\exp(- \theta S)]} \big) 
\\&\leq& 1+ r \exp(2 \theta) \Gamma(r) \big(1 - \E[\exp(- \theta S)]\big)^{-r},
\end{eqnarray*}
with the final inequality following from the fact that $\log(\frac{1}{x}) \geq 1 - x$ for all $x \in (0,1)$, and the fact that Jensen's inequality implies $\big(\E[\exp(- \theta S)]\big)^{-1} \leq \exp(\theta)$ since $\E[S] = 1$.  Combining with the fact that $\Gamma(1 + x) \leq 1.5 x^x$ for all $x \geq 0$ (which follows from the bounds of \citet{Batir17} Theorem 2.3) completes the proof.  $\halmos$
\endproof
To bound those terms involving $\theta$ in Lemma\ \ref{renewalbound00a}, we now prove the following result, which allows us to bound those terms (for an appropriate choice of $\theta$) purely in terms of the moments of $S$.  The proof uses a similar strategy to that used in our proof of Lemma\ \ref{kingmanmart2}.
\begin{lemma}\label{mgfboundlem}
Suppose that $\E[S] = 1$, and that $\E[S^2] < \infty.$  Then for $\theta = (2 \E[S^2])^{-1}$, $\big(1 - \E[\exp(- \theta S)]\big)^{-1} \leq 4 \E[S^2].$
\end{lemma}
\proof{Proof : }
Note that for all $\theta > 0$, w.p.1, $\exp(\theta S) \geq 1 + \theta S$ (by the exponential inequality), and hence 
\begin{eqnarray*}
\exp(- \theta S) &\leq& \frac{1}{1 + \theta S}
\\&=& 1 - \theta S + \frac{\theta^2 S^2}{1 + \theta S}
\\&\leq& 1 - \theta S + \theta^2 S^2.
\end{eqnarray*}
It follows that for all $\theta > 0$,
\begin{equation}\label{byeexpo1}
1 - \E[\exp(- \theta S)] \geq \theta \E[S] - \theta^2 \E[S^2].
\end{equation}
Taking $\theta = \frac{\E[S]}{2 \E[S^2]}$ and recalling that $\E[S] = 1$, we find that
\begin{equation}\label{byeexpo2}
\frac{1}{1 - \E[\exp(- \theta S)]} \leq 4 \E[S^2],
\end{equation}
completing the proof.  $\Halmos$
\endproof
Let us note that since in general $P(S = 0)$ may be positive (in which case $\E[\exp(-\theta S)]$ will not shrink to 0 as $\theta \uparrow \infty$), it seems dificult to derive generic bounds using large values of $\theta$ without introducing additional parameters related to the distribution of $S$, and thus we have instead taken advantage of the fact that Taylor series may be applied when $\theta$ is small.
\\\\Combining Lemmas\ \ref{renewalbound00a} and \ref{mgfboundlem} with some straightforward algebra, we come to the following corollary.

\begin{corollary}\label{renewalboundcccora}
Suppose that $E[S] = 1$.  Then for all $r \geq 1$,
$$\E\Big[ \big(  N_o(1) + 1 \big)^{r} \Big] \leq 4.4 \times \big( 4 \E[S^2] r)^r.$$
\end{corollary}

Further combining with Lemma\ \ref{yeq10}, we come to the following corollary.

\begin{corollary}\label{renewalboundcccor}
Suppose that $E[S] = 1$.  Then for all $r \geq 1$ and $t \geq 1$,
$$\E\Big[ \big(  N_o(t) + 1 \big)^{r} \Big] \leq 4.4 \times \big( 8 \E[S^2] r)^r \times t^r.$$
\end{corollary}

As it will arise in some of our calculations, let us also state a tighter known explicit bound for the first moment of $N_o(t) + 1$, which follows directly from \citet{Daley.80}.

\begin{lemma}[\citet{Daley.80}]\label{renewalboundcccor22}
Suppose that $E[S] = 1$.  Then for all $t \geq 1$,
$$\E\Big[ N_o(t) + 1 \Big] \leq (1 + \E[S^2]) t.$$
\end{lemma}
\subsubsection{Proof of Lemma\ \ref{explicitrenewaltheorem1}.} We now complete the proof of Lemma\ \ref{explicitrenewaltheorem1}.
\proof{Proof of Lemma\ \ref{explicitrenewaltheorem1}:}
By definition (as is well-known), $N_o(t) + 1 = \min\lbrace n \geq 1 : \sum_{i=1}^n S_i > t \rbrace$ is a stopping time w.r.t. the natural filtration generated by $\lbrace S_i, i \geq 1 \rbrace$.  By the triangle inequality, w.p.1
\begin{eqnarray}
\big| N_o(t)- t\big| &=& \big| \big(N_o(t) + 1\big) - t - 1\big| \nonumber
\\&\leq& \big| \big(N_o(t) + 1\big) - t\big| + 1 \nonumber
\\&\leq& \bigg|\sum_{i=1}^{N_o(t)+1} S_i - \big(N_o(t) + 1\big) \bigg| + \bigg| \sum_{i=1}^{N_o(t)+1} S_i - t \bigg| + 1. \label{yeq2}
\end{eqnarray}
It then follows from (\ref{convexxx1}) and (\ref{yeq2}) that 
\begin{eqnarray}
\ &\ &\ \E\bigg[\bigg|N_o(t) - t\bigg|^r\bigg]\nonumber
\\&\leq&\ 3^{r-1} \E\Bigg[\bigg|\sum_{i=1}^{N_o(t)+1} S_i - \big(N_o(t) + 1\big) \bigg|^r\Bigg]\label{yeq3}
\\&\ &\ \ \ +\ \ \ 3^{r-1}\E\Bigg[\bigg| \sum_{i=1}^{N_o(t)+1} S_i - t \bigg|^r\Bigg]\label{yeq4}
\\&\ &\ \ \ +\ \ \ 3^{r-1}.\label{yeq4b}
\end{eqnarray}
We next bound 
\begin{equation}\label{yeq5}
\E\Bigg[\bigg|\sum_{i=1}^{N_o(t)+1} S_i - \big(N_o(t) + 1\big) \bigg|^r \Bigg],
\end{equation}
and proceed by applying the Burkholder-Rosenthal Inequality.  In particular, we will use Corollary\ \ref{rosenthalcor1} to bound (\ref{yeq5}).  First, we rewrite (\ref{yeq5}) in terms of an appropriate martingale difference sequence.  Namely, note that (\ref{yeq5}) equals
\begin{eqnarray}
\ &\ &\ \E\Bigg[\bigg|\sum_{i=1}^{\infty} S_i I\big( N_o(t) + 1 \geq i \big) - \sum_{i=1}^{\infty} I\big( N_o(t) + 1 \geq i\big) \bigg|^r \Bigg]\nonumber
\\&\ &\ \ =\ \ \E\Bigg[\bigg|\sum_{i=1}^{\infty} (S_i - 1) I\big( N_o(t) + 1 \geq i \big) \bigg|^r \Bigg]. \label{yeq6}
\end{eqnarray}
We now prove that $\lbrace  (S_i - 1) I\big( N_o(t) + 1 \geq i \big), i \geq 1 \rbrace$ is a martingale difference sequence w.r.t. the filtration $\lbrace \sigma(S_1,\ldots,S_i), i \geq 1 \rbrace$.  Finite expectations and measurability are trivial.  Furthermore, since $I\big( N_o(t) + 1 \geq i \big)$ is $\sigma(S_1,\ldots,S_{i-1})$-measurable (due to the greater than or equal to sign), it follows from independence and the basic properties of conditional expectation that w.p.1
\begin{eqnarray*}
\ &\ &\ \ \E\bigg[ (S_i - 1) I\big( N_o(t) + 1 \geq i \big) | \sigma(S_1,\ldots,S_{i-1}) \bigg]
\\&\ &\ \ \ =\ \ \ I\big( N_o(t) + 1 \geq i \big) \E\bigg[ (S_i - 1) | \sigma(S_1,\ldots,S_{i-1}) \bigg]
\\&\ &\ \ \ =\ \ \ I\big( N_o(t) + 1 \geq i \big) \E\big[ S_i - 1 \big]\ \ \ =\ \ \ 0.
\end{eqnarray*}
Thus we find that the conditions of Corollary\ \ref{rosenthalcor1} are satisfied with $X_i =  (S_i - 1) I\big( N_o(t) + 1 \geq i \big), {\mathcal F}_i = \sigma(S_1,\ldots,S_i)$.  Before stating the given implication, we first show that several resulting terms can be simplified.  First, note that
\begin{eqnarray}
\ &\ &\ \E\Bigg[ \Bigg( \sum_{i=1}^{\infty} \E\bigg[ \bigg( (S_i - 1) I\big( N_o(t) + 1 \geq i \big) \bigg)^2 \bigg| \sigma(S_1,\ldots,S_{i-1}) \bigg] \Bigg)^{\frac{r}{2}} \Bigg]\nonumber
\\&\ &\ \ \ =\ \ \ \E\Bigg[ \Bigg( \sum_{i=1}^{\infty} \E\bigg[ (S_i - 1)^2 I\big( N_o(t) + 1 \geq i \big) \bigg| \sigma(S_1,\ldots,S_{i-1}) \bigg] \Bigg)^{\frac{r}{2}} \Bigg]\nonumber
\\&\ &\ \ \ =\ \ \ \E\Bigg[ \Bigg( \sum_{i=1}^{\infty} I\big( N_o(t) + 1 \geq i \big) \E\bigg[ (S_i - 1)^2 \bigg| \sigma(S_1,\ldots,S_{i-1}) \bigg] \Bigg)^{\frac{r}{2}} \Bigg]\nonumber
\\&\ &\ \ \ =\ \ \ \E\Bigg[ \Bigg( \sum_{i=1}^{\infty} I\big( N_o(t) + 1 \geq i \big) \E\big[ (S - 1)^2 \big] \Bigg)^{\frac{r}{2}} \Bigg]\nonumber
\\&\ &\ \ \ =\ \ \ \bigg(\E\big[ (S - 1)^2 \big]\bigg)^{\frac{r}{2}} \E\Bigg[ \Bigg( \sum_{i=1}^{\infty} I\big( N_o(t) + 1 \geq i \big) \Bigg)^{\frac{r}{2}} \Bigg]\nonumber
\\&\ &\ \ \ \leq\ \ \ \big(\E[S^2] + 1\big)^{\frac{r}{2}}  \E\Big[ \big(  N_o(t) + 1 \big)^{\frac{r}{2}} \Big].\label{yeq7}
\end{eqnarray}
Second, note that
\begin{eqnarray}
\ &\ &\ \E\Bigg[\sup_{i \geq 1} \bigg| \bigg( (S_i - 1) I\big( N_o(t) + 1 \geq i \big) \bigg)^r \bigg| \Bigg] \nonumber
\\&\ &\ \ \ \leq\ \ \ \E\Bigg[\sum_{i = 1}^{\infty} I\big( N_o(t) + 1 \geq i \big) \big| S_i - 1 \big|^r \Bigg] \nonumber
\\&\ &\ \ \ =\ \ \ \sum_{i = 1}^{\infty}\E\Bigg[I\big( N_o(t) + 1 \geq i \big) \big| S_i - 1 \big|^r \Bigg] \nonumber
\\&\ &\ \ \ =\ \ \ \sum_{i = 1}^{\infty}\pr\big( N_o(t) + 1 \geq i \big) \E\big[\big|S_i - 1\big|^r \big]\ \ \textrm{by independence}\nonumber
\\&\ &\ \ \ =\ \ \ \E\big[N_o(t) + 1\big] E\big[\big| S - 1 \big|^r\big], \label{yeq8}
\end{eqnarray}
the final inequality following since $\lbrace S_i, i \geq 1 \rbrace$ is i.i.d.  Combining (\ref{yeq7}) and (\ref{yeq8}) with the fact that the conditions of Corollary\ \ref{rosenthalcor1} are satisfied with $X_i =  (S_i - 1) I\big( N_o(t) + 1 \geq i \big), {\mathcal F}_i = \sigma(S_1,\ldots,S_i)$, we conclude that (\ref{yeq5}) is at most
\begin{equation}\label{yeq9}
(20 r)^r \big(\E[S^2] + 1\big)^{\frac{r}{2}}  \E\Big[ \big(  N_o(t) + 1 \big)^{\frac{r}{2}} \Big] + (20 r)^r \E\big[N_o(t) + 1\big] E\big[\big| S - 1 \big|^r\big]
\end{equation}
Combining (\ref{yeq9}), Corollary\ \ref{renewalboundcccor}, and Lemma\ \ref{renewalboundcccor22}, we conclude (after some straightforward algebra) that (\ref{yeq5}) is at most
\begin{equation}
4.4 \times \big(\E[S^2] + 1\big)^{\frac{r}{2}}\times \big( 40 \sqrt{\E[S^2]} r^{1.5})^{r} \times t^{\frac{r}{2}}
+
(20 r)^r (1 + \E[S^2]) E\big[\big| S - 1 \big|^r\big] t
. \label{yeq11}
\end{equation}
\ \\\\We next bound (\ref{yeq4}), by bounding
\begin{equation}\label{yeq12}
\E\Bigg[\bigg| \sum_{i=1}^{N_o(t)+1} S_i - t \bigg|^r\Bigg].
\end{equation}
By definition, $\sum_{i=1}^{N_o(t)+1} S_i - t$ is the residual life of the renewal process $N_o$ at time $t$, i.e. the remaining time until the next renewal (at time $t$), and it follows that w.p.1
\begin{eqnarray*}
\big|\sum_{i=1}^{N_o(t)+1} S_i - t\big|^r &\leq& S^r_{N_o(t)+1}
\\&\leq& \sum_{i=1}^{N_o(t) + 1} S^r_i.
\end{eqnarray*}
Combining with Wald's identity, we conclude that (\ref{yeq12}) is at most $\E\big[N_o(t) + 1\big] \E[S^r]$, also providing a bound for (\ref{yeq4}).  Again using Corollary\ \ref{renewalboundcccor} to bound $\E\big[N_o(t) + 1\big],$ applying (\ref{yeq11}) to bound (\ref{yeq5}), and combining with some straightforward algebra completes the proof.
$\halmos$
\endproof
\subsubsection{Extend Lemma\ \ref{explicitrenewaltheorem1} to the corresponding equilibrium renewal process.} We now extend Lemma\ \ref{explicitrenewaltheorem1} to the corresponding equilibrium renewal process.  We note that given the results of Lemma\ \ref{explicitrenewaltheorem1}, such an extension follows nearly identically to the proof of Lemma\ 8 of \citet{GG13}, although we include a self-contained proof for completeness.

\begin{corollary}\label{explicitrenewalcor1}
Suppose that $\E[S] = 1$, and $\E[S^r] < \infty$ for some $r \geq 2$.  Then for all $t \geq 1$, $\E\bigg[\bigg|N_e(t) - t\bigg|^r\bigg]$ is at most
$$
\bigg( .76 \times \big( \E[S^2] + 1 \big)^r \times 240^r \times r^{1.5 r} + .21 \times (\E[S^2] + 1) \times 120^r \times r^r \times (\E[S^r] + 1) \bigg) \times t^{\frac{r}{2}}.
$$
\end{corollary}
\proof{Proof:} Let $S^e$ denote the first renewal interval in ${\mathcal N}_e$, and $f_{S^e}$ its density function, whose existence is guaranteed by the basic properties of the equilibrium distribution.  Observe that we may construct ${\mathcal N}_e$ and ${\mathcal N}_o$ on the same probability space so that ${\mathcal N}_o$ is independent of $S^e$, and for all $t \geq 0$, w.p.1 
$$N_{e}(t) - t =  \bigg( N_o\big( ( t - S^e)^+ \big) - \big( t - S^e \big)^+ \bigg) I( S^e < t ) + \bigg( I(S^e \leq t) - \big(t - (t - S^e)^+ \big) \bigg).$$
Fixing some $t \geq 1$, it follows from (\ref{convexxx1}) and the triangle inequality that 
\begin{eqnarray}
\ &\ &\ \E\big[| N_{e}(t) - t |^r\big]\nonumber
\\&\leq& 2^{r - 1} \E\bigg[\big| N_o\big( ( t - S^e)^+ \big) - \big( t - S^e \big)^+\big|^r I( S^e < t ) \bigg]\label{verylarge2}
\\&\indent& + 2^{r-1} \E\big[ | I\big( S^e \leq t \big) - \big(t - (t - S^e)^+\big) |^r \big].\label{verylarge3}
\end{eqnarray}
We now bound the term $\E\bigg[\big| N_o\big( ( t - S^e)^+ \big) - \big( t - S^e \big)^+\big|^r I( S^e < t ) \bigg]$ appearing in (\ref{verylarge2}), which equals
\begin{eqnarray}
&\ & \int_0^{t-1} \E\big[| N_o\big( t - s \big) - \big( t - s \big) |^r\big] f_{S^e}(s) ds
+ \int_{t-1}^{t} \E\big[| N_o\big( t - s \big) - \big( t - s \big) |^r\big] f_{S^e}(s) ds. \label{verylarge2b}
\end{eqnarray}
Lemma\ \ref{explicitrenewaltheorem1} and Markov's inequality (after raising both sides to the $r$th power), combined with our assumptions on $r$ and $t$,  implies that the first summand of (\ref{verylarge2b}) is at most  
\begin{eqnarray*}
\ &\ &\ \bigg( 1.5 \times \big( \E[S^2] + 1 \big)^r \times 120^r \times r^{1.5 r} + .4 \times (\E[S^2] + 1) \times 60^r \times r^r \times ( \E[S^r] + 1) \bigg) \times 
\int_0^{t-1}  (t-s)^{\frac{r}{2}} f_{S^e}(s) ds
\\&\ &\ \ \ \leq\ \ \ \bigg( 1.5 \times \big( \E[S^2] + 1 \big)^r \times 120^r \times r^{1.5 r} + .4 \times (\E[S^2] + 1) \times 60^r \times r^r \times ( \E[S^r] + 1) \bigg) \times t^{\frac{r}{2}} \int_0^{t-1}  f_{S^e}(s) ds 
\\&\ &\ \ \ \leq\ \ \ \bigg( 1.5 \times \big( \E[S^2] + 1 \big)^r \times 120^r \times r^{1.5 r} + .4 \times (\E[S^2] + 1) \times 60^r \times r^r \times ( \E[S^r] + 1) \bigg) \times t^{\frac{r}{2}}.
\end{eqnarray*}
Since $t - s \leq 1$ implies w.p.1 $|N_o(t-s) - (t-s)|^r \leq \big(N_o(1) + 1\big)^r$, it follows from Corollary\ \ref{renewalboundcccora} that the second summand of (\ref{verylarge2b}) is at most 
\begin{eqnarray*}
4.4 \times \big(4 \E[S^2] r\big)^r.
\end{eqnarray*}
Combining the above, we find that (\ref{verylarge2}) is at most
\begin{eqnarray}
\ &\ &\ 2^{r-1} \times \bigg( 1.5 \times \big( \E[S^2] + 1 \big)^r \times 120^r \times r^{1.5 r} + .4 \times (\E[S^2] + 1) \times 60^r \times r^r \times ( \E[S^r] + 1) \bigg) \times t^{\frac{r}{2}}\nonumber
\\&\ &\ \ \ +\ \ \ 2^{r-1} \times 4.4 \times \big(4 \E[S^2] r\big)^r\nonumber
\\&\leq& \bigg( .76 \times \big( \E[S^2] + 1 \big)^r \times 240^r \times r^{1.5 r} + .2 \times (\E[S^2] + 1) \times 120^r \times r^r \times ( \E[S^r] + 1) \bigg) \times t^{\frac{r}{2}} \label{verylarge2b2}.
\end{eqnarray}
We now bound (\ref{verylarge3}), which is at most
\begin{eqnarray}
2^{2r - 2} \bigg(1 + \E\big[ | \big(t - (t - S^e)^+\big)|^r \big] \bigg)
&\leq& 2^{2r - 2}\bigg(1 + \big( \int_0^t s^r f_{S^e}(s) ds + \int_t^{\infty} t^r f_{S^e}(s) ds \big) \bigg).\label{verylarge33a}
\end{eqnarray}
It follows from the basic properties of the equilibrium distribution and Markov's inequality that for all $s \geq 0$,
$$f_{S^e}(s) = \pr(S > s) \leq \E[S^r] s^{-r}.$$
Thus the term $\int_0^t s^r f_{S^e}(s) ds + \int_t^{\infty} t^r f_{S^e}(s) ds$ appearing in (\ref{verylarge33a}) is at most
\begin{eqnarray}
\int_0^t s^r \big( \E[S^r] s^{-r} \big) ds + t^r \int_t^{\infty} \big( \E[S^r] s^{-r} \big) ds
&=& 
\E[S^r]\bigg( \int_0^t ds + t^r \int_t^{\infty} s^{-r} ds \bigg) \nonumber
\\&=&
\E[S^r] \bigg( t + t^r (r-1)^{-1} t^{1-r} \bigg) \nonumber
\\&\leq&
2 \E[S^r] t.\label{verylarge3b3}
\end{eqnarray}
Using (\ref{verylarge2b2}) to bound (\ref{verylarge2}), and (\ref{verylarge3b3}) and (\ref{verylarge33a}) to bound (\ref{verylarge3}), and combining with some straightforward algebra completes the proof.  $\halmos$
\endproof
\subsubsection{Result from literature to convert bounds for moments of zero-mean r.v.s to bounds for moments of sums of zero-mean r.v.s.} Before completing the proof of Lemma\ \ref{poolerbound1}, we recall the celebrated Marcinkiewicz-Zygmund inequality, a close relative of the Rosenthal inequality.  The precise result which we will use follows immediately from \citet{Ren.01} Theorem 2, and we refer the interested reader to \citet{F97} for a further overview of related results.  We note that for several results which we will state, it is not required that the r.v.s be identically distributed, although we only state the results for that setting.

\begin{lemma}[\citet{Ren.01} Theorem 2]\label{csumbound}
Suppose that for some $r \geq 2$, $\lbrace X_i, i \geq 1 \rbrace$ is a collection of i.i.d. zero-mean r.v.s. s.t. $\E[|X_1|^r] < \infty$.  Then for all $k \geq 1$,
$$\E\big[ \big| \sum_{i=1}^k X_i \big|^r \big] \leq \big(4.3 r^{\frac{1}{2}}\big)^r \E[|X_1|^r] k^{\frac{r}{2}}.$$
\end{lemma}

We note that the bounds of Lemma\ \ref{csumbound}, in particular the $r^{\frac{r}{2}}$ scaling, are tight even in the i.i.d. case (\citet{Ren.01}).
\subsubsection{Proof of Lemma\ \ref{poolerbound1}.} With Corollary\ \ref{explicitrenewalcor1} and Lemma\ \ref{csumbound} in hand, we now complete the proof of Lemma\ \ref{poolerbound1}.

\proof{Proof of Lemma\ \ref{poolerbound1}:}
Applying Lemma\ \ref{csumbound} with $X_i = N_{e,i}(t) - t$, we find that 
$$\E\big[ \big|\sum_{i=1}^n N_{e,i}(t) - n t \big|^r\big] \leq \big(4.3 r^{\frac{1}{2}}\big)^r E\big[|N_e(t) - t|^r\big] n^{\frac{r}{2}}.$$
Combining with Corollary\ \ref{explicitrenewalcor1} and some straightforward algebra completes the proof.  $\halmos$
\endproof

\subsection{Proof of Lemma\ \ref{binomial2}.}\label{binomial2sec}
We will prove Lemma\ \ref{binomial2} by breaking up the term $\sum_{i=1}^n N_{e,i}(t)$ in such a manner that our bounds for moments of sums, such as Lemma\ \ref{csumbound}, can be applied to achieve the desired $\frac{1}{1-\rho}$ scaling.  This requires us to show that $\E\big[ \big|\sum_{i=1}^n N_{e,i}(t) - n t\big|^r \big]$ scales (jointly in $n,t$) as $(nt)^{\frac{r}{2}}$ when $t$ may be near 0, and may also scale non-trivially with $n$.  Note that our analysis for the $t \geq 1$ case used heavily results for the scaling of the higher central moments of renewal processes for $t \geq 1$, quantifying explicitly the asymptotic scaling for large $t$, and those results are no longer applicable when $t$ is very small.  Instead, we will proceed as follows.  We rewrite $\E\big[ \big|\sum_{i=1}^n N_{e,i}(t) - n t\big|^r \big]$ as a double-sum (plus remainder term), in such a way that two essential properties hold.  Let $n_1(t) \stackrel{\Delta}{=} \lfloor n t \rfloor$, and $n_2(t) \stackrel{\Delta}{=} \lfloor \frac{n}{n_1(t)} \rfloor$.  Then we rewrite $\sum_{i=1}^n N_{e,i}(t) - n t$ as $$\sum_{m = 1}^{n_1(t)} \sum_{l=1}^{ n_2(t) }  \big( N_{(m-1)n_2(t) + l } ( t ) - t \big)  +  \sum_{l = n_1(t) n_2(t) + 1}^n \big( N_l ( t ) - t \big).$$  The two essential properties are the following.  First, $n_1(t)$ scales roughly as $n t$, which means that if we apply Lemma\ \ref{csumbound} to $\sum_{m = 1}^{n_1(t)} \sum_{l=1}^{ n_2(t) }  \big( N_{(m-1)n_2(t) + l } ( t ) - t \big)$ by thinking of each $\sum_{l=1}^{ n_2(t) }  \big( N_{(m-1)n_2(t) + l } ( t ) - t \big)$ term as its own r.v. $Y_m$ (and thus thinking of the overall sum as $\sum_{m=1}^{n_1(t)} Y_m$ with $\lbrace Y_m, m = 1,\ldots,n_1(t) \rbrace$ i.i.d.), Lemma\ \ref{csumbound} would yield a $(n t)^{\frac{r}{2}}$ scaling as long as we could sufficiently control the moments of each $Y_m$.  Second, $n_2(t)$ scales roughly as $\frac{1}{t}$, which will allow us to think of each such $Y_m$ term as the sum of $\frac{1}{t}$ independent terms each of which is 0 with probability roughly $1 - t$, and some modest value with probability $t$, since an equilibrium renewal process over a small interval $t$ has no events with probability roughly $1-t$.  This intuition will allow us to tightly bound the moments of each $Y_m$ roughly by thinking of $Y_m$ as a modified Binomial r.v. (with mean which can be bounded independent of $n,t$), in such a way that the higher moments of each $Y_m$ will not scale with $nt$.  By showing the remainder term $\sum_{l = n_1(t) n_2(t) + 1}^k \big( N_l ( t ) - t \big)$ consists of so few terms that its moments can also be sufficiently bounded, combining all of the above will yield the desired $(nt)^{\frac{r}{2}}$ scaling.
\\\\To implement the above approach and prove the desired $(n t)^{\frac{r}{2}}$ scaling and Lemma\ \ref{binomial2}, we first prove a bound for $\E\big[ \big|\sum_{i=1}^k N_{e,i}(t) - k t\big|^r \big]$ which will be applicable for the inner terms of the aforementioned double sum, for which $k$ scales roughly as $\frac{1}{t}$.  As mentioned above, the proof proceeds by interpreting $\sum_{i=1}^k N_{e,i}(t)$ as a modified binomial random variable.  To make the overall proof of Lemma\ \ref{binomial2} more readable, we defer the proof to the end of this section, here only stating the relevant result.

\begin{lemma}\label{binomial1}
Suppose that $\E[S] = 1$.  Then for all $k \geq 1, r \geq 2, t \in [0, \min(\frac{2}{k},1)],$ 
\begin{equation}\label{binomial1b}
\E\big[ \big|\sum_{i=1}^k N_{e,i}(t) - k t\big|^r \big] \leq 9.1 \times \big( \E[S^2] + 1 \big)^r \times 4^r \times r^{2 r}.
\end{equation}
\end{lemma}

\proof{Proof of Lemma\ \ref{binomial2}:}
Let $n_1(t) \stackrel{\Delta}{=} \lfloor n t \rfloor$.  Noting that $t \geq \frac{2}{n}$\ implies\ $n_1(t) > 0$, in this case we may define $n_2(t) \stackrel{\Delta}{=} \lfloor \frac{n}{n_1(t)} \rfloor$.  Then the left-hand-side of (\ref{focus1}) equals
\begin{eqnarray}
\ &\ &
\E\big[\big| \sum_{m = 1}^{n_1(t)} \sum_{l=1}^{ n_2(t) }  \big( N_{e,(m-1)n_2(t) + l } ( t ) - t \big)  +  \sum_{l = n_1(t) n_2(t) + 1}^k \big( N_{e,l} ( t ) - t \big)
 \big|^r\big] \nonumber
\\&\leq& 2^{r - 1} \E\big[ \big| \sum_{m = 1}^{n_1(t)} \sum_{l=1}^{ n_2(t)} \big( N_{e, (m-1)n_2(t) + l } ( t ) - t \big) \big|^r \big] \label{maxme1}
\\&\indent& + 2^{r-1} \E\big[ \big| \sum_{l = n_1(t) n_2(t) + 1}^k \big( N_{e,l} ( t ) - t \big) \big|^{r}\big].\label{maxme2}
\end{eqnarray}
We now bound (\ref{maxme1}).  First, let us apply Lemma\ \ref{csumbound} to conclude that (\ref{maxme1}) is at most
\begin{equation}\label{maxme11a}
2^{r - 1} \times \big(4.3 r^{\frac{1}{2}}\big)^r \times \E\big[ \big| \sum_{l=1}^{n_2(t)} \big( N_{e,l}(t) - t \big) \big|^r \big] \times \big( n_1(t) \big)^{\frac{r}{2}}.
\end{equation}
Next, we show that we may apply Lemma\ \ref{binomial1} to $\E\big[ \big| \sum_{l=1}^{n_2(t)} \big( N_{e,l}(t) - t \big) \big|^r \big]$, by arguing that $t \leq \frac{2}{n_2(t)}$.  In particular, 
\begin{equation}\label{abcd2}
t n_2(t)\ \ \ =\ \ \ t \lfloor \frac{n}{ \lfloor n t \rfloor } \rfloor\ \ \ \leq\ \ \ \frac{nt}{nt - 1}.
\end{equation}
But since $t \geq \frac{2}{n}$ implies $nt \geq 2$, and $g(z) \stackrel{\Delta}{=}\frac{z}{z-1}$ is a decreasing function of $z$ on $(1,\infty)$, it follows from (\ref{abcd2}) that $t n_2(t) \leq 2$.  Thus we may apply Lemma\ \ref{binomial1} (with $k = n_2(t)$, along with the fact that $n_1(t) \leq n t, t n_2(t) \leq 2$, and $t \leq 1$) to conclude that (\ref{maxme1}) is at most
\begin{eqnarray}
&\ & 2^{r - 1} \times \big(4.3 r^{\frac{1}{2}}\big)^r \times 
9.1 \times \big( \E[S^2] + 1 \big)^r \times 4^r \times r^{2 r} \times (n t)^{\frac{r}{2}} \nonumber
\\&\leq& 4.6 \times \big( 35 (1 + \E[S^2]) \big)^r \times r^{2.5 r} \times (nt)^{\frac{r}{2}}.
\label{abcd1}
\end{eqnarray}
We now bound (\ref{maxme2}).  Note that the sum $\sum_{l = n_1(t) n_2(t) + 1}^n \big( N_l ( t ) - t \big)$ appearing in (\ref{maxme2}) is taken over $n - n_1(t) n_2(t)$ terms.  Furthermore,
\begin{eqnarray*}
n - n_1(t) n_2(t) &=& n - n_1(t) \lfloor \frac{n}{n_1(t)} \rfloor 
\\&\leq& n - n_1(t) \big( \frac{n}{n_1(t)} - 1 \big) 
\\&=& n_1(t).
\end{eqnarray*}
As $n_1(t) \leq n t$, it thus follows from Lemma\ \ref{csumbound} that (\ref{maxme2}) is at most
\begin{equation}\label{abcd6aa0}
2^{r-1} \times \big(4.3 r^{\frac{1}{2}}\big)^r \times (n t)^{\frac{r}{2}} \times \E\big[| N_e(t) - t |^r \big].
\end{equation}
Noting that $|N_e(t) - t|$ is stochastically dominated by $N_o(1) + 1$ for all $t \leq 1$ (using the basic relationship between equilibrium and ordinary renewal processes), we may use Corollary\ \ref{renewalboundcccora} to bound 
$\E\big[| N_e(t) - t |^r\big]$ by $4.4 \times \big( 4 \E[S^2] r)^r.$  We conclude that (\ref{maxme2}) is at most 
$$2.2 \times (35 \E[S^2])^r \times r^{1.5 r} \times (nt)^{\frac{r}{2}}.$$
Combining the above with some straightforward algebra completes the proof.  $\halmos$
\endproof

\subsubsection{Proof of Lemma\ \ref{binomial1}.}
We now complete the proof of Lemma\ \ref{binomial1}.  First, let us state a bound for the uncentered moments of sums of i.i.d. non-negative random variables from the literature, which implies a simple and explicit bound on the moments of a binomially distributed r.v..  The result follows immediately from the results of \citet{Berend10}.

\begin{lemma}[\citet{Berend10}]\label{csumboundbb2}
Suppose that for some $r \geq 2$, $\lbrace X_i, i \geq 1 \rbrace$ is a collection of i.i.d. non-negative r.v.s. s.t. $\E[X^r_1] < \infty$.  Then for all $k \geq 1$,
$$\E\big[ \big( \sum_{i=1}^k X_i \big)^r \big] \leq r^r \max\big( (k \E[X_1])^r , k \E[X^r_1] \big).$$
\end{lemma}

We note that the bounds of \citet{Berend10} in fact show that the $r^r$ scaling of Lemma\ \ref{csumboundbb2} can be improved slightly (for large $r$) to $(\frac{r}{\log(r)})^r$, but that this $(\frac{r}{\log(r)})^r$ scaling is essentially tight, even for the moments of a binomial distribution (\citet{Ahle22}).  For simplicity, and as it will not substantially change the asymptotics of our final bounds since either way a term scaling as $r^r$ persists, here we use the simpler bound $r^r$.

\proof{Proof of Lemma\ \ref{binomial1}:}
Since $|a-b|^r \leq a^r + b^r$ for any $a,b \in {\mathcal R}^+$, the left-hand-side of (\ref{binomial1b}) is at most 
\begin{equation}
\E\big[\big(\sum_{i=1}^k N_{e,i}(t)\big)^r\big] + (k t)^r.\label{referenceeq1}
\end{equation}
We now bound the term $\E\big[\big(\sum_{i=1}^k N_{e,i}(t)\big)^r\big]$ appearing in (\ref{referenceeq1}).  Let $\lbrace B_i , i \geq 1\rbrace$ denote a sequence of i.i.d. Bernoulli r.v. s.t $\pr(B_i = 1) = p_t \stackrel{\Delta}{=} \pr( R(S) \leq t)$, and $\pr(B_i = 0) = 1 - p_t$.  Note that we may construct $\lbrace N_{e,i}(t) , i \geq 1 \rbrace$, $\lbrace N_{o,i}(t) , i \geq 1 \rbrace, \lbrace B_i , i \geq 1 \rbrace$ on the same probability space s.t. w.p.1 $N_{e,i}(t) \leq B_i \big(1 + N_{o,i}(t) \big)$ for all $i \geq 1$, with $\lbrace N_{o,i}(t) , i \geq 1 \rbrace, \lbrace B_i , i \geq 1 \rbrace$ mutually independent.
Let $M_t \stackrel{\Delta}{=} \sum_{i=1}^k B_i$, i.e. $M_t$ is the corresponding binomially distributed r.v.  Then it follows from Lemma\ \ref{csumboundbb2}, Corollary\ \ref{renewalboundcccor}, Lemma\ \ref{renewalboundcccor22}, the fact that $t \leq 1$, and Jensen's inequality that 
\begin{eqnarray*}
\E\big[ \big( \sum_{i=1}^k N_{e,i}(t) \big)^r \big] &\leq& \E\big[ \bigg( \sum_{i=1}^{ M_t } \big(1 + N_{o,i}(t) \big) \bigg)^r \big] \nonumber
\\&=& \E\bigg[ \E\big[ \bigg( \sum_{i=1}^{ M_t } \big(1 + N_{o,i}(t) \big) \bigg)^r | M_t \big] \bigg]\nonumber
\\&\leq& \E\bigg[ r^r \max\bigg( \big(M_t \E[1 + N_o(t)]\big)^r , M_t \E[\big(1 + N_o(t)\big)^r] \bigg) \bigg]
\\&\leq& r^r \bigg( \E\big[(M_t)^r\big] \times \big(\E[1 + N_o(1)]\big)^r + \E[M_t] \E\big[\big(1 + N_o(1)\big)^r] \bigg)
\\&\leq& r^r \E\big[ (M_t)^r \big] \times \big( 1 + \E[S^2] \big)^r
\\&\ &\ \ \ \ +\ \ r^r \E[M_t] \times 4.4 \times \big( 4 \E[S^2] r)^r.
\end{eqnarray*}
It follows from Lemma\ \ref{csumboundbb2} that
$$\E[ M^r_t ] \leq r^r \max\big( (k p_t)^r, k p_t \big),$$
and $\E[M_t] = k p_t$ since $M_t$ has a binomial distribution.  We may combine the above and find that $\E\big[ \big( \sum_{i=1}^k N_{e,i}(t) \big)^r \big]$ is at most
\begin{eqnarray*}
\ &\ & r^{2 r} \times \big( 1 + \E[S^2] \big)^r \max\big( (k p_t)^r, k p_t \big)
\\&\ &\ \ \ +\ \ 4.4 \times r^r \times \big( 4 \E[S^2] r)^r \times k p_t.
\end{eqnarray*}
Since it follows from the definition of the equilibrium distribution and $p_t$ that $p_t \leq t$ (as here we are assuming $\E[S] = 1$), and as $k t \leq 2$, the desired result then follows from straightforward algebra.  $\Halmos$ \endproof 

\subsection{Sketch of plausible approach generalizing our results to queueing networks.}\label{appnetworksec}
It is natural and interesting to ask whether our approach extends to more complex queueing systems, such as queueing networks.  Here we explore this question at a somewhat informal / conjectural level, sketching a possible approach, and leaving a formal investigation as an interesting direction for future research.  For simplicity, let us restrict our discussion to a tandem system of two $n$-server queues, one upstream and one downstream, although note that (as we will see) this setting already captures many of the complexities of such an extension.  Furthermore, suppose all external arrivals are to the upstream queue, that arrival process is Markovian with rate $\lambda$, and all service times are i.i.d. with mean one.  In this setting, the natural extension of our approach would be to consider a modified system in which an extra arrival is added to the upstream queue whenever a server would otherwise have gone idle in the upstream system, and an extra arrival is added to the downstream queue whenever a server would otherwise have gone idle in the downstream system.  It seems likely our approach could be extended to this setting to yield bounds in terms of certain suprema of processes which are the difference of pooled renewal processes, or for general networks the splitting and merging of appropriate renewal processes.  Intuitively, the ``input" at certain queues would now be the splitting and merging of pooled renewal processes representing the departures from other queues.  The relevant monotonocities necessary for such a modification to yield upper bounds (on e.g. the total number in system) should follow from known results for stochastic comparison of queueing networks (\citet{Shanthikumar89, Chen01}).  \\\indent However, a naive implementation of such a bounding methodology does not work in the network setting, as we now explain.  In particular, our bounding methodology would cause the departure process from the upstream system to become the pooling of $n$ renewal processes.  However, in general this will drive the downstream system into instability.  Indeed, the SLLN for renewal processes implies that the long-run rate at which work departs the upstream system (and heads to the downstream system) under the modifications required for our approach will be $n$, not $\lambda$ as in the original system, which will overload (critically load, to be more precise) the downstream system.  We note that this ``overloading phenomena" only arises in the network setting, since when there is a single multi-server queue the ``extra arrivals" only occur when a server would have anyways gone idle.  Alternatively, in the network setting, there is not ``sufficient coordination" between the ``extra arrivals" at different queues to prevent instability.  
\\\indent Perhaps surprisingly, this issue is not insurmountable.  As explored in \citet{Chang92,Chang94} for networks of single-server queues, a viable approach to overcome this problem is as follows.  First, one ``slows down" the service times at the upstream station (e.g. by simply multiplying all service times at the first station by some constant inflation factor greater than one), in such a way that stability at the upstream station is maintained.  Then, on this modified system (in which service times at the upstream station are now stochastically larger than at the downstream system), one implements our approach.  With the upstream station services ``slowed down", the departure process from the upstream station (under the modifications required for our approach) will no longer induce instability at the downstream station.  Interestingly, it is shown in \citet{Chang92,Chang94} that for a broad class of single-server queueing networks, it is always possible to implement such an approach (i.e. slowing down service times at each queue by an appropriate factor) such that under this construction stability is maintained for the overall network.  Although those works considered networks of single-server queues, the general methodology seems likely to directly extend to multi-server queues.  Furthermore, such a transformation will again lead to an upper bound, using the same standard results for comparison of queueing networks (\citet{Shanthikumar89, Chen01}).  However, those works only show that such a transformation can be implemented to preserve stability, without studying how this would effect the scaling of queue lengths.  
\\\indent We conjecture that such an approach can indeed be implemented to yield general and explicit bounds with an appropriate analogue of $\frac{1}{1-\rho}$ scaling for a broad range of queuing networks.  Although we are not aware of any simple and explicit analogues of Kingman's bound for queueing networks conjectured in the literature, we note that past work on heavy-traffic in queueing networks suggests that the number in queue at each station $i$ should scale as $\frac{1}{1-\rho_i}$ with $\rho_i$ the effective traffic intensity at that station (as dictated by the so-called traffic equations, see e.g. \citet{Reiman84, MMR.98, GZ06, Dai20}).  For example, for the simple 2-queue tandem queue described above, such a scaling can be acheived by ``slowing down" service times at the upstream station by multiplying service times at that station by $\sqrt{\frac{n}{\lambda}}.$  Under such a transformation, the upstream station becomes an $n$-server queueing system with arrival rate $\lambda' = \lambda$ and service rate $\mu' = \sqrt{\frac{\lambda}{n}}$ (with additional arrivals as appropriate when a server would go idle), and the downstream station becomes an $n$-server queue with arrival rate $n \mu'$ (coming from the departure process at the upstream station) and service rate 1 (again with additional arrivals when servers would go idle).  Thus both the upstream and downstream queues effectively become n-server queues with traffic intensity $\frac{\lambda}{n\sqrt{\frac{\lambda}{n}}} = \frac{n\sqrt{\frac{\lambda}{n}}}{n} = \sqrt{\frac{\lambda}{n}}$, with extra arrivals when a server would otherwise go idle.  But as the effective traffic intensity at both stations in the original system is easily seen to be $\frac{\lambda}{n}$, and as it is easily verified that $\frac{1}{1 - \sqrt{\frac{\lambda}{n}}} \leq 2 \times \frac{1}{1 - \frac{\lambda}{n}}$ for all $n \geq 0$ and $\lambda  \in (0,n),$ we find that such a transformation indeed preserves the desired $\frac{1}{1-\rho}$ scaling at each station in an appropriate sense.  We leave a formal investigation along these lines as an interesting direction for future research, and point out that for more general (e.g. multi-class) queueing networks questions of stability in networks can indeed be quite subtle (\citet{Dai20}).

\begin{thebibliography}{}

\bibitem[{Abate et al. (1991)}]{Abate94}
Abate, J., L. C. Gagan, W. Whitt. "Waiting-time tail probabilities in queues with long-tail service-time distributions." Queueing systems 16 (1994): 311-338.

\bibitem[{Aghajani and Ramanan (2020)}]{AR16}
Aghajani, R., K. Ramanan. ``The limit of stationary distributions of many-server queues in the Halfin–Whitt regime." Mathematics of Operations Research 45, no. 3 (2020): 1016-1055.

\bibitem[{Ahle (2022)}]{Ahle22}
Ahle, T. "Sharp and simple bounds for the raw moments of the binomial and Poisson distributions." Statistics and Probability Letters 182 (2022): 109306.

\bibitem[{Allen (2014)}]{Allen.14}
Allen, A.O. Probability, statistics, and queueing theory. Academic Press, 2014.

\bibitem[{Arjas et al. (1978)}]{Arjas78}
Arjas, E., T. Lehtonen. ``Approximating many server queues by means of single server queues." Mathematics of Operations Research 3.3 (1978): 205-223.

\bibitem[{Asmussen et al. (1996)}]{AS96}
Asmussen, S., O. Nerman, M. Olsson. ``Fitting phase-type distributions via the EM algorithm." Scandinavian Journal of Statistics (1996): 419-441.

\bibitem[{Asmussen (2008)}]{AS08}
Asmussen, S. Applied probability and queues. Vol. 51. Springer Science and Business Media, 2008.

\bibitem[{Atar et al. (2011)}]{Atar11}
Atar, R., N. Solomon. ``Asymptotically optimal interruptible service policies for scheduling jobs in a diffusion regime with nondegenerate slowdown." Queueing Systems 69.3 (2011): 217-235.

\bibitem[{Baccelli et al. (2013)}]{Bac13}
Baccelli, F., P. Bremaud. Elements of queueing theory: Palm Martingale calculus and stochastic recurrences. Vol. 26. Springer Science and Business Media, 2013.

\bibitem[{Bandi et al. (2015)}]{Bandi15}
Bandi, C., D. Bertsimas, N. Youssef. ``Robust queueing theory." Operations Research 63.3 (2015): 676-700.

\bibitem[{Batir (2017)}]{Batir17}
Batır, N. "Bounds for the Gamma Function." Results in Mathematics 72 (2017): 865-874.

\bibitem[{Bazhba et al. (2019)}]{Bazhba19}
Bazhba, M., J. Blanchet, C. Rhee, B. Zwart. "Queue length asymptotics for the multiple-server queue with heavy-tailed Weibull service times." Queueing Systems 93 (2019): 195-226.

\bibitem[{Beesack (1969)}]{Beesack69}
Beesack, P. "Improvements of Stirling's formula by elementary methods." Publikacije Elektrotehničkog fakulteta. Serija Matematika i fizika 274/301 (1969): 17-21.

\bibitem[{Berger et al. (1992)}]{Berger92}
Berger, A., W. Whitt. ``Comparisons of multi-server queues with finite waiting rooms." Stochastic Models 8, no. 4 (1992): 

\bibitem[{Berend et al. (2010)}]{Berend10}
Berend, D., T. Tassa. "Improved bounds on Bell numbers and on moments of sums of random variables." Probability and Mathematical Statistics 30, no. 2 (2010): 185-205.

\bibitem[{Bertsimas (1995)}]{Bertsimas95}
Bertsimas, D., D. Nakazato. ``The distributional Little's law and its applications." Operations Research 43.2 (1995): 298-310.

\bibitem[{Bhattacharya et al. (1991)}]{Bhat91}
Bhattacharya, P., A. Ephremides. ``Stochastic monotonicity properties of multiserver queues with impatient customers." Journal of Applied Probability 28, no. 3 (1991): 673-682.

\bibitem[{Borovkov (1965)}]{Borovkov65}
Borovkov, A.A. ``Some Limit Theorems in the Theory of Mass Service, II Multiple Channels Systems." Theory of Probability and Its Applications 10.3 (1965): 375-400.

\bibitem[{Borst et al. (2004)}]{Borst04}
Borst, S., A. Mandelbaum, M. Reiman. ``Dimensioning large call centers." Operations research 52, no. 1 (2004): 17-34.

\bibitem[{Braverman et al. (2017)}]{BDF15}
Braverman, A., J. G. Dai, and J. Feng. ``Stein’s method for steady-state diffusion approximations: an introduction through the Erlang-A and Erlang-C models." Stochastic Systems 6, no. 2 (2017): 301-366.

\bibitem[{Braverman and Dai (2015)}]{BD15}
Braverman, A., J. G. Dai. ``Stein's method for steady-state diffusion approximations of $ M/Ph/n+ M $ systems." Annals of applied probability 27(1) : 550 - 581.

\bibitem[{Braverman and Dai (2016)}]{BD16}
Braverman, A., J. G. Dai. ``High order steady-state diffusion approximation of the Erlang-C system." arXiv preprint arXiv:1602.02866 (2016).

\bibitem[{Braverman and Dai (2020)}]{BD20}
Braverman, A., J. G. Dai, X. Fang. ``High order steady-state diffusion approximations." arXiv preprint arXiv:2012.02824 (2020).

\bibitem[{Brumelle (1973)}]{Brumelle73}
Brumelle, S. ``Bounds on the Wait in a GI/M/k Queue." Management Science 19.7 (1973): 773-777.  

\bibitem[{Brumelle et al. (1975)}]{Brumelle75}
Brumelle, S. L., R. G. Vickson. ``A unified approach to stochastic dominance." In Stochastic optimization models in finance, pp. 101-113. Academic Press, 1975.

\bibitem[{Burman et al. (1983)}]{Burman83}
Burman, D., Donald. Smith. ``A light-traffic theorem for multi-server queues." Mathematics of Operations Research 8, no. 1 (1983): 15-25.

\bibitem[{Chang et al. (1992)}]{Chang92}
Chang, C. ``Stability, queue length and delay. I. Deterministic queueing networks." In [1992] Proceedings of the 31st IEEE Conference on Decision and Control, pp. 999-1004. IEEE, 1992.

\bibitem[{Chang et al. (1994)}]{Chang94}
Chang, C., J. Thomas, S. Kiang. ``On the stability of open networks: a unified approach by stochastic dominance." Queueing systems 15, no. 1 (1994): 239-260.**

\bibitem[{Chao et al. (1979)}]{Chao79}
Chao, Y., C. Hsiung, T. Lai, ``Extended renewal theory and moment convergence in Anscombe's theorem."  The Annals of Probability 7 (1979), no. 2, 304-318.

\bibitem[{Chawla et al. (2017)}]{Chawla17}
Chawla, S., N. Devanur, A. Holroyd, A. Karlin, J. Martin, B. Sivan. ``Stability of Service under Time-of-Use Pricing."  In Proceedings of the 49th Annual ACM SIGACT Symposium on Theory of Computing, pp. 184-197. 2017.

\bibitem[{Chen et al. (1986)}]{Chen86}
Chen, J., H. Rubin. ``Bounds for the difference between median and mean of gamma and Poisson distributions." Statistics and probability letters 4, no. 6 (1986): 281-283.

\bibitem[{Chen et al. (2001)}]{Chen01}
Chen, H., D. Yao. Fundamentals of queueing networks: Performance, asymptotics, and optimization. Vol. 4. New York: Springer, 2001.

\bibitem[{Dai et al. (2010)}]{Dai10}
Dai, J. G., S. He, T. Tezcan. ``Many-server diffusion limits for G/Ph/n+ GI queues." The Annals of Applied Probability 20, no. 5 (2010): 1854-1890.

\bibitem[{Dai et al. (2014)}]{DDG14}
Dai, J. G., A. B. Dieker, X. Gao. ``Validity of heavy-traffic steady-state approximations in many-server queues with abandonment." Queueing Systems 78.1 (2014): 1-29.

\bibitem[{Dai and He (2013)}]{DH13}
Dai, J. G., S. He. ``Many-server queues with customer abandonment: Numerical analysis of their diffusion model." Stochastic Systems 3.1 (2013): 96-146.

\bibitem[{Dai et al. (2020)}]{Dai20}
Dai, J. G., J.M. Harrison. Processing Networks: Fluid Models and Stability. Cambridge University Press, 2020.

\bibitem[{Daley (1977)}]{Daley.77}
Daley, D. ``Inequalities for moments of tails of random variables, with a queueing application." Probability Theory and Related Fields 41.2 (1977): 139-143.

\bibitem[{Daley (1978)}]{Daley.78}
Daley, D. "Bounds for the variance of certain stationary point processes." Stochastic Processes and their Applications 7, no. 3 (1978): 255-264.

\bibitem[{Daley (1980)}]{Daley.80}
Daley, D. "Tight bounds for the renewal function of a random walk." The Annals of Probability 8, no. 3 (1980): 615-621.

\bibitem[{Daley and Rolski (1984)}]{Daley.84}
Daley, D., T. Rolski. ``Some Comparability Results for Waiting Times in Single- and Many-Server Queues.” Journal of Applied Probability, vol. 21, no. 4, 1984, pp. 887–900.

\bibitem[{Daley and Rolski (1992)}]{Daley.92}
Daley, D., T. Rolski. ``Light traffic approximations in many-server queues." Advances in applied probability 24, no. 1 (1992): 202-218.

\bibitem[{Daley (1997)}]{Daley.97}
Daley, D.: Some results for the mean waiting-time and workload in GI/GI/k queues. In: Dshalalow,
J.H. (ed.) Frontiers in Queueing: Models and Applications in Science and Engineering, Boca
Raton, FL, USA, pp. 35–59 (1997).

\bibitem[{Davis et al. (1995)}]{Davis.95}
Davis, J.L., W. Massey, W. Whitt. ``Sensitivity to the service-time distribution in the nonstationary Erlang loss model." Management Science 41, no. 6 (1995): 1107-1116.

\bibitem[{Doig (1957)}]{Doig57}
Doig, A. ``A bibliography on the theory of queues." Biometrika (1957): 490-514.

\bibitem[{Downey (1991)}]{Downey.91}
Downey, P. "Bounding Synchronization Overhead for Parallel Iteration." ORSA Journal on Computing 3, no. 4 (1991): 288-298.

\bibitem[{Eick et al. (1993)}]{Eick93}
Eick, S., W. Massey, W. Whitt. ``The physics of the Mt/G/$\infty$ queue." Operations Research 41, no. 4 (1993): 731-742.

\bibitem[{Figiel et al. (1997)}]{F97}
Figiel, T., P. Hitczenko, W. Johnson, G. Schechtman, and J. Zinn.  ``Extremal properties of Rademacher functions with applications to the Khintchine and Rosenthal inequalities." Transactions of the American Mathematical Society 349.3 (1997): 997-1027.

\bibitem[{Franken et al. (1982)}]{Franken82}
Franken, P., D. Konig, U. Arndt, V. Schmidt. ``Queues and Point Processes." John Wiley and Sons, Inc., 1 Wiley Drive, Somerset, N.J. 08873, 1982, 230 (1982).

\bibitem[{Gamarnik and Goldberg (2013)}]{GG13}
Gamarnik, D., D.A. Goldberg. ``Steady-state $ GI/G/n $ queue in the Halfin–Whitt regime." The Annals of Applied Probability 23.6 (2013): 2382-2419.

\bibitem[{Gamarnik and Momcilovic (2008)}]{GM08}
Gamarnik, D., P. Momcilovic. ``Steady-state analysis of a multiserver queue in the Halfin-Whitt regime." Advances in Applied Probability 40.2 (2008): 548-577.

\bibitem[{Gamarnik and Stolyar (2012)}]{GS12}
Gamarnik, D., A. Stolyar. ``Multiclass multiserver queueing system in the Halfin–Whitt heavy traffic regime: asymptotics of the stationary distribution." Queueing Systems 71.1-2 (2012): 25-51.

\bibitem[{Gamarnik and Zeevi (2006)}]{GZ06}
Gamarnik, D., A. Zeevi. ``Validity of heavy traffic steady-state approximations in generalized Jackson networks." The Annals of Applied Probability (2006): 56-90.

\bibitem[{Gaunt and Walton (2020)}]{Gaunt20}
Gaunt, R., N. Walton. ``Stein’s method for the single server queue in heavy traffic." Statistics and Probability Letters 156 (2020): 108566.

\bibitem[{Glynn (1987)}]{Glynn87}
Glynn, P. ``Upper bounds on Poisson tail probabilities." Operations research letters 6, no. 1 (1987): 9-14.

\bibitem[{Goldberg (2016)}]{G16}
Goldberg, D.A. ``On the steady-state probability of delay and large negative deviations for the $ GI/GI/n $ queue in the Halfin-Whitt regime." Under revision, previous version available at arXiv preprint arXiv:1307.0241 (2016).  \URL{https://arxiv.org/abs/1307.0241v2}  

\bibitem[{Goldberg (2016)}]{G16b}
Goldberg, D.A., D. Katz-Rogozhnikov, Y. Lu, M. Sharma, and M. Squillante. "Asymptotic optimality of constant-order policies for lost sales inventory models with large lead times." Mathematics of Operations Research 41, no. 3 (2016): 898-913.

\bibitem[{Goldberg (2017)}]{G17}
Goldberg, D.A., Y. Li. ``Heavy-tailed queues in the Halfin-Whitt regime." Under revision, previous version available at arXiv preprint arXiv:1707.07775 (2017).  \URL{https://arxiv.org/abs/1707.07775}

\bibitem[{Goldberg (2017)}]{G17c}
Goldberg, D.A., Y. Li. "Simple and explicit bounds for multi-server queues with universal 1/(1-rho) scaling." arXiv preprint arXiv:1706.04628 (2017).

\bibitem[{Goldberg (2018)}]{G18}
Goldberg, D.A., D. Mukherjee, Y. Li. ``Large deviations analysis for the $ M/H_2/n+ M $ queue in the Halfin-Whitt regime." Under revision, previous version available at arXiv preprint arXiv:1803.01082 (2018).  \URL{https://arxiv.org/abs/1803.01082}

\bibitem[{Gong et al. (1992)}]{Gong92}
Gong, W., J. Hu. "The MacLaurin series for the GI/G/1 queue." Journal of Applied Probability 29, no. 1 (1992): 176-184.

\bibitem[{Grosof et al. (2018)}]{grosof2018srpt}
Grosof, I., Z. Scully, M. Harchol-Balter. ``SRPT for multiserver systems." Performance Evaluation 127 (2018): 154-175.

\bibitem[{Grosof et al. (2021)}]{grosof2021finite}
Grosof, I., M. Harchol-Balter, A. Scheller-Wolf. ``The Finite-Skip Method for Multiserver Analysis." arXiv preprint arXiv:2109.12663 (2021).

\bibitem[{Gross et al. (2011)}]{Gross.11}
Gross, D., J. Shortle, J. Thompson, C. Harris. Fundamentals of Queueing Theory. Vol. 627. John Wiley and Sons, 2011.

\bibitem[{Grubel et al. (1981)}]{Grubel.81}
Grubel, R., and U. Jensen. "On the moments of the number of renewal epochs." ZAMM? Zeitschrift f r Angewandte Mathematik und Mechanik 61, no. 10 (1981): 531-532.

\bibitem[{Gupta et al. (2010)}]{Gupta10}
Gupta, V., M. Harchol-Balter, J. G. Dai, B. Zwart.  ``On the inapproximability of M/G/K: why two moments of job size distribution are not enough." Queueing Systems 64.1 (2010): 5-48.

\bibitem[{Gupta et al. (2011)}]{Gupta11}
Gupta, V., T. Osogami. ``Tight moments-based bounds for queueing systems." In Proceedings of the ACM SIGMETRICS joint international conference on Measurement and modeling of computer systems, pp. 133-134. 2011.

\bibitem[{Gurvich et al. (2013)}]{Gur13}
Gurvich, I., J. Huang, A. Mandelbaum. ``Excursion-based universal approximations for the Erlang-A queue in steady-state." Mathematics of Operations Research 39.2 (2013): 325-373.

\bibitem[{Gurvich et al. (2014)}]{Gur14}
Gurvich, I. ``Diffusion models and steady-state approximations for exponentially ergodic Markovian queues." The Annals of Applied Probability 24.6 (2014): 2527-2559.

\bibitem[{Gut (2009)}]{Gut09}
Gut, A. Stopped random walks. Springer-Verlag New York Incorporated, 2009.

\bibitem[{Halfin and Whitt (1981)}]{HW.81}
Halfin, S., W. Whitt. ``Heavy-traffic limits for queues with many exponential servers." Operations research 29.3 (1981): 567-588.

\bibitem[{Harchol-Balter et al. (2009)}]{HB.09}
Harchol-Balter, M., A. Scheller-Wolf, A. Young. ``Surprising results on task assignment in server farms with high-variability workloads." ACM SIGMETRICS Performance Evaluation Review 37.1 (2009): 287-298.

\bibitem[{Hariel (1988)}]{Hariel.88}
Harel, A. ``Sharp bounds and simple approximations for the Erlang delay and loss formulas." Management Science 34, no. 8 (1988): 959-972.

\bibitem[{Hariel (2010)}]{Hariel.10}
Harel, A. ``Sharp and simple bounds for the Erlang delay and loss formulae." Queueing Systems 64, no. 2 (2010): 119-143.

\bibitem[{Heyman et al. (1982)}]{Heyman82}
D. P. Heyman and M. J. Sobel, Stochastic Models in Operations Research, Vol. I, New York: McGraw-Hili, 1982.

\bibitem[{Hitczenko (1990)}]{hitczenko1990best}
Hitczenko, P. ``Best constants in martingale version of Rosenthal's inequality." The Annals of Probability (1990): 1656-1668.

\bibitem[{Hokstad (1985)}]{hokstad85}
Hokstad, P. ``Relations for the workload of the GI/G/s queue." Advances in applied probability 17, no. 4 (1985): 887-904.

\bibitem[{Hong et al. (2021)}]{hong2021sharp}
Hong, Y., W. Wang. ``Sharp Waiting-Time Bounds for Multiserver Jobs." arXiv preprint arXiv:2109.05343 (2021).
 
\bibitem[{Huang et al. (2016)}]{huang16}
Huang, J., I. Gurvich. ``Beyond heavy-traffic regimes: Universal bounds and controls for the single-server queue." Operations Research 66, no. 4 (2018): 1168-1188.

\bibitem[{Hunter (1969)}]{hunter69}
Hunter, J. "On the moments of Markov renewal processes." Advances in Applied Probability 1, no. 2 (1969): 188-210.

\bibitem[{Iglehart and Whitt (1970)}]{IW.70b}
Iglehart, D., W. Whitt. ``Multiple channel queues in heavy traffic. II: Sequences, networks, and batches." Advances in Applied Probability 2.02 (1970): 355-369.

\bibitem[{Janssen et al. (2008)}]{Janssen08}
Janssen, A., J. S. H. Van Leeuwaarden. ``Back to the roots of the M/D/s queue and the works of Erlang, Crommelin and Pollaczek." Statistica Neerlandica 62.3 (2008): 299-313.

\bibitem[{Janssen et al. (2008)}]{Janssen08b}
Janssen, A., J. S. H. Van Leeuwaarden, B. Zwart. ``Corrected asymptotics for a multi-server queue in the Halfin-Whitt regime." Queueing Systems 58.4 (2008): 261.

\bibitem[{Janssen et al. (2008)}]{Janssen08c}
Janssen, A., J. S. H. Van Leeuwaarden, B. Zwart. ``Gaussian expansions and bounds for the Poisson distribution applied to the Erlang B formula." Advances in Applied Probability 40, no. 1 (2008): 122-143.

\bibitem[{Janssen et al. (2008)}]{Janssen11}
Janssen, A., J.S.H. Van Leeuwaarden, B. Zwart. ``Refining square-root safety staffing by expanding Erlang C." Operations Research 59.6 (2011): 1512-1522.

\bibitem[{Jin et al. (2021)}]{Jin21}
Jin, X., G. Pang, L. Xu, X. Xu. ``An approximation to steady-state of M/Ph/n+ M queue." arXiv preprint arXiv:2109.03623 (2021).

\bibitem[{Johnson et al. (1985)}]{J.85}
Johnson, W., G. Schechtman, J. Zinn. ``Best constants in moment inequalities for linear combinations of independent and exchangeable random variables." The Annals of Probability (1985): 234-253.

\bibitem[{Kennedy (1972)}]{Kennedy72}
Kennedy, D. ``Rates of convergence for queues in heavy traffic. II: Sequences of queueing systems." Advances in Applied Probability 4.02 (1972): 382-391.

\bibitem[{Kiefer and Wolfowitz (1956)}]{Kiefer56}
Kiefer, J., and J. Wolfowitz. "On the characteristics of the general queueing process, with applications to random walk." The Annals of Mathematical Statistics (1956): 147-161.

\bibitem[{Kingman (1962)}]{Kingman62a}
Kingman, J. F. C. ``On queues in heavy traffic." Journal of the Royal Statistical Society. Series B (Methodological) (1962): 383-392.

\bibitem[{Kingman (1964)}]{K64}
Kingman, J.F.C. "A martingale inequality in the theory of queues." In Mathematical Proceedings of the Cambridge Philosophical Society, vol. 60, no. 2, pp. 359-361. Cambridge University Press, 1964.

\bibitem[{Kingman (1966)}]{K65}
Kingman, J. F. C. ``The heavy traffic approximation in the theory of queues." Proceedings of the Symposium on Congestion Theory. No. 2. University of North Carolina Press, Chapel Hill, NC, 1965.

\bibitem[{Kingman (1970)}]{K70}
Kingman, J. F. C. ``Inequalities in the Theory of Queues.” Journal of the Royal Statistical Society. Series B (Methodological), vol. 32, no. 1, 1970, pp. 102–110., www.jstor.org/stable/2984406.

\bibitem[{Kollerstrom (1974)}]{Kollerstrom74}
Kollerstrom, J. ``Heavy Traffic Theory for Queues with Several Servers. I." Journal of Applied Probability (1974): 544-552.

\bibitem[{Kollerstrom (1979)}]{Kollerstrom79}
Kollerstrom, J. ``Heavy Traffic Theory for Queues with Several Servers. II." Journal of Applied Probability (1979): 393-401.

\bibitem[{Kollerstrom (1981)}]{Kollerstrom81}
Kollerstrom, J. "A second-order heavy traffic approximation for the queue GI/G/1." Advances in Applied Probability 13, no. 1 (1981): 167-185.

\bibitem[{Leadbetter (1963)}]{Leadbettter63}
Leadbetter, M. R. "On series expansions for the renewal moments." Biometrika 50, no. 1-2 (1963): 75-80.

\bibitem[{Lehman (1963)}]{Lehman63}
Lehman, E. "Shapes, moments and estimators of the Weibull distribution." IEEE Transactions on Reliability 12, no. 3 (1963): 32-38.

\bibitem[{Longnecker and Serfling (1977)}]{LS77}
Longnecker, M., R. J. Serfling. ``General moment and probability inequalities for the maximum partial sum." Acta Mathematica Hungarica 30.1-2 (1977): 129-133.

\bibitem[{Loulou (1973)}]{Loulou73}
Loulou, R. ``Multi-channel queues in heavy traffic." Journal of Applied Probability 10.04 (1973): 769-777.

\bibitem[{Maglaras et al. (2018)}]{Mag18}
Maglaras, C., J. Yao, A. Zeevi. ``Optimal price and delay differentiation in large-scale queueing systems." Management science 64, no. 5 (2018): 2427-2444.

\bibitem[{Makino (1969)}]{Makino.69}
Makino, T. ``Investigation of the mean waiting time for queueing system with many servers." Annals of the Institute of Statistical Mathematics 21.1 (1969): 357-366.

\bibitem[{Mandelbaum et al. (1998)}]{MMR.98}
Mandelbaum, A., W. Massey, M. Reiman. ``Strong approximations for Markovian service networks." Queueing Systems 30.1 (1998): 149-201.

\bibitem[{Mandelbaum et al. (2012)}]{Mandelbaum2012}
Mandelbaum, A., P. Momcilovic. ``Queues with many servers and impatient customers." Mathematics of Operations Research 37, no. 1 (2012): 41-65.

\bibitem[{Mori (1975)}]{Mori.75}
Mori, M.  ``Some bounds for queues." J. Operat. Res. Soc. Japan 18 (1975): 152-181.

\bibitem[{Nadarajah et al. (2022)}]{Nad.22}
Nadarajah, Saralees, and Idika E. Okorie. "On the tail integral formulae for real-valued random variables." The Mathematical Gazette 106, no. 567 (2022): 487-493.

\bibitem[{Nagaev (1970a)}]{Nagaev.70}
Nagaev, S. ``On the speed of convergence in a boundary problem. I." Theory of Probability and Its Applications 15.2 (1970): 163-186.

\bibitem[{Nagaev (1970b)}]{Nagaev.70b}
Nagaev, S. ``On the speed of convergence in a boundary problem. II." Theory of Probability and Its Applications 15.3 (1970): 403-429.

\bibitem[{Oliver (1974)}]{Oliver.74}
Oliver, S.Y. ``Stochastic bounds for heterogeneous-server queues with Erlang service times." Journal of Applied Probability 11.04 (1974): 785-796.

\bibitem[{Olvera-Cravioto et al. (2011)}]{Olvera.11}
Olvera-Cravioto, M., P. Glynn. "Uniform approximations for the M/G/1 queue with subexponential processing times." Queueing Systems 68, no. 1 (2011): 1-50.

\bibitem[{Olvera-Cravioto et al. (2011b)}]{Olvera.11b}
Olvera-Cravioto, M., J. Blanchet, P. Glynn. "On the transition from heavy traffic to heavy tails for the M/G/1 queue: the regularly varying case." (2011): 645-668.

\bibitem[{Ovuworie (1980)}]{Ovuworie80}
Ovuworie, G. ``Multi-channel queues: a survey and bibliography." International Statistical Review/Revue Internationale de Statistique (1980): 49-71.

\bibitem[{Reed (2009)}]{R.09}
Reed, J. ``The G/GI/N queue in the Halfin–Whitt regime." The Annals of Applied Probability 19.6 (2009): 2211-2269.

\bibitem[{Reiman (1984)}]{Reiman84}
Reiman, M. ``Open queueing networks in heavy traffic." Mathematics of operations research 9, no. 3 (1984): 441-458.

\bibitem[{Rein and Scheller-Wolf (2013)}]{RSW13}
Vesilo, R., A. Scheller-Wolf. ``Delay Moment Bounds for Multiserver Queues with Infinite Variance Service Times." INFOR: Information Systems and Operational Research 51.4 (2013): 161-174.

\bibitem[{Ren and Liang (2001)}]{Ren.01}
Ren, Y., H. Liang. ``On the best constant in Marcinkiewicz–Zygmund inequality." Statistics and probability letters 53.3 (2001): 227-233.

\bibitem[{Rolski and Stoyan (1976)}]{Rol.76}
Rolski, T., D. Stoyan, ``On the comparison of waiting times in GI/G/1 queues". Oper. Res. 24, 197-200 (1976).

\bibitem[{Sadowsky (1991)}]{Sadowsky91}
Sadowsky, J. ``Large deviations theory and efficient simulation of excessive backlogs in a GI/GI/m queue." IEEE Transactions on Automatic Control 36.12 (1991): 1383-1394.

\bibitem[{Scheller-Wolf (2003)}]{Scheller03}
Scheller-Wolf, A. ``Necessary and sufficient conditions for delay moments in FIFO multiserver queues with an application comparing s slow servers with one fast one." Operations Research 51.5 (2003): 748-758.

\bibitem[{Scheller-Wolf et al. (2006)}]{Scheller06}
Scheller-Wolf, A., R. Vesilo.  ``Structural interpretation and derivation of necessary and sufficient conditions
for delay moments in FIFO multiserver queues."  Queueing Syst. 54(3), 221–232 (2006).

\bibitem[{Scully et al. (2020)}]{scully2020gittins}
Scully, Z., I. Grosof, M. Harchol-Balter. ``The Gittins policy is nearly optimal in the M/G/k under extremely general conditions." Proceedings of the ACM on Measurement and Analysis of Computing Systems 4, no. 3 (2020): 1-29.

\bibitem[{Seshadri (1996)}]{Seshadri.96}
Seshadri, S.  ``A sample path analysis of the delay in the M/G/C system."  J. Appl. Prob 33 (1996): 256-266.

\bibitem[{Sevastyanov (1957)}]{Sev57}
Sevastyanov, B.A. ``An ergodic theorem for Markov processes and its application to telephone systems with refusals." Theory of Probability and Its Applications 2, no. 1 (1957): 104-112.

\bibitem[{Shanthikumar et al. (1989)}]{Shanthikumar89}
Shanthikumar, J.G., D. Yao. ``Stochastic monotonicity in general queueing networks." Journal of Applied Probability 26, no. 2 (1989): 413-417.

\bibitem[{Smith (1959)}]{Smith59}
Smith, W. "On the cumulants of renewal processes." Biometrika 46, no. 1/2 (1959): 1-29.

\bibitem[{Smith and Whitt (1981)}]{Smith81}
Smith, D., W. Whitt. ``Resource sharing for efficiency in traffic systems." Bell System Technical Journal 60.1 (1981): 39-55.

\bibitem[{Suzuki et al. (1970)}]{Suzuki70}
Suzuki, T., Y. Yoshida.  ``Inequalities for many-server queue and other queues." J. Oper. Res. Soc. Japan 13 (1970): 59-77.

\bibitem[{Szczotka (1999)}]{Sz.99}
Szczotka, W. ``Tightness of the stationary waiting time in heavy traffic." Advances in Applied Probability (1999): 788-794.

\bibitem[{Szczotka and Woyczynski (2004)}]{Sz.04}
Szczotka, W., W. A. Woyczynski. ``Heavy-tailed dependent queues in heavy traffic." Probability and Mathematical Statistics - Wroclaw University.  24.1 (2004): 67.

\bibitem[{Taga (1963)}]{Taga63}
Taga, Y. "On high order moments of the number of renewals." Annals of the Institute of Statistical Mathematics 15 (1963): 187-196.

\bibitem[{Takacs (1962)}]{Tak62}
Takacs, L. "A single-server queue with Poisson input." Operations research 10, no. 3 (1962): 388-394.

\bibitem[{Vera et al. (2021)}]{vera2021}
Vera, A., and S. Banerjee. ``The bayesian prophet: A low-regret framework for online decision making." Management Science 67, no. 3 (2021): 1368-1391.

\bibitem[{Wang et al. (2021)}]{wang2021zero}
Wang, W., Q. Xie, M. Harchol-Balter. ``Zero queueing for multi-server jobs." In Abstract Proceedings of the 2021 ACM SIGMETRICS/International Conference on Measurement and Modeling of Computer Systems, pp. 13-14. 2021.

\bibitem[{Whitt (1974)}]{W74}
Whitt, W. ``The continuity of queues." Advances in Applied Probability 6, no. 1 (1974): 175-183.

\bibitem[{Whitt (1980)}]{W80}
Whitt, W. ``The effect of variability in the GI/G/s queue".  J. Appl. Probab. 17, 1062-1071 (1980).

\bibitem[{Whitt (1981)}]{W81}
Whitt, W. ``Comparing counting processes and queues." Advances in Applied Probability 13, no. 1 (1981): 207-220.

\bibitem[{Whitt (1982a)}]{Whitt.82}
Whitt, W. ``Refining diffusion approximations for queues." Operations Research Letters 1.5 (1982): 165-169.

\bibitem[{Whitt (1982b)}]{Whitt.82b}
Whitt, W. ``The Marshall and Stoyan bounds for IMRL/G/1 queues are tight." Operations Research Letters 1.6 (1982): 209-213.

\bibitem[{Whitt (1982c)}]{Whitt.82d}
Whitt, W. ``On the heavy-traffic limit theorem for GI/G/$\infty$ queues." Advances in Applied Probability 14, no. 1 (1982): 171-190.

\bibitem[{Whitt (1983a)}]{Whitt.83}
Whitt, W. ``Comparison conjectures about the M/G/s queue." Operations Research Letters 2.5 (1983): 203-209.

\bibitem[{Whitt (1983b)}]{Whitt.83c}
Whitt, W. ``The queueing network analyzer." The bell system technical journal 62, no. 9 (1983): 2779-2815.

\bibitem[{Whitt (1984a)}]{Whitt.84}
Whitt, W. ``On approximations for queues, I: Extremal distributions." ATT Bell Laboratories Technical Journal 63.1 (1984): 115-138.

\bibitem[{Whitt et al. (1984b)}]{Whitt.84b}
Whitt, W., J.G. Kuncewicz. ``On approximations for queues, II: Shape constraints." ATT Bell Laboratories Technical Journal 63.1 (1984): 139-161.

\bibitem[{Whitt (1984)}]{Whitt.84c}
Whitt, W. ``On approximations for queues, III: Mixtures of exponential distributions." ATT Bell Laboratories Technical Journal 63.1 (1984): 163-175.

\bibitem[{Whitt (1993)}]{Whitt.93}
Whitt, W. ``Approximations for the GI/G/m queue." Production and Operations Management 2.2 (1993): 114-161.

\bibitem[{Whitt (2000)}]{Whitt.00}
Whitt, W. ``The impact of a heavy-tailed service-time distribution upon the M/GI/s waiting-time distribution." Queueing Systems 36, no. 1 (2000): 71-87.

\bibitem[{Wolff (1987)}]{Wolff87}
Wolff, R.W. ``Upper bounds on work in system for multichannel queues." Journal of applied probability 24.02 (1987): 547-551.

\bibitem[{Wolff (2011)}]{Wolff2011}
Wolff, R.W. ``Little’s law and related results." Wiley encyclopedia of operations research and management science 4 (2011): 2828-2841.

\bibitem[{Worthington (2009)}]{Worth09}
Worthington, D. ``Reflections on queue modelling from the last 50 years." Journal of the Operational Research Society 60.1 (2009): S83-S92.

\bibitem[{Xin et al. (2016)}]{Xin16}
Xin, L., D.A. Goldberg. ``Optimality gap of constant-order policies decays exponentially in the lead time for lost sales models." Operations Research 64, no. 6 (2016): 1556-1565.

\end{thebibliography}
\end{document}